\documentclass[12pt,reqno]{amsart}
\usepackage{cooltooltips}
\usepackage[all,arc]{xy}
\usepackage{enumerate}
\usepackage{mathrsfs}
\usepackage{euscript}
\usepackage{hyperref}
\usepackage{amsthm}
\usepackage{amsmath}
\usepackage{comment} 
\usepackage{tikz-cd}
\usepackage{xspace}

\tikzset{
labl1/.style={anchor=north, rotate=90, inner sep=1.2mm}
}

\usepackage[belowskip=-15pt,aboveskip=0pt]{caption}

 \usepackage{amsmath}
   \usepackage{amssymb}
   \usepackage{bm}

\usepackage{esvect}

\usepackage[latin1]{inputenc}

\usepackage{textcomp}

\usepackage{usnomencl}

\usepackage{graphicx}
\usepackage{color}
\usepackage{eso-pic}

\usepackage{mathtools}
\usepackage{multicol}
\usepackage{paralist}
\usepackage{geometry}
\usepackage{algorithmic}
\usepackage[ruled,vlined]{algorithm2e}
\usepackage{comment}
\usepackage{tikz-cd}
\usepackage{url}
\usepackage{graphicx,scalerel}
\usetikzlibrary{decorations.pathmorphing}
\usepackage{float}
\usepackage{tabulary}
\usepackage{booktabs}
\usepackage{mathabx}

\usepackage[sort,numbers]{natbib}
\usepackage{filecontents}

\newtheorem{theorem}{Theorem}[section]

\newtheorem{lemma}[theorem]{Lemma}
\newtheorem{definition}[theorem]{Definition}
\newtheorem{remark}[theorem]{Remark}

\newcommand{\Ops}[1]{\mathbf{Ops}_{\substack{\text{\scalebox{0.8}{$#1$}}}}}
\newcommand{\Psh}[2]{\mathbf{Psh}_{\scalebox{0.5}{$\mathbf{#1},\! \mathbf{#2}$}}}
\newcommand{\Sh}[2]{\mathbf{Sh}_{\scalebox{0.5}{$\mathbf{#1},\! \mathbf{#2}$}}}
\newcommand{\iuv}[2]{\mathfrak{i}_{\scalebox{0.8}{$\scriptstyle #1$}\mkern-0.7mu\scalebox{0.8}{$\scriptstyle , $}\mkern-1mu\scalebox{0.8}{$\scriptstyle #2$}}}
\newcommand{\iuvop}[3]{\mathfrak{i}_{\scalebox{0.8}{$\scriptstyle #1$}\mkern-0.7mu\scalebox{0.8}{$\scriptstyle , $}\mkern-1mu\scalebox{0.8}{$\scriptstyle #2$}}^{\operatorname{#3}}}
\newcommand{\rt}{\mathbf{RTS}}
\newcommand{\lrt}{\mathbf{LRTS}}
\newcommand{\topo}{\mathbf{(o)Top}}
\newcommand{\Uij}[2]{U_{\substack{\text{\scalebox{0.7}{$\mkern-1mu #1 \! \wedge\! #2$}}}}}

\newcommand{\rgdt}[1]{\mathfrak{#1}_{\substack{\text{\scalebox{0.7}{$\mkern-1mu \mathbf{Top}$}}}}}
\newcommand{\rgds}[1]{\mathfrak{#1}_{\substack{\text{\scalebox{0.7}{$\mkern-1mu \mathbf{Sh}$}}}}}
\newcommand{\subsm}[2]{{#1}_{\substack{\text{\scalebox{0.8}{$\mkern-0.8mu #2$}}}}}
\newcommand{\subsms}[2]{{#1}_{{}_{\substack{\text{\scalebox{0.8}{$\mkern-0.8mu #2$}}}}}}
\newcommand{\subsmt}[3]{{#1}_{\substack{\text{\scalebox{0.8}{$\mkern-1mu #2 \mkern-1mu , \! #3$}}}}}
\newcommand{\gdt}[1]{\mathbf{#1}_{\substack{\text{\scalebox{0.7}{$\mkern-1mu \mathbf{Top}$}}}}}
\newcommand{\gds}[1]{\mathbf{#1}_{\substack{\text{\scalebox{0.7}{$\mkern-1mu \mathbf{Sh}$}}}}}
\newcommand{\dt}[1]{{#1}_{\substack{\text{\scalebox{0.7}{$\mkern-1mu \mathbf{Top}$}}}}}
\newcommand{\ds}[1]{{#1}_{\substack{\text{\scalebox{0.7}{$\mkern-1mu \mathbf{Sh}$}}}}}

\newcommand{\uiuv}[3]{\underbar{$\mathfrak{i}$}_{{\substack{\text{\scalebox{0.8}{$#1$}}}}_{\scalebox{0.6}{$\scriptstyle #2$}\mkern-0.7mu\scalebox{0.6}{$\scriptstyle , $}\mkern-1mu\scalebox{0.6}{$\scriptstyle #3$}}}}

\newcommand{\op}{{\substack{\text{\scalebox{0.7}{$\operatorname{op}$}}}}}
\newcommand{\sch}{\mathbf{Sch}}
\newcommand{\subop}{{\subseteq}_{\operatorname{op}}}
\newcommand{\Rel}[1]{\mathcal{R}_{\substack{\text{\scalebox{0.8}{$#1$}}}}}

\newcommand{\Co}[1]{{#1}_{\substack{\text{\scalebox{0.7}{$\mkern-1.5mu 0$}}}}}
\newcommand{\Cn}[1]{{#1}_{\substack{\text{\scalebox{0.7}{$ \mkern-1.5mu 1$}}}}}

\newcommand{\Cos}[1]{{#1}_{{}_{\substack{\text{\scalebox{0.7}{$\mkern-1.5mu 0$}}}}}}
\newcommand{\Cns}[1]{{#1}_{{}_{\substack{\text{\scalebox{0.7}{$ \mkern-1.5mu 1$}}}}}}
\newcommand{\mbf}[1]{\mathbf{#1}}
\newcommand{\etaij}[2]{\eta_{\substack{\text{\scalebox{0.7}{${#1,\![\mkern-1mu #1\mkern-1mu ,\mkern-1mu #2\mkern-1mu]}$}}}}}
\newcommand{\etaji}[2]{\eta_{\substack{\text{\scalebox{0.7}{${#2,\![\mkern-1mu #2\mkern-1mu ,\mkern-1mu #1\mkern-1mu ]}$}}}}}
\newcommand{\tauij}[2]{\tau_{\substack{\text{\scalebox{0.7}{$\mkern-1.2mu [\mkern-1mu #1\mkern-1mu,\mkern-1mu #2\mkern-1mu]$}}}}}
\newcommand{\tauji}[2]{\tau_{\substack{\text{\scalebox{0.7}{$\mkern-1.2mu [\mkern-1mu #2 \mkern-1mu,\mkern-1mu #1\mkern-1mu]$}}}}}
\newcommand{\etaijk}[4]{\eta_{\substack{\text{\scalebox{0.7}{${[\mkern-1mu #2 \mkern-1mu,\mkern-1mu #1\mkern-1mu ],\![\mkern-1mu #2 \mkern-1mu ,\mkern-1mu #3 \mkern-1mu,\mkern-1mu #4\mkern-1mu]}$}}}}}
\newcommand{\etaiijk}[3]{\eta_{\substack{\text{\scalebox{0.7}{${\mkern-1mu #1 \mkern-1mu,\![\mkern-1mu #1 \mkern-1mu ,\mkern-1mu #2 \mkern-1mu,\mkern-1mu #3\mkern-1mu]}$}}}}}
\newcommand{\tauijk}[4]{\tau_{{\substack{\text{\scalebox{0.7}{\!#1}}}}_{\substack{\text{\scalebox{0.8}{$[\mkern-1mu #2 \mkern-1mu,\mkern-1mu #3\mkern-1mu]$}}}}}}
\newcommand{\limi}[1]{\operatorname{lim}\!#1}
\newcommand{\phijk}[3]{\varphi_{{\substack{\text{\scalebox{0.7}{$\! #1$}}}}_{\substack{\text{\scalebox{0.8}{$(\mkern-1mu #2 \mkern-1mu,\mkern-1mu #3\mkern-1mu)$}}}}}}

\newcommand{\Goi}[2]{\Co{#1} \mkern-1.7mu (\mkern-1mu #2 \mkern-1mu)}
\newcommand{\Goij}[3]{\Co{#1} \mkern-1.7mu (\mkern-1mu #2 \mkern-1mu ,\mkern-1mu #3 \mkern-1mu )}
\newcommand{\Goijk}[4]{\Co{#1} \mkern-1.7mu (\mkern-1mu #2 \mkern-1mu ,\mkern-1mu #3 \mkern-1mu ,\mkern-1mu #4 \mkern-1mu )}

\newcommand{\Coi}[2]{#1 \mkern-1.7mu (\mkern-1mu #2 \mkern-1mu)}
\newcommand{\Coij}[3]{#1\mkern-1.7mu (\mkern-1mu #2 \mkern-1mu ,\mkern-1mu #3 \mkern-1mu )}
\newcommand{\Coijk}[4]{#1 \mkern-1.7mu (\mkern-1mu #2 \mkern-1mu ,\mkern-1mu #3 \mkern-1mu ,\mkern-1mu #4 \mkern-1mu )}
\newcommand{\glI}[1]{\mathbf{GLD}(\mathrm{#1})}
\newcommand{\Gnij}[2]{\Cn{#1} \mkern-1.7mu (\mkern-1.2mu #2  \mkern-1.2mu )}
\newcommand{\Cnij}[2]{#1 \! ( \mkern-1.2mu #2  \mkern-1.2mu )}
\newcommand{\subsc}[2]{{}_{\substack{\text{\scalebox{0.6}{$\mkern-1.8mu #1$}}}} \mkern-1.8mu #2}
\newcommand{\randi}[2]{#1_{\substack{\text{\scalebox{0.7}{$\mkern-1mu \mkern-1.2mu #2\mkern-1.2mu$}}}}}
\newcommand{\randij}[3]{#1_{\substack{\text{\scalebox{0.7}{$\mkern-1mu [\mkern-1.2mu #2\mkern-1.2mu, \mkern-1.2mu #3 \mkern-1.2mu]$}}}}}
\newcommand{\randijk}[4]{#1_{\substack{\text{\scalebox{0.7}{$\mkern-1mu [\mkern-1.2mu #2\mkern-1.2mu, \mkern-1.2mu #3 \mkern-1.2mu,\mkern-1.2mu #4 \mkern-1.2mu]$}}}}}
\newcommand{\randji}[3]{#1_{\substack{\text{\scalebox{0.7}{$\mkern-1mu [\mkern-1.2mu #2\mkern-1.2mu, \mkern-1.2mu #3 \mkern-1.2mu]$}}}}}
\newcommand{\idij}[2]{{\operatorname{id}}_{\substack{\text{\scalebox{0.7}{$\mkern-1mu [\mkern-1.2mu #1\mkern-1.2mu, \mkern-1.2mu #2 \mkern-1.2mu]$}}}}}
\newcommand{\idijk}[3]{{\operatorname{id}}_{\substack{\text{\scalebox{0.7}{$\mkern-1mu [\mkern-1.2mu #1\mkern-1.2mu, \mkern-1.2mu #2 \mkern-1.2mu, \mkern-1.2mu #3 \mkern-1.2mu]$}}}}}
\newcommand{\dindi}[3]{{#1}_{{\substack{\text{\scalebox{0.8}{$\mkern-1mu #2$}}}}_{{\mkern-1mu}_{\substack{\text{\scalebox{0.7}{$\mkern-1mu #3$}}}}}}}
\newcommand{\dindij}[4]{{#1}_{{\substack{\text{\scalebox{0.8}{$\mkern-1mu #2$}}}}_{{\mkern-1mu}_{\substack{\text{\scalebox{0.7}{$\mkern-1mu [\mkern-1.2mu #3\mkern-1.2mu, \mkern-1.2mu #4 \mkern-1.2mu]$}}}}}}}

\newcommand{\dindijk}[5]{{#1}_{{\substack{\text{\scalebox{0.8}{$\mkern-1mu #2$}}}}_{{\mkern-1mu}_{\substack{\text{\scalebox{0.7}{$\mkern-1mu [\mkern-1.2mu #3\mkern-1.2mu, \mkern-1.2mu #4 \mkern-1.2mu, \mkern-1.2mu #5 \mkern-1.2mu]$}}}}}}}
\newcommand{\sta}{ \mkern-1mu {}_{\substack{\text{\scalebox{0.8}{$\mkern-1mu *$}}}}}

\newcommand{\dindiv}[4]{{#1}_{{\substack{\text{\scalebox{0.8}{$\mkern-1mu #2$}}}}_{{\mkern-1mu}_{\substack{\text{\scalebox{0.7}{$\mkern-1mu #3$}}}_{{}_{\substack{\text{\scalebox{0.6}{$\mkern-1mu #4$}}}}}}}}}

\newcommand{\dindgij}[5]{{#1}_{{\substack{\text{\scalebox{0.8}{$\mkern-1mu #2$}}}_{{\mkern-1mu}_{{\substack{\text{\scalebox{0.6}{$\mkern-1mu #3$}}}_{{}_{\substack{\text{\scalebox{0.6}{$\mkern-1mu [\mkern-1.2mu #4\mkern-1.2mu, \mkern-1.2mu #5 \mkern-1.2mu]$}}}}}}}}}}}

\newcommand{\dindgijk}[6]{{#1}_{{{\substack{\text{\scalebox{0.8}{$\mkern-1mu #2$}}}_{{\mkern-1mu}_{{\substack{\text{\scalebox{0.6}{$\mkern-1mu #3$}}}_{{}_{\substack{\text{\scalebox{0.6}{$\mkern-1mu [\mkern-1.2mu #4\mkern-1.2mu, \mkern-1.2mu #5 \mkern-1.2mu, \mkern-1.2mu #6 \mkern-1.2mu]$}}}}}}}}}}}}

\newcommand{\dindgiv}[5]{{#1}_{{\substack{\text{\scalebox{0.6}{$\mkern-1mu #2$}}}_{{{}_{\substack{\text{\scalebox{0.5}{$\mkern-1mu #3$}}}_{{{}_{\substack{\text{\scalebox{0.6}{$\mkern-1.2mu #4\mkern-1.2mu$}}}}}_{{}_{\substack{\text{\scalebox{0.5}{$\mkern-1mu #5$}}}}}}}}}}}}

\newcommand{\dindgijv}[6]{{#1}_{{\substack{\text{\scalebox{0.6}{$\mkern-1mu #2$}}}_{{{}_{\substack{\text{\scalebox{0.5}{$\mkern-1mu #3$}}}_{{{}_{\substack{\text{\scalebox{0.5}{$\mkern-1mu [\mkern-1.2mu #4\mkern-1.2mu, \mkern-1.2mu #5 \mkern-1.2mu]$}}}}}_{{}_{\substack{\text{\scalebox{0.5}{$\mkern-1mu #6$}}}}}}}}}}}}

\newcommand{\dindgijkv}[7]{{#1}_{{{\substack{\text{\scalebox{0.6}{$\mkern-1mu #2$}}}_{{\mkern-1mu}_{{}_{\substack{\text{\scalebox{0.5}{$\mkern-1mu #3$}}}_{{{}_{\substack{\text{\scalebox{0.5}{$\mkern-1mu [\mkern-1.2mu #4\mkern-1.2mu, \mkern-1.2mu #5 \mkern-1.2mu, \mkern-1.2mu #6 \mkern-1.2mu]$}}}}}_{{}_{\substack{\text{\scalebox{0.6}{$\mkern-1mu #7$}}}}}}}}}}}}}

\newcommand{\dindpi}[3]{{#1}_{{\substack{\text{\scalebox{0.8}{$\mkern-1mu #2$}}}}_{{{\mkern-1mu}_{\substack{\text{\scalebox{0.7}{$\mkern-1mu #3$}}}}}}'}}
\newcommand{\dindpij}[4]{{#1}_{{\substack{\text{\scalebox{0.8}{$\mkern-1mu #2$}}}}_{{\mkern-1mu}_{\substack{\text{\scalebox{0.7}{$\mkern-1mu [\mkern-1.2mu #3\mkern-1.2mu, \mkern-1.2mu #4 \mkern-1.2mu]$}}}}}'}}
\newcommand{\dindpijk}[5]{{\mkern-1mu #1}_{{\substack{\text{\scalebox{0.8}{$\mkern-1mu #2$}}}}_{{\mkern-1mu}_{\substack{\text{\scalebox{0.7}{$\mkern-1mu [\mkern-1.2mu #3\mkern-1.2mu, \mkern-1.2mu #4 \mkern-1.2mu, \mkern-1.2mu #5 \mkern-1.2mu]$}}}}}'}}

\newcommand{\dindsi}[3]{{#1}_{{\substack{\text{\scalebox{0.8}{$\mkern-1mu #2$}}}}_{{}_{\substack{\text{\scalebox{0.7}{$#3$}}}}}}}

\newcommand{\gcov}[1]{\bm{\mathcal C}_{\substack{\text{\scalebox{0.8}{$#1$}}}}}
\newcommand{\sho}[2]{{#1}_{\substack{\text{\scalebox{0.8}{$\mkern-1mu 0$}}}}\mkern-1mu (\mkern-1mu #2 \mkern-1mu)}

\newcommand{\circc}{\! \circ  \!}

\newcommand{\iotaij}[2]{\iota_{\substack{\text{\scalebox{0.7}{$\mkern-1.2mu [\mkern-1mu #1\mkern-1mu,\mkern-1mu #2\mkern-1mu]$}}}}}

\newcommand{\damas}[1]{{\color{blue}\textsf{$ \spadesuit \spadesuit \spadesuit$}  damas: [#1]}}
\newcommand\restr[2]{{
\left.\kern-\nulldelimiterspace 
#1 
\vphantom{\big|} 
\right|_{#2} 
}}

\newcommand{\suchthat}{\;\ifnum\currentgrouptype=16 \middle\fi|\;}



\title{\bfseries Gluing Data categories and gluing data functors}
\author{Sophie Marques and Damas Mgani}

\begin{document}

\begin{abstract}
We present a novel approach to the concept of gluing in mathematics by introducing the notions of a gluing data category and a  gluing data functor. Our work provides a formal categorical characterization of the notion of gluing in algebraic geometry. By using this characterization, we are able to describe gluing in a unified way that applies to a wide range of mathematical structures, including topological spaces, presheaves, sheaves, ringed topological spaces, locally ringed topological spaces, and schemes. Our results provide a fresh perspective on gluing that is both abstract and formal, offering a deeper understanding of this fundamental concept in mathematics.


\noindent \textbf{2020 Math. Subject Class:} 14A15, 18F20, 18F60.

\noindent \begin{center}

\rm e-mail: \href{mailto:d smarques@sun.ac.za}{ smarques@sun.ac.za}

\it
Department of Mathematical Sciences, 
University of Stellenbosch, \\
Stellenbosch, 7600, 
South Africa\\ 
\&
NITheCS (National Institute for Theoretical and Computational Sciences), \\
South Africa \\ \bigskip

\rm e-mail: \href{mailto:d.mgani99@gmail.com}{d.mgani99@gmail.com}

\it
Department of Mathematical Sciences, 
University of Stellenbosch, \\
Stellenbosch, 7600, 
South Africa\\ 
\end{center} 

\end{abstract}

\setcounter{tocdepth}{3}
\maketitle
  \tableofcontents

\newpage

\section*{Introduction} 

Mathematics, as a field, thrives on exploration and the relentless quest for unification. Novelty emerges when we dare to redefine and reinterpret well-established ideas from multiple vantage points. In our endeavor, we turn our focus to the act of "gluing," a fundamental technique in algebraic geometry. Gluing allows us to elegantly "patch together" smaller spaces to craft larger, more intricate ones. This technique forms the very fabric of algebraic geometry, and herein lies the novelty - not in introducing mere terminology but in offering a unified perspective on gluing.

At the heart of our work lies a fundamental question:

\begin{center}
{\bf What if we could define gluing uniformly across categories, including topological spaces, as a limit within the category, allowing us to determine in which categories the notion of gluing holds significance? }
\end{center}

The significance of our work lies in our ability to construct a "gluing data functor" tailored to any category.

This work resonates on multiple levels. Firstly, it empowers us to define gluing as a checkable property of a category, thus paving the way for deeper comprehension. One could identify which categories possess this defining property, and, significantly, it permits us to regard gluing as a limit within a certain category in each presented case, offering new insights into its essence.

In the ever-evolving landscape of mathematical research, the ability to unify, reinterpret, and offer fresh insights into established concepts stands as an essential pillar of progress. 

Every proof presented here is meticulously tailored and adapted to align with the present point of view and the overarching unifying context of this paper. The proof may resemble well-known ones, but we have included it here for the sake of completeness, providing all the essential details that are sometimes omitted in existing literature.

Our work is bridging algebraic geometry and formal category theory; for this reason, we offer an index of notation in the first section to help readers keep up with the uniformized notations. 

In the second section, we introduce the main contributions of this paper, including the definitions of a Gluing Data category (see Definition \ref{gic}), a gluing data functor, and a glued-up object (see Definition \ref{limext}). We then simplify the description of cones over a gluing data functor in two different ways (see Lemma \ref{cocone} and Theorem \ref{preglued}).

Section 3  systematically characterizes $\topo^{\op}$-glued-up objects (as demonstrated in Theorem \ref{gluingtop}). This section applies our framework to topological spaces.

In Section 4, we further expand our scope to encompass $\mathbf{E}\Psh{\Ops{X}}{Ab}$ and $\mathbf{E}\Sh{\Ops{X}}{Ab}$ categories. Here, we introduce the $\mathbf{E}\Psh{\Ops{X}}{Ab}$ (resp. $\mathbf{E}\Sh{\Ops{X}}{Ab}$)-gluing data functor, meticulously defined in Definition \ref{sheffnct}, and provide a comprehensive characterization of $\mathbf{E}\Psh{\Ops{X}}{Ab}$ (resp. $\mathbf{E}\Sh{\Ops{X}}{Ab}$)-glued-up objects, elucidated in Theorem \ref{rivegluedshf}.

Finally, in the last section, our exploration is based on topological spaces endowed with the (pre)sheaf structure but in particular we focus on defining $\rt$ (resp. $\lrt$, resp. $\sch$)-gluing data functor, meticulously outlined in Definition \ref{shcmct}. These definitions provide a gateway to characterizing $\rt$ (resp. $\lrt$, resp. $\sch$)-glued-up objects, an endeavor culminating in Theorem \ref{gluedsch}.

\section{Index of notation}
 
\begin{center} 
  \begin{tabular}{l p{11cm} }
  \textbf{CATEGORIES} & \\[0.5cm]
    $\mathbf{C}$ &   Category with the following components \\
    & $( \Co{\mbf{C}} ,\Cn{\mbf{C}}, \mathbf{d}_{\substack{\text{\scalebox{0.7}{$\mathbf{C}$}}}},\mathbf{c}_{\substack{\text{\scalebox{0.7}{$\mathbf{C}$}}}},\mathbf{e}_{\substack{\text{\scalebox{0.7}{$\mathbf{C}$}}}},\mathbf{m}_{\substack{\text{\scalebox{0.7}{$\mathbf{C}$}}}})$ where $ \Co{\mbf{C}}$ is a collection of objects in $\mathbf{C}$, $ \Cn{\mbf{C}}$ is a collection of morphisms in $\mathbf{C}$, $\mathbf{d}_{\substack{\text{\scalebox{0.7}{$\mathbf{C}$}}}}$ and $\mathbf{c}_{\substack{\text{\scalebox{0.7}{$\mathbf{C}$}}}}$ are maps from $ \Cn{\mbf{C}}$ to $ \Co{\mbf{C}}$, called domain and codomain respectively,  $\mathbf{e}_{\substack{\text{\scalebox{0.7}{$\mathbf{C}$}}}}$ is a map from $ \Co{\mbf{C}}$ to $ \Cn{\mbf{C}}$, called the identity, $\mathbf{m}_{\substack{\text{\scalebox{0.7}{$\mathbf{C}$}}}}$ is a map from the set
     $\{(f,g)\in  \Cn{\mbf{C}}\times  \Cn{\mbf{C}}\;|\;\mathbf{d}_{\substack{\text{\scalebox{0.7}{$\mathbf{C}$}}}}(f)=\mathbf{c}_{\substack{\text{\scalebox{0.7}{$\mathbf{C}$}}}}(g)\}$ to $ \Cn{\mbf{C}}$, called composition. Categories are defined in this way in \cite{maclane} and \cite{Joji};\\
    
    $\glI{I}$ &  Gluing Data category of type $\rm{I}$ (see Definition \ref{gic});\\ 
    
$\mathbf{Ab}$ &   Category of abelian groups where morphisms are \ morphisms of abelian groups; \\

  $\mathbf{Rings}$ & Category of rings without unity where morphisms are ring homomorphisms;\\
  
 $\mathbf{CURings}$ & Subcategory of $\mathbf{Rings}$ whose objects are commutative rings with unity and whose morphisms are unital ring homomorphisms;\\
 
$\Ops{X}$ & Category of open sets of a topological space $X$ where morphisms are inclusion maps;\\

$\Ops{X\!,\!x}$ & Full subcategory of $\Ops{X}$ whose object are open subsets of $X$ containing $x$ an element of $X$;\\

 $\topo$& Subcategory of topological spaces whose morphisms are (open)continuous maps;\\

$\mathbf{E}\Psh{\mathbf{C}}{\mathbf{D}}$ & Category of enriched presheaves on $\mathbf{C}$ with values in $\mathbf{D}$ (see Definition \ref{catshf}); \\

$\mathbf{E}\Sh{\mathbf{C}}{\mathbf{D}}$ &  Full subcategory of $\mathbf{E}\Psh{\mathbf{C}}{\mathbf{D}}$  whose objects are enriched sheaves on $\mathbf{C}$ with values in $\mathbf{D}$; \\

$\rt$ & Category of ringed topological spaces where morphisms are morphisms of ringed topological spaces;\\

$\lrt$ &  Subcategory of $\rt$ whose objects are locally ringed topological spaces where morphisms are morphisms of locally ringed topological spaces;\\

$\sch$ &  Full subcategory of $\lrt$ whose objects are schemes.\\ \\
\textbf{OBJECTS} &\\ \\ 

 $\Uij{i,j}{i,k}$& $\subsm{U}{i,j} \cap \subsm{U}{i,k}$ where $\subsm{U}{i,j}$ and $\subsm{U}{i,k}$ are open sets over some topological space $X$;\\

   $\subsm{\coprod\nolimits}{i\! \in \! \rm{I}} \subsm{A}{i}$ &  Coproduct of the family $\subsm{(\subsm{A}{i})}{i \! \in \! \rm{I}}$ where $A_i \in \Co{\mathbf{C}}$ for all $i\in \mathrm{I}$ and $\mathbf{C}$ is a category;

\end{tabular}
\end{center}

\begin{center} 
  \begin{tabular}{l p{10cm} }

   $\rgdt{R}$ & Underlying topological space for some ringed topological space $\mathfrak{R}$;\\
  $\rgds{R}$ & Underlying sheaf over  $\rgdt{R}$ for some ringed topological space $\mathfrak{R}$;\\
 $\limi{\mathbf{G}}$ & An arbitrarily chosen terminal cone via the axiom of choice in the category of cones over some diagram $\mathbf{G}$;\\
  $(\subsm{Q}{\mathbf{G}}, \subsm{\iota}{\subsm{Q}{\mathbf{G}}})$ & Standard representative of the limit of $\mathbf{G}$ where $\mathbf{G}$ is an $\topo^{\op}$-gluing data functor and $\subsm{Q}{\mathbf{G}}$ is a glued-up $\topo^{\op}$-object along $\mathbf{G}$ through $\subsm{\iota}{\subsm{Q}{\mathbf{G}}}^{\op}$ (see Definition \ref{deflema}); \\
  $((U, \subsm{\mathcal{L}}{\mathbf{G}}) , ( \subsm{\mathfrak{i}}{\mathcal{U}}^{\op}, 
  \subsm{\pi}{\subsm{\mathcal{L}}{\mathbf{G}}}))$ & Standard representative of the limit of $\mathbf{G}$ where $\mathbf{G}$ is a $\mathbf{E}\Psh{\Ops{X}}{Ab}$ (resp. $\mathbf{E}\Sh{\Ops{X}}{Ab}$)-gluing data functor and $(U , \subsm{\mathcal{L}}{\mathbf{G}})$ is a glued-up $\mathbf{E}\Psh{\Ops{X}}{Ab}$ (resp. $\mathbf{E}\Sh{\Ops{X}}{Ab}$)-object along $\mathbf{G}$ through $( \subsm{\mathfrak{i}}{\mathcal{U}}^{\op}, 
  \subsm{\pi}{\subsm{\mathcal{L}}{\mathbf{G}}})$ (see Definition \ref{stlemma});\\
  $(( \subsm{Q}{\gdt{G}} , \subsm{\mathcal{L}}{\gds{G}}) , {(\mathfrak{i}_{\mathcal{U}}^{\op}, 
  \subsm{\pi}{\subsm{\mathcal{L}}{\gds{G}}}}))$ & Standard representative of the limit of $\mathbf{G}$ where $\mathbf{G}$ is a $\rt$ (or $\lrt$ or $\sch$)-gluing data functor and $( \subsm{Q}{\gdt{G}} , \subsm{\mathcal{L}}{\gds{G}}) $ is a glued-up $\rt$ (or $\lrt$ or $\sch$)-object along $\mathbf{G}$ through ${(\mathfrak{i}_{\mathcal{U}}^{\op}, 
  \subsm{\pi}{\subsm{\mathcal{L}}{\gds{G}}}})$ (see Definition \ref{glusch}).
  \\\\
\textbf{MORPHISMS}&\\\\
$V\subsm{\subseteq}{\operatorname{op}}U$ &  $V$ open subset of $U$;\\

 $\iuv{V}{U}$ &  Inclusion map from open set $V$ to $U$; \\
 
 $\dt{\Phi}$ & Underlying morphism in ${\Cn{\topo}^{\op}}$ of some ringed topological space morphism $\Phi$;  \\
 
 $\ds{\Phi}$ & Underlying morphism of sheaves of some ringed topological space morphism $\Phi$; \\
 
$ \uiuv{\mathfrak{R}}{U}{V}$ & $\big(\iuv{V}{U}^{\op}, \subsm{(\dindsi{\mathfrak{R}}{\mathbf{Sh}}{1}\!(\iuv{W\cap V}{W\cap U}^{\op}))}{W\!\in\! \dindi{\mbf{Ops}}{U}{0}}\big)$ where $V\subsm{\subseteq}{\operatorname{op}} U\subsm{\subseteq}{\operatorname{op}} \rgdt{R}$ and $\mathfrak{R}$ ringed topological space;\\

   $\mathcal{F}$ & Functor from $\mbf{C}$ to $\mbf{D}$ where $\mbf{C}$ and $\mbf{D}$ are categories. That is, a pair $(\Co{\mathcal{F}},\Cn{\mathcal{F}})$ where $\Co{\mathcal{F}}$ is the map from $\Co{\mbf{C}}$ to $\Co{\mathbf{D}}$ and $\Cn{\mathcal{F}}$ is the map from $\Cn{\mbf{C}}$ to $\Cn{\mbf{D}}$ satisfying $\text{dom}( \Cn{\mathcal{F}}(\iuvop{V}{U}{\op})=\Co{\mathcal{F}}(U)$ and $\text{codm}( \Cn{\mathcal{F}}(\iuvop{V}{U}{\op})=\Co{\mathcal{F}}(V)$, $\Cn{\mathcal{F}} ( \subsm{\operatorname{id}}{A}) =\subsm{\operatorname{id}}{\Co{\mathcal{F}}(A)}$ and $\Cn{\mathcal{F}} ( f \circ g )=\Cn{\mathcal{F}} (f)\circ \Cn{\mathcal{F}} ( g)$ where $A\in \Co{\mbf{C}}$ and $f,g \in \Cn{\mbf{C}}$;\\
   
 $\subsm{\mathcal{F}|}{U}$ & Restriction to $\Ops{U}$ of a sheaf $\mathcal{F}$ on a topological space $X$, where $U \subsm{\subseteq}{\operatorname{op}} X$;\\
 
 $ \subsm{\mathfrak{R}|}{U}$ & Ringed topological space $(U , \subsm{\rgds{R}|}{U})$ where $\mathfrak{R}$ is a ringed  topological space and $U \subsm{\subseteq}{\operatorname{op}}{\rgdt{R}}$;\\
 
$f^{\op}$ & The morphism corresponding to a morphism $f$ of $\mathbf{C}$ in the opposite category $\mathbf{C}^{\op}$;\\

$\subsm{s|}{V}$ & $ \subsm{\mathcal{F}}{1}(\iuv{V}{ U}^{\op})(s) $ where $s \in \Co{\mathcal{F}}(U)$, $V\subop U \subop X$ and $\mathcal{F}$ is a presheaf over a topological space $X$ (when the sheaf $\mathcal{F}$ and $U$ are clear from the context).
\end{tabular}
\end{center}

\newpage

\section{Gluing objects categorically} 

In this section, we introduce the Gluing Data Categories and the Gluing Data Functors, novel constructions that can serve as tools for understanding and formalizing gluing operations across various mathematical contexts. Encapsulating gluing data within the gluing data functor framework reveals that achieving a limit over this gluing data functor, which represents the gluing data, directly corresponds to the concept of a glued-up object. 

This newfound clarity in understanding gluing data emerges as we establish subsequent lemmas while developing each definition. These lemmas gradually refine our comprehension, enabling us to encapsulate the core essence of gluing objects in a category in a more abstract and unified manner. 

The insights obtained from these constructions will find practical application in the subsequent section, particularly within the context of topological spaces, (pre)sheaves, topological spaces endowed with the structure of a (pre)sheaves and schemes. We refer to  \cite{adamek2009abstract},  \cite{Joji} and \cite{maclane} for categorical background.

\subsection{Gluing Data Categories and Gluing Data Functors}

In the definition of Gluing data functor, pushouts will have a special place and we will need the concept of pushout morphisms as defined below.

\begin{definition}
We say that a morphism is a \textbf{\textit{pushout morphism}} to refer to any canonical morphism from a component of a pushout to this pushout, giving the pushout its structure of cone.
\end{definition}

In the following definition, we present the Gluing Data category, denoted as $\glI{I}$. This category encapsulates fundamental components required to establish a flexible gluing framework applicable to diverse categories. This conceptualization is inspired by the notion of topological gluing data. Our aim is to comprehensively describe gluing data in a categorical manner. Our specific objective is to define a functor from the Gluing Data category with certain categorical properties. This functor will be constructed in such a way that limits consistently exist over these functors in the categories where gluing is usually used. Furthermore, this definition of a functor will be precisely equivalent to providing a gluing data within the conventional categories where gluing is defined.

In order to gain some intuition, one can think 
\begin{itemize} 
\item each individual index corresponds to an open set, 
\item pairs of indices represent opens contained within the single index open that serves as the patching for the gluing process, 
\item triples of indices denote the intersection of these double index opens, and 
\item The relation among these indices is designed to eliminate unnecessary canonical isomorphisms.
\end{itemize} 
The morphisms within $\glI{I}$ align with the respective inclusion maps, and their uniqueness guarantees the satisfaction of essential gluing data relations. Moreover, the equivalence relation applied to the indices within $\glI{I}$ reflects analogous relations encountered in the category of topological gluing data. Further insight into the more intricate and technical properties of these gluing index categories can be gleaned by carefully examining the retrieval of the gluing data through the framework of diagrams.

\begin{definition}\label{gic}
	Let ${\mathrm{I}}$ be an index set. We define the \textbf{\textit{Gluing Data category of type}} $\mathrm{I}$ denoted $\glI{I}$ as follows:
	
	\begin{enumerate}
		\item \textbf{Objects}: The equivalence classes of elements in $\mathrm{I} \cup \mathrm{I}^2\cup \mathrm{I}^3$ with respect to an equivalence relation $\Rel{\mathrm{I}}$ where $\Rel{\mathrm{I}}$ is the equivalence relation, generated by the relations $i \Rel{\mathrm{I}} (i,i)$, $i \Rel{\mathrm{I}} (i,i,i)$, $(i,j, k)\Rel{\mathrm{I}} (i, k, j)$, $(i,i, j) \Rel{\mathrm{I}} (i, j)$ and $(i,j, j) \Rel{\mathrm{I}} (i, j)$. We denote the equivalence class of the element $i$ as $[i]$, the equivalence class of the element $(i,j,k)$ as $[i,j,k] $ and the equivalence class of the element $(i,j)$ as $[i,j] $, for all $i,j,k \in \mathrm{I}$.
		\item \textbf{Morphisms}: Morphisms in $\glI{I}$ are structured as follows:
		
		\begin{enumerate}
			\item For each $a,b\in \Co{\glI{I}}$, a morphism from $a$ to $b$ is unique when it exists.
			
			\item Each object $a \in \Co{\glI{I}}$ has an associated identity map.
			
			\item Two morphisms in $\glI{I}$, $f: a \rightarrow b$ and $g : c \rightarrow d$, are composable as usual.
			
			\item For each $(i,j) \in \mathrm{I}^2$, we have:
			
			\begin{enumerate}
				\item A unique morphism from $i$ to $[i,j]$, denoted $\etaij{i}{j}$.
				
				\item A unique morphism from $[j,i]$ to $[i,j]$, denoted $\tauij{i}{j}$.
			\end{enumerate}
			
			\item For each $(i,j,k)\in \mathrm{I}^3$ and $n \in \{j,k\}$, we have:
			
			\begin{enumerate}
				\item A unique morphism from $[i,n]$ to $[i,j,k]$, denoted $\etaijk{n}{i}{j}{k}$.
				
				\item A unique morphism from $[j,i,k]$ to $[i,j,k]$, denoted $\tauijk{k}{i}{j}{k}$.
			\end{enumerate}
		\end{enumerate}
	\end{enumerate}
\end{definition}

\begin{figure}[H]
$$ {\tiny\xymatrixrowsep{0.2in}
\xymatrixcolsep{0.1in} \xymatrix{ 
&&&&& i\ar[rrdd]|-{\etaij{i}{k}} \ar[lldd]|-{\etaij{i}{j}}&&&&& \\
&&&&&&&&&& \\
&&&[i,j]  \ar@/^1pc/[ddddll]|-{\tauji{j}{i}} \ar[rdd]|-{\tiny{\etaijk{j}{i}{j}{k}}} & &\circlearrowleft&   & [i,k]  \ar@/^1pc/[ddddrr]|-{\tauji{k}{i}} \ar[ldd]|-{\tiny{\etaijk{k}{i}{k}{j}}}&&& \\
&&&&&&&&&& \\
&&&\circlearrowleft&[i,j,k]\ar@/_1pc/[ldd]|-{\tauijk{k}{j}{i}{j}}   \ar@{=}[rr] & & [i,k,j] \ar@/_1pc/[rdd]|-{\tauijk{j}{k}{i}{j}}  & \circlearrowleft & &&\\  
&& &&&  &&&& \\
  &  [j,i]\ar@/^1pc/[uuuurr]|-{\tauji{i}{j}} \ar[rr]|-{\scalebox{0.6}{$\etaijk{i}{j}{i}{k}$} }&&\ar@/_1pc/[ruu]|-{\tauijk{k}{i}{j}{k}} [j,i,k]    &  &\circlearrowleft & & [k,i,j] \ar@/_1pc/[luu]|-{\tauijk{j}{i}{k}{j}}  
&&  [k,i] \ar@/^1pc/[uuuull]|-{\tauji{i}{k}} \ar[ll]|- {\scalebox{0.6}{$\etaijk{i}{k}{i}{j}$} }&  \\
&&&&&&&&&& \\
&&\circlearrowleft&& [j,k,i]\ar@{=}[uul]  \ar@/_1pc/[rr]|-{\tauijk{i}{k}{j}{i}} && \ar@/_1pc/[ll]|-{\tauijk{i}{j}{k}{i}} [k,j,i] \ar@{=}[uur] &&\circlearrowleft& &\\
j\ar[uuur]|-{\etaij{j}{i}}\ar[rrrd]|-{\etaij{j}{k}}&&&&&\circlearrowleft&&&&&k \ar[llld]|-{\etaij{k}{j}} \ar[uuul]|-{\etaij{k}{i}} \\
&&&[j,k]\ar[ruu]|-{\etaijk{k}{j}{k}{i}} \ar@/^1pc/[rrrr]|-{{\begin{array}{@{}c@{}} \\ 
\tauji{k}{j}  \end{array}}}& &&&\ar@/^1pc/[llll]|-{\tauji{j}{k}} [k,j]\ar[luu]|-{\etaijk{j}{k}{j}{i}} &&&\\
}}$$
\caption{Diagram representation of a gluing Data category
 of type $ \mathrm{I}=\{ i, j, k\}$. The identity maps on every element are not presented in this diagram for the sake of clarity. The arrows in both directions each compose into the identity map.\\}
\end{figure}

In light of the preceding definition, we encounter fundamental properties of the Gluing Data category, $\glI{I}$, which plays a pivotal role in understanding its structure. These properties can be succinctly summarized by observing that for any two objects $a$ and $b$ within $\glI{I}$, the presence of a unique morphism from $a$ to $b$ is a hallmark characteristic.  Within this context, we uncover a series of intricate relationships that link various morphisms and objects in $\glI{I}$, revealing the rich interplay between its elements. These relationships are encapsulated in the ensuing Remark \ref{reglue}, which offers a deeper insight into the consequences of uniqueness within $\Cn{\glI{I}}$.

\begin{remark}\label{reglue}
	The construction of $\glI{I}$ ensures that for any two objects $a,b \in \Co{\glI{I}}$, a morphism from $a$ to $b$ is unique. This uniqueness leads to several essential properties:
	
	\begin{enumerate}
		\item For all $i,j,k \in \mathrm{I}$:
		\begin{enumerate}
		\item $\etaij{i}{i}=\tauij{i}{i}= \subsm{\operatorname{id}}{i}$;
			\item $\tauij{i}{j}\circc \tauji{i}{j} = \idij{i}{j}$ and $\tauji{i}{j}\circc \tauij{i}{j} =  \idij{j}{i}$;
			\item $\tauijk{k}{i}{j}{k} \circc \tauijk{i}{j}{k}{i} = \tauijk{j}{i}{k}{j}$ and $\tauijk{k}{i}{j}{k} \circc \tauijk{k}{j}{i}{k} =\idijk{i}{j}{k}$;
			\item $\etaijk{j}{i}{j}{k} \circc\etaij{i}{j}= \etaijk{k}{i}{j}{k} \circc \etaij{i}{k}$;
			\item $\tauijk{k}{i}{j}{k} \circc \etaijk{i}{j}{i}{k}=  \etaijk{j}{i}{j}{k} \circc \tauij{i}{j}$.
		\end{enumerate}
		
		\item For all $i,j,k \in \mathrm{I}$,  $[i,j,k] $ is a pushout with respect to $\etaij{i}{j}$ and $\etaij{i}{k}$. This follows from the equality $\etaijk{j}{i}{j}{k}  \circc\etaij{i}{j}= \etaijk{k}{i}{j}{k} \circc \etaij{i}{k}$ and the uniqueness of morphisms.
	\end{enumerate}
\end{remark}

We introduce the concept of a $\mathbf{C}$-gluing data functor, a fundamental notion that plays a pivotal role in understanding gluing operations categorically. The existence of a limit determines whether a $\mathbf{C}$-gluing data functor is sufficient to construct a glued-up object within $\mathbf{C}$. Consequently, we introduce this concept. This definition lays the foundation for understanding how gluing operations are realized within a given category.

\begin{definition}\label{limext}
	Let $\mathrm{I}$ be an index set and $\mathbf{C}$ be a category.
	\begin{enumerate} 
	\item We define a $\mathbf{C}$-\textbf{\textit{gluing data functor}} $\mathbf{G}$ of type $\mathrm{I}$ to be a functor from $\glI{I}$ to $\mathbf{C}$ sending, for any distinct indices $i,j, k \in \mathrm{I}$,
	\begin{itemize} 
	\item the pushout $(i,j,k)$ to a corresponding pushout in $\mathbf{C}$, and 
	\item the image of pushout morphisms (i.e $\etaijk{k}{i}{j}{k}$, for any $i,j, k \in \mathrm{I}$) is the corresponding pushout morphism in $\Cn{\mathbf{C}}$. 
	\end{itemize} 
	For any distinct indices $i,j, k \in \mathrm{I}$, we denote $\Goij{\mathbf{G}}{i}{j} := \Co{\mathbf{G}}([i,j])$ and $\Goijk{\mathbf{G}}{i}{j}{k} := \Co{\mathbf{G}}([i,j,k])$.
	
	\item If $\limi{\mathbf{G}}$ exists, then we say that $\mathbf{G}$ is a \textbf{\textit{gluable data functor}}. In this case, we say that $L$ is a \textbf{\textit{glued-up $\mathbf{C}$-object along $\mathbf{G}$ through}} $\subsm{\pi}{L}$ if $(L, \subsm{\pi}{L})$ is a cone over $\mathbf{G}$ that is isomorphic to $\limi{\mathbf{G}}$.
	\end{enumerate}
\end{definition}

It's important to emphasize the significance of the uniqueness of morphisms within the gluing Data category $\glI{I}$. This uniqueness is pivotal, as it enables us to establish crucial equations that underpin the foundations of our framework. As we examine the implications of these definitions and remarks, we gain valuable insights into how gluing data functors serve as a vital bridge between the abstract realm of category theory and the practical application of gluing operations within diverse mathematical structures.

\begin{remark}\label{rem1}
	\begin{enumerate}
		\item Since each morphism between any two objects in $\Co{\glI{I}}$ is unique, any gluing data functor $\mathbf{G}$ is faithful.
		
		\item Given two distinct indices $i,j \in \mathrm{I}$, from Remark \ref{reglue} $1.$ $a)$, applying a $\mathbf{C}$-gluing data functor $\mathbf{G}$ to both equations yields $\Gnij{\mbf{G}}{\tauij{i}{j}}\circ \Gnij{\mbf{G}}{\tauij{j}{i}}=\subsm{\operatorname{id}}{\Goij{\mbf{G}}{i}{j}}$ and $\Gnij{\mbf{G}}{\tauij{j}{i}}\circ \Gnij{\mbf{G}}{\tauij{i}{j}}=\subsm{\operatorname{id}}{\Goij{\mbf{G}}{j}{i}}$. Therefore, $\Gnij{\mbf{G}}{\tauij{i}{j}}$ is indeed an isomorphism.
		\item When $\mathrm{I}=\{i\}$ we have $\glI{I}=\{i\}$, a $\mathbf{C}$-gluing data functor is a functor of the form $$\begin{array}{cccl}\mathbf{G}&: \glI{I} & \rightarrow & \mathbf{C}\\& i& \mapsto & \Goi{\mathbf{G}}{i}. \end{array}$$ Therefore a glued-up object is simply an object of $\mathbf{C}$ isomorphic to $\Goi{\mathbf{G}}{i}$.
	\end{enumerate}
\end{remark}
In the following definition, we introduce a category type: categories that admit glued-up objects, allowing us to explore their properties and characteristics.
\begin{definition}
We say that a \textbf{category admits a glued-up object} if every gluing data functor admits a glued-up object.
\end{definition}
We will prove that $\topo^\op$ is a category admitting glued-up objects.

\begin{figure}[H]
\includegraphics[scale=0.3]{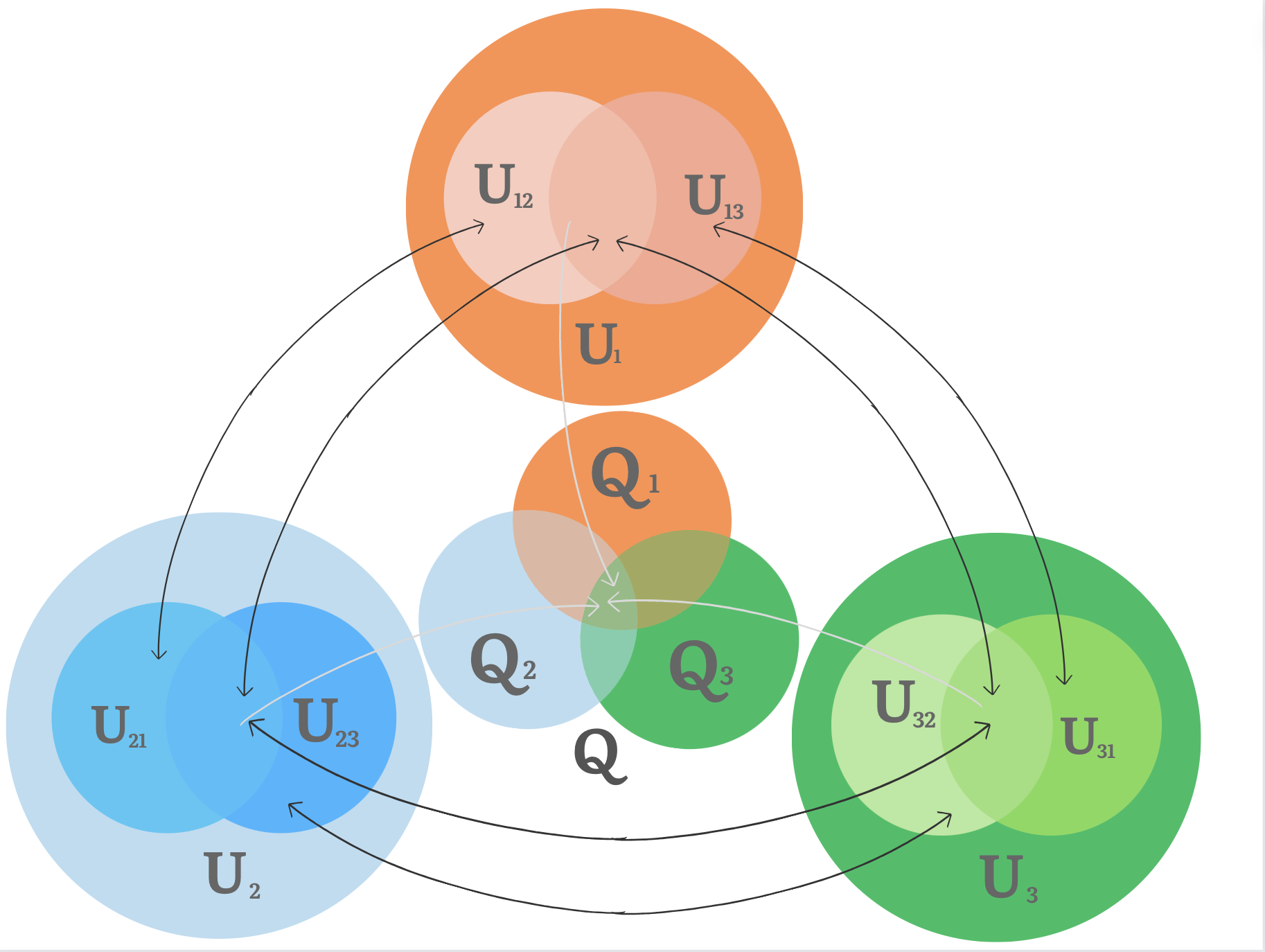}\\ \bigskip
 \caption{Representation for gluing and the glued-up object in $\mathbf{Top}^\op$ is as follows: The glued-up object $Q$ is situated in the center and is obtained from the three topological spaces, namely $U_1$, $U_2$, and $U_3$, which are mapped via the limit maps to $Q_1$, $Q_2$, and $Q_3$ respectively, forming a covering of $Q$. Moreover, $U_{12}$ is glued with $U_{21}$ and mapped to $Q_1\cap Q_2$, $U_{13}$ is glued with $U_{31}$ and mapped to $Q_1\cap Q_3$, and $U_{32}$ is glued with $U_{23}$ and mapped to $Q_2\cap Q_3$. Furthermore, all the double intersections, namely $U_{12}\cap U_{13}$, $U_{21}\cap U_{23}$, and $U_{31}\cap U_{32}$, are all mapped to the triple intersection $Q_1\cap Q_2 \cap Q_3$.
\\}
\end{figure}

\subsection{Characterizing cones and limits over gluing data functors}

In the upcoming lemma, we present a synthetic characterization of cones over the gluing data functor $\mathbf{G}$. Our goal is to establish the equivalence among three conditions that not only simplify verification but also aid in identifying limits more conveniently. These conditions play a pivotal role in understanding the properties and behavior of cones within the context of a gluing data functor.

\begin{lemma}\label{cocone} 
	Let $\mathbf{G}$ be a $\mathbf{C}$-gluing data functor. Let $N\in \Co{\mathbf{C}}$ and $\subsc{N}{\psi} : N \rightarrow \mathbf{G}$ is a family $\subsm{\left(\subsm{\subsc{N}{\psi}}{a}:N \rightarrow  \Goi{\mbf{G}}{a}\right)}{a\!\in\!  \Co{\glI{I}}}$ of morphisms in $\Cn{\mathbf{C}}$. The following statements are equivalent: 
	
	\begin{enumerate}
		\item $(N, \subsc{N}{\psi})$ is a cone over the underlying diagram of $\mathbf{G}$;
		\item $(N, \subsc{N}{\psi})$ makes the following diagrams commute, for all $i,j ,k \in \mathrm{I} $ and $n \in \{ j,k\}$:
		
		\begin{figure}[H]
			\begin{center}
				\begin{tabular}{cccl}
					{\tiny	\begin{tikzcd}[column sep=normal]&  N \arrow[]{dr}{ \randji{\subsc{N}{\psi}}{j}{i}}\arrow[swap]{dl}{\randij{\subsc{N}{\psi}}{i}{j}}  \\ \Goij{\mbf{G}}{i}{j}   && \Goij{\mbf{G}}{j}{i} \arrow[]{ll}{ \Gnij{\mbf{G}}{\tauij{i}{j}}}	\end{tikzcd}}			  	
					&	{\tiny\begin{tikzcd}[column sep=normal]&  N\arrow[swap]{dl}{\randij{\subsc{N}{\psi}}{i}{j}}  \arrow[]{dr}{ \randi{\subsc{N}{\psi}}{i}}\\ \Goij{\mbf{G}}{i}{j}   && \Goi{\mbf{G}}{i}\arrow[]{ll}{ \Gnij{\mbf{G}}{\etaij{i}{j}}}	\end{tikzcd}}	
					&   {\tiny \begin{tikzcd}[column sep=normal]&  N\arrow[]{dr}{ \randij{\subsc{N}{\psi}}{i}{n}}	\arrow[swap]{dl}{\randijk{\subsc{N}{\psi}}{i}{j}{k}}   \\\Goijk{\mbf{G}}{i}{j}{k}  && \Goij{\mbf{G}}{i}{n}\arrow[]{ll}{\Gnij{\mbf{G}}{\etaijk{n}{i}{j}{k}}} \end{tikzcd}}			  	
					\\ \tiny(a) & \tiny(b) &\tiny (c)				  	
				\end{tabular}
				\caption{Diagrams for Condition (2)}
				\label{topol10}
			\end{center}
		\end{figure}
		\item $(N, \randi{\subsc{N}{\psi}}{})$ makes the following diagrams commute, for all $i,j ,k \in \mathrm{I} $ and $n \in \{ j,k\}$:
		
		\begin{figure}[H]\begin{center}
			\begin{tabular}{cccr}
				{\tiny	\begin{tikzcd}[column sep=normal]&  N \arrow[]{dr}{ {\randi{\subsc{N}{\psi}}{j}}}	\arrow[swap]{dl}{\randij{\subsc{N}{\psi}}{i}{j}}  \\ \Goij{\mbf{G}}{i}{j}      && \Goi{\mbf{G}}{j} \arrow[]{ll}{\Gnij{\mbf{G}}{\tauij{i}{j}\circ \etaji{i}{j}}}\end{tikzcd}}			  	
				&	{\tiny\begin{tikzcd}[column sep=normal]&  N\arrow[]{dr}{ {\randi{\subsc{N}{\psi}}{i}}} \arrow[swap]{dl}{\randij{\subsc{N}{\psi}}{i}{j}} \\\Goij{\mbf{G}}{i}{j}     &&  \Goi{\mbf{G}}{i}\arrow[]{ll}{\Gnij{\mbf{G}}{\etaij{i}{j}}}	\end{tikzcd}}	
				& 	{\tiny\begin{tikzcd}[column sep=normal]&  N\arrow[]{dr}{ \randij{\subsc{N}{\psi}}{i}{n}}  \arrow[swap]{dl}{\randijk{\subsc{N}{\psi}}{i}{j}{k}}\\\Goijk{\mbf{G}}{i}{j}{k}    && \Goij{\mbf{G}}{i}{n} 	\arrow[]{ll}{\Gnij{\mbf{G}}{\etaijk{n}{i}{j}{k}}}\end{tikzcd}}
				\\ \tiny (a) &\tiny (b) & \tiny(c)
			\end{tabular}\end{center}
			\caption{Diagrams for Condition (3)}
			\label{topol11}
		\end{figure}
	\end{enumerate}
	
	\begin{proof}
		\begin{enumerate}
			\item[$(1)\Rightarrow (2) $]   This is clear. 
			
			\item[$(2)\Rightarrow (3)$] Assume that $(2)$ is satisfied. We only need to prove that $\Gnij{\mbf{G}}{\tauij{i}{j}\circ \etaji{i}{j}}\circ {\randi{\subsc{N}{\psi}}{j}}={\randij{\subsc{N}{\psi}}{i}{j}}$ for all $i,j\in \mathrm{I}$. This is a consequence of the commutativity of the diagrams in Figure \ref{topol10} $(2) (a),(b)$.  
			
			\item[$(3)\Rightarrow (1)$] Suppose that $(N, \randi{\subsc{N}{\psi}}{})$ makes the diagrams in $(3)$ commute. We want to prove that $(N, \randi{\subsc{N}{\psi}}{})$ is a cone over the underlying diagram of $\mathbf{G}$. 
			
			To show that $(N, \randi{\subsc{N}{\psi}}{})$ is a cone over the underlying diagram of $\mathbf{G}$, we need to prove that for any $a, b \in \Co{\glI{I}}$ and morphism $f: a \rightarrow b$ in $\Cn{\glI{I}}$, we have $\Gnij{\mbf{G}}{f}\circc \randi{\subsc{N}{\psi}}{a} = \randi{\subsc{N}{\psi}}{b}$. Since $(3)$ is satisfied, it suffices to prove that this equality holds for $f$ equal to $\tauij{i}{j}$ and $\tauijk{k}{i}{j}{k}$ for all $i,j,k \in \mathrm{I}$. 
			
			Let $i,j,k\in \mathrm{I}$. Since by assumption we have that $\Gnij{\mbf{G}}{\etaij{i}{j}} \circc \randi{\subsc{N}{\psi}}{i}=\randij{\subsc{N}{\psi}}{i}{j}$, then 
			\begin{align*}
				&\Gnij{\mbf{G}}{\tauij{i}{j} \circc \etaji{i}{j}}\circc  \randi{\subsc{N}{\psi}}{j}\\
				&=\Gnij{\mbf{G}}{\tauij{i}{j}}\circc \Gnij{\mbf{G}}{\etaji{i}{j}}\circc  \randi{\subsc{N}{\psi}}{j}\\
				& = \Gnij{\mbf{G}}{\tauij{i}{j}}\circc \randji{\subsc{N}{\psi}}{j}{i}. 
			\end{align*}
			
			Therefore, since $ \Gnij{\mbf{G}}{\tauij{i}{j}\circ \etaji{i}{j}}\circc  \randi{\subsc{N}{\psi}}{j}=\randij{\subsc{N}{\psi}}{i}{j}$, we obtain $\Gnij{\mbf{G}}{\tauij{i}{j}}\circc \randji{\subsc{N}{\psi}}{j}{i}=\randij{\subsc{N}{\psi}}{i}{j}$ as required.
			
			On the other hand, applying the functor $\mathbf{G}$ to the equality $\tauijk{k}{i}{j}{k} \circc \etaijk{i}{j}{i}{k}=   \etaijk{j}{i}{j}{k}  \circc \tauij{i}{j}$, we obtain $\Gnij{\mbf{G}}{\tauijk{k}{i}{j}{k} \circc \etaijk{i}{j}{i}{k}}=\Gnij{\mbf{G}}{\etaijk{j}{i}{j}{k} \circc \tauij{i}{j}}$. Therefore, we have
			\begin{align*}
				& \Gnij{\mbf{G}}{\etaijk{j}{i}{j}{k} \circc \tauij{i}{j}}\circc\randij{\subsc{N}{\psi}}{i}{j} \\
				&=(\Gnij{\mbf{G}}{\tauijk{k}{i}{j}{k}} \circc \Gnij{\mbf{G}}{\etaijk{i}{j}{i}{k}})\circc \randij{\subsc{N}{\psi}}{i}{j}\\
				&=\Gnij{\mbf{G}}{\tauijk{k}{i}{j}{k}} \circc (\Gnij{\mbf{G}}{\etaijk{i}{j}{i}{k}}\circc \randij{\subsc{N}{\psi}}{i}{j})\\
				&=\Gnij{\mbf{G}}{\tauijk{k}{i}{j}{k}} \circc \randijk{\subsc{N}{\psi}}{j}{i}{k}.
			\end{align*}
			
			Moreover, since by assumption $\Gnij{\mbf{G}}{\etaijk{j}{i}{j}{k} \circ \tauij{i}{j}}\circ \randji{\subsc{N}{\psi}}{j}{i}=\randijk{\subsc{N}{\psi}}{i}{j}{k}$, we obtain $\Gnij{\mbf{G}}{\tauijk{k}{i}{j}{k}} \circ \randijk{\subsc{N}{\psi}}{j}{i}{k}=\randijk{\subsc{N}{\psi}}{i}{j}{k}$ as required. This proves that $(N, \randi{\subsc{N}{\psi}}{})$ is a cone over the underlying diagram of $\mathbf{G}$.
		\end{enumerate}
	\end{proof}
\end{lemma}

The upcoming theorem provides a characterization of a glued-up object in the form of a pullback diagram. It serves as a fundamental result that elucidates the criteria for designating an object within a category as a "glued-up object"  along a specific gluing data functor.

\begin{theorem}   \label{preglued}
	Let $\mathbf{G}$ be a gluing data functor, $L\in \Co{\mathbf{C}}$ and $\subsm{\pi}{L}$ is a family $\subsm{(\dindi{\pi}{L}{a})}{a\! \in\! \Co{\glI{I}}}$ of morphisms $\dindi{\pi}{L}{a}: L\rightarrow \Goi{\mbf{G}}{a}$ for all $a \in \Co{\glI{I}}$. $L$ is a  glued-up $\mathbf{C}$-object  along $\mathbf{G}$ through $ \subsm{\pi}{L}$ if and only if for all $i,j, k \in \mathrm{I}$ and  $n \in \{ j,k\}$, the following properties are satisfied:
	\begin{enumerate}
		\item $\dindij{\pi}{L}{i}{j} = \Gnij{\mbf{G}}{\etaij{i}{j}}\circc \dindi{\pi}{L}{i}$ ;  
		\item $\dindijk{\pi}{L}{i}{j}{k}= \Gnij{\mbf{G}}{\etaijk{n}{i}{j}{k}}\circc \dindij{\pi}{L}{i}{n} $;
		\item $(L, \subsm{(\dindi{\pi}{L}{i})}{i\!\in\! \mathrm{I}})$ is the limit over the pullback diagram defined by the morphisms
		$$\Gnij{\mbf{G}}{\etaij{i}{j}}: \Goi{\mbf{G}}{i} \rightarrow \Goij{\mbf{G}}{i}{j} \text{ and } \Gnij{\mbf{G}}{\tauij{i}{j}\circc \etaji{i}{j}}: \Goi{\mbf{G}}{j} \rightarrow \Goij{\mbf{G}}{i}{j}.$$
	\end{enumerate}  
	\begin{proof}
		Let $L$ be a glued-up $\mathbf{C}$-object constructed by gluing along $\mathbf{G}$ through the morphisms $\dindi{\pi}{L}{}$. According to Definition \ref{limext}, we know that $\limi{\mathbf{G}}$ exists and is isomorphic to $(L, \subsm{\pi}{L})$. Our objective is to demonstrate the satisfaction of properties $(1)$, $(2)$, and $(3)$ as outlined in the theorem.
		
		The diagrams in Lemma \ref{cocone} $(3)$ $(b)$ and $(c)$ are equivalent to properties $(1)$ and $(2)$, respectively. To establish property $(3)$ of the theorem, we start by proving that the following diagram for arbitrary $i$ and $j$ in $\mathrm{I}$:
		\begin{figure}[H] 
			\begin{center}
				{\tiny	\begin{tikzcd}[column sep=normal]
						&\Goij{\mbf{G}}{i}{j}  &&& \Goi{\mbf{G}}{j}\arrow[swap]{lll}{\Gnij{\mbf{G}}{\tauij{i}{j}\circ \etaji{i}{j}}}\\
						&\Goi{\mbf{G}}{i} \arrow[]{u}{\Gnij{\mbf{G}}{\etaij{i}{j}}}&& & L\arrow[]{lll}{ \dindi{\pi}{L}{i}}\arrow[swap]{u} { \dindi{\pi}{L}{j}}
			\end{tikzcd}}\end{center} 	 \caption{}\label{topol13}	
		\end{figure}    
		\noindent commutes. This holds by combining the commutativity of the diagram in Figure (3) (a) and (3) (b) of Lemma \ref{cocone}. 
		
		Let us now consider the scenario where we have a pair $(L',\subsm{\pi}{L'})$, with $L'\in \Co{\mathbf{C}}$ and $\subsm{\pi}{L'}$ representing a family of maps $\dindpi{\pi}{L}{i} : L'\rightarrow \Goi{\mbf{G}}{a}$ where $a\in  \Co{\glI{I}}$, such that the following diagram commutes for all $i,j\in \mathrm{I}$:
		
		\begin{figure}[H] 
			\begin{center}
				{\tiny	
					\begin{tikzcd}[column sep=normal]
						&\Goij{\mbf{G}}{i}{j} & &&  \Goi{\mbf{G}}{j}\arrow[swap]{lll}{{\Gnij{\mbf{G}}{\tauij{i}{j}\circ \etaji{i}{j}}}}\\
						& \Goi{\mbf{G}}{i} \arrow[]{u}{\Gnij{\mbf{G}}{\etaij{i}{j}}}&& &  L'\arrow[]{lll}{\dindpi{\pi}{L}{i}}\arrow[swap]{u} {\dindpi{\pi}{L}{j}}
					\end{tikzcd}
				}
			\end{center} 
			\caption{}
			\label{topol116}
		\end{figure}
		
		Now, our aim is to establish that the pair $(L',\widetilde{\subsm{\pi}{L'}})$ forms a cone over $\mathbf{G}$, where $\widetilde{\subsm{\pi}{L'}}$ is a family of maps $\widetilde{\dindpi{\pi}{L}{a}}: L'\rightarrow \Goi{\mbf{G}}{a}$, and each $\widetilde{\dindpi{\pi}{L}{i}}$ coincides with $\dindpi{\pi}{L}{i}$, $\widetilde{\dindpij{\pi}{L}{i}{j}}$ is given by ${\Gnij{\mbf{G}}{\etaij{i}{j}}}\circc \dindpi{\pi}{L}{i}$, and $ \widetilde{\dindpijk{\pi}{L}{i}{j}{k}}$ is computed as $\Gnij{\mbf{G}}{\etaijk{n}{i}{j}{k}}\circc {\dindpij{\pi}{L}{i}{n}}$, for all $i,j, k \in \mathrm{I}$ and  $n \in \{ j,k\}$.
		
		In order to prove that $(L',\widetilde{\subsm{\pi}{L'}})$ is a cone over $\mathbf{G}$, according to the definition of $\widetilde{\subsm{\pi}{L'}}$ and Lemma \ref{cocone}, it is enough to prove that the following diagram commutes, for all $i,j\in \mathrm{I}$: 
		\begin{figure}[H]
			\begin{center}
				{\tiny	\begin{tikzcd}[column sep=normal]&  L' \arrow[]{dr}{\widetilde{\dindpi{\pi}{L}{j}}}	\arrow[swap]{dl}{\widetilde{\dindpij{\pi}{L}{i}{j}}}  \\ \Goij{\mbf{G}}{i}{j}    && \Goi{\mbf{G}}{j} \arrow[]{ll}{\Gnij{\mbf{G}}{\tauij{i}{j}\circ \etaji{i}{j}}}\end{tikzcd}}
			\end{center}
			\caption{}\label{topol126}
		\end{figure}
		
		\noindent 
		The commutativity of this diagram follows directly from the definition of $\widetilde{\subsm{\pi}{L'}}$ and the commutativity of the diagram in Figure \ref{topol116}. Therefore, the pair $(L',\widetilde{\subsm{\pi}{L'}})$ is a cone over $\mathbf{G}$. By assumption, $\limi{\mathbf{G}}\simeq (L,\subsm{\pi}{L})$, so there exists a unique morphism, say $\mu: L'\rightarrow L$, making each of the following diagrams commute, for all $i,j, k \in \mathrm{I}$ and  $n \in \{ j,k\}$:
		\begin{figure}[H]\begin{center}	\begin{tabular}{cccl}
				{\tiny\begin{tikzcd}[column sep=normal] &L'\arrow[bend left, labels=description]{dddr}{\widetilde{\dindpi{\pi}{L}{j}}} \arrow[dashed]{dd}{\exists ! \mu}\arrow[bend right, labels=description]{dddl}{\widetilde{\dindpij{\pi}{L}{i}{j}}} \\\\ & L\arrow[labels=description]{dr}{{\subsm{\pi}{\subsm{L}{j}}}} \arrow[labels=description]{dl}{\dindij{\pi}{L}{i}{j}}    \\ \Goij{\mbf{G}}{i}{j} & & \Goi{\mbf{G}}{j}\arrow[]{ll}{{\Gnij{\mbf{G}}{\tauij{i}{j}\circ \etaji{i}{j}}}}   
				\end{tikzcd}} &	{\tiny\begin{tikzcd}[row sep=normal] & L'\arrow[bend left, labels=description]{dddr}{\widetilde{\dindpi{\pi}{L}{i}}} \arrow[bend right, labels=description]{dddl}{\widetilde{\dindpij{\pi}{L}{i}{j}}} \arrow[dashed]{dd}{\exists ! \mu}  \\\\ & L \arrow[labels=description]{dr}{\dindi{\pi}{L}{i}}\arrow[labels=description]{dl}{\dindij{\pi}{L}{i}{j}}    \\ \Goij{\mbf{G}}{i}{j}  & & \Goi{\mbf{G}}{i} \arrow[]{ll}{\Gnij{\mbf{G}}{\etaij{i}{j}}}
				\end{tikzcd}} &  {\tiny\begin{tikzcd}[row sep=normal] & L'\arrow[bend left, labels=description]{dddr}{\widetilde{\dindpij{\pi}{L}{i}{n}}} \arrow[bend right, labels=description]{dddl}{\widetilde{\dindpijk{\pi}{L}{i}{j}{k}}} \arrow[dashed]{dd}{\exists ! \mu}     \\\\ &L \arrow[labels=description]{dr}{\subsm{\pi}{\subsm{L}{[i,n]}}} \arrow[labels=description]{dl}{\dindijk{\pi}{L}{i}{j}{k}} \\ \Goijk{\mbf{G}}{i}{j}{k} & & \Goij{\mbf{G}}{i}{n} \arrow[]{ll}{\Gnij{\mbf{G}}{\etaijk{n}{i}{j}{k}}}
				\end{tikzcd}}\\ \tiny (a) &\tiny (b) &\tiny (c) \end{tabular} \end{center}\caption{}\label{topol118}           \end{figure} 
		
		We now pick such a $\mu$ and by construction of $\mu$ the following diagram

		\begin{figure}[H]
			\begin{center}
				{\tiny\begin{tikzcd}[column sep=normal]
						\Goij{\mbf{G}}{i}{j}  && & 
						\Goi{\mbf{G}}{j} \arrow[swap]{lll}{{\Gnij{\mbf{G}}{\tauij{i}{j}\circ \etaji{i}{j}}}}  \\
						\Goi{\mbf{G}}{i}\arrow[]{u}{\Gnij{\mbf{G}}{\etaij{i}{j}}}     && & 
						L\arrow[swap]{u}{\subsm{\pi}{\subsm{L}{j}}}\arrow[]{lll}{\dindi{\pi}{L}{i}} \\
						&&&& L' \arrow[swap,dashed]{ul}{ \mu}\arrow[bend left=10]{ullll}{{ {\dindpi{\pi}{L}{i}}}} \arrow[swap,bend right=15]{uul}{ \dindpi{\pi}{L}{j}}
				\end{tikzcd} }
			\end{center}\caption{}\label{topol117}
		\end{figure}
		\noindent also commutes, for all $i,j\in \mathrm{I}$. This concludes the proof. 
		
		Conversely, suppose that properties $(1)$, $(2)$ and $(3)$ of the theorem are satisfied. We want to prove that $L$ is a glued-up $\mathbf{C}$-object along $\mathbf{G}$ through $\subsm{\pi}{L}$. From the properties $(1)$, $(2)$ and $(3)$ of the theorem, we have that $(L , \subsm{\pi}{L})$ is a cone over $\mathbf{G}$.  
		Suppose that $(L' , \subsm{\pi}{L'})$ is another cone over $\mathbf{G}$. That is, the following diagrams
		
		\begin{figure}[H] \begin{center}\begin{tabular}{cccr}
				{\tiny	\begin{tikzcd}[column sep=normal]&  L' \arrow[]{dr}{\dindpi{\pi}{L}{j}}	\arrow[swap]{dl}{\dindpij{\pi}{L}{i}{j}}  \\ \Goij{\mbf{G}}{i}{j}    && \Goi{\mbf{G}}{j} \arrow[]{ll}{{\Gnij{\mbf{G}}{\tauij{i}{j}\circ \etaji{i}{j}}}}\end{tikzcd}}			  	
				&	{\tiny\begin{tikzcd}[column sep=normal]&  L'\arrow[]{dr}{{\dindpi{\pi}{L}{i}}} \arrow[swap]{dl}{\dindpij{\pi}{L}{i}{j}} \\\Goij{\mbf{G}}{i}{j}    &&\Goi{\mbf{G}}{i}\arrow[]{ll}{\Gnij{\mbf{G}}{\etaij{i}{j}}}	\end{tikzcd}}				  				         & 	{\tiny\begin{tikzcd}[column sep=normal]&  L'\arrow[]{dr} {\dindpij{\pi}{L}{i}{n}}  \arrow[swap]{dl}{{\dindpijk{\pi}{L}{i}{j}{k}}}\\\Goijk{\mbf{G}}{i}{j}{k}    && \Goij{\mbf{G}}{i}{n} 	\arrow[]{ll}{\Gnij{\mbf{G}}{\etaijk{n}{i}{j}{k}}}\end{tikzcd}} \\\tiny (a) &\tiny (b) & \tiny(c)		\end{tabular}	\end{center}	  	
			\caption{}\label{topol120}				  	\end{figure}
		\noindent   commute, for all $i,j, k \in \mathrm{I}$ and  $n \in \{ j,k\}$. Now, we want to show that the pair $(L' , \subsm{\pi}{L'})$ makes the following diagram
		\begin{figure}[H] 
			\begin{center}
				{\tiny	\begin{tikzcd}[column sep=normal]
						&\Goij{\mbf{G}}{i}{j}  & && \Goi{\mbf{G}}{j}\arrow[swap]{lll}{{\Gnij{\mbf{G}}{\tauij{i}{j}\circ \etaji{i}{j}}}}\\
						& \Goi{\mbf{G}}{i} \arrow[]{u}{\Gnij{\mbf{G}}{\etaij{i}{j}}}&& & L'\arrow[]{lll}{{ \dindpi{\pi}{L}{i}}}\arrow[swap]{u} {{ \dindpi{\pi}{L}{j}}}
			\end{tikzcd}}\end{center} 	 \caption{}\label{topol127}	
		\end{figure}
		\noindent commute, for all $i,j \in \mathrm{I}$. By Figure \ref{topol120} $(a)$ and $(b)$, we have 
		$$\Gnij{\mbf{G}}{\tauij{i}{j}\circ \etaji{i}{j}}\circ \dindpi{\pi}{L}{j}={\Gnij{\mbf{G}}{\tauij{i}{j}}\circ\big(\Gnij{\mbf{G}}{\etaji{i}{j}}}\circ \dindpi{\pi}{L}{j}\big)=\Gnij{\mbf{G}}{\tauij{i}{j}}\circ\dindpij{\pi}{L}{j}{i} =\dindpij{\pi}{L}{i}{j}=\Gnij{\mbf{G}}{\etaij{i}{j}}\circ {\dindpi{\pi}{L}{i}}.$$ Hence the diagram in Figure \ref{topol127} commutes.   Then, by Property $(3)$ of the theorem, there exists a unique morphism $\mu: L'\rightarrow L$ making the following diagram 
		\begin{figure}[H]
			\begin{center}
				{\tiny\begin{tikzcd}[column sep=normal]
						\Goij{\mbf{G}}{i}{j}  && & 
						\Goi{\mbf{G}}{j} \arrow[swap]{lll}{{\Gnij{\mbf{G}}{\tauij{i}{j}\circ \etaji{i}{j}}}}  \\
						\Goi{\mbf{G}}{i}\arrow[]{u}{\Gnij{\mbf{G}}{\etaij{i}{j}}}     & && 
						L\arrow[swap]{u}{{\subsm{\pi}{\subsm{L}{j}}}}\arrow[]{lll}{\dindi{\pi}{L}{i}} \\
						&&&& L' \arrow[swap,dashed, near end]{ul}{ \exists ! \mu}\arrow[bend left=15]{ullll}{{\dindpi{\pi}{L}{\!i}}} \arrow[swap,bend right=15]{uul}{\dindpi{\pi}{L}{\!j}}
				\end{tikzcd} }
			\end{center}\caption{}\label{topol128}
		\end{figure}
		\noindent   commute, for all $i,j\in \mathrm{I}$. 
		We choose such a $\mu$. Thanks to the commutativity of the diagram as shown in Figure \ref{topol128}, combined with Properties $(1)$ and $(2)$ of the theorem, which are assumed to be satisfied, we obtain that each of the following diagrams 
		\begin{figure}[H]\begin{center}	\begin{tabular}{cccl}
				{\tiny\begin{tikzcd}[column sep=normal] &L'\arrow[bend left, labels=description]{dddr}{\dindpi{\pi}{L}{j}} \arrow[dashed]{dd}{\exists ! \mu}\arrow[bend right, labels=description]{dddl}{\dindpij{\pi}{L}{i}{j}} \\\\ & L\arrow[labels=description]{dr}{\dindi{\pi}{L}{j}} \arrow[labels=description]{dl}{\dindij{\pi}{L}{i}{j}}    \\ \Goij{\mbf{G}}{i}{j} & & \Goi{\mbf{G}}{j}\arrow[]{ll}{{\Gnij{\mbf{G}}{\tauij{i}{j}\circ \etaji{i}{j}}}}   
				\end{tikzcd}} &	{\tiny\begin{tikzcd}[row sep=normal] & L'\arrow[bend left, labels=description]{dddr}{{\dindpi{\pi}{L}{i}}} \arrow[bend right, labels=description]{dddl}{\dindpij{\pi}{L}{i}{j}} \arrow[dashed]{dd}{\exists ! \mu}  \\\\ & L \arrow[labels=description]{dr}{\dindi{\pi}{L}{i}}\arrow[labels=description]{dl}{\dindij{\pi}{L}{i}{j}}    \\ \Goij{\mbf{G}}{i}{j}  & & \Goi{\mbf{G}}{i} \arrow[]{ll}{\Gnij{\mbf{G}}{\etaij{i}{j}}}
				\end{tikzcd}} &  {\tiny\begin{tikzcd}[row sep=normal] & L'\arrow[bend left, labels=description]{dddr}{\dindpij{\pi}{L}{i}{n}} \arrow[bend right, labels=description]{dddl}{{\dindpijk{\pi}{L}{i}{j}{k}}} \arrow[dashed]{dd}{\exists ! \mu}     \\\\ &L \arrow[labels=description]{dr}{{\subsm{\pi}{\subsm{L}{[i,n]}}}} \arrow[labels=description]{dl}{{\dindijk{\pi}{L}{i}{j}{k}}} \\ \Goijk{\mbf{G}}{i}{j}{k}  & & \Goij{\mbf{G}}{i}{n} \arrow[]{ll}{\Gnij{\mbf{G}}{\etaijk{n}{i}{j}{k}}}
				\end{tikzcd}}\\ \tiny (a) &\tiny (b) &\tiny (c) \end{tabular} \end{center}\caption{}\label{topol121}           \end{figure} 
		\noindent commutes, for all $i,j, k \in \mathrm{I}$ and  $n \in \{ j,k\}$.  
		With this, we have shown that $(L,\subsm{\pi}{L})$ is a cone over $\mathbf{G}$, and any other cone $(L',\subsm{\pi}{L'})$ can be uniquely mapped to $(L,\subsm{\pi}{L})$ via the unique morphism $\mu: L' \to L$ satisfying the commutativity of the diagram in Figure \ref{topol116}. 
		Therefore, $\limi{\mathbf{G}}$ exists and is isomorphic to $(L,\subsm{\pi}{L})$, which implies that $L$ is a glued-up $\mathbf{C}$-object along $\mathbf{G}$ through $\subsm{\pi}{L}$.
	\end{proof}
\end{theorem}

In the following remark, we translate the properties of a glued-up object within the context of a contravariant functor.
\begin{remark}\label{oppc}
\hspace{2em}
\begin{enumerate}
	\item If $\mathbf{G}$ is a gluing data functor of  type $\mathrm{I}$ from $\glI{I}$ to $\mathbf{C}^{\op}$. 
	Then, $Q$ is a glued-up $\mathbf{C}^{\op}$-object along $\mathbf{G}$ through $\subsm{\iota}{Q}^{\op}$ if and only if the following properties hold in $\mathbf{C}$, for all $i,j, k \in \mathrm{I}$ and  $n \in \{ j,k\}$:
	\begin{itemize}
		\item[(a)] $\dindij{\iota}{Q}{i}{j}=\dindi{\iota}{Q}{i} \circc \Gnij{\mbf{G}}{\etaij{i}{j}}^{\op}$;  
		\item[(b)] $\dindijk{\iota}{Q}{i}{j}{k} =\dindij{\iota}{Q}{i}{n}  \circc \Gnij{\mbf{G}}{\etaijk{n}{i}{j}{k}}^{\op} $;
		
		\item[(c)] $(Q, \subsm{(\dindi{\iota}{Q}{i})}{i \! \in \! \mathrm{I}})$ is the limit over the pushout diagram defined by the morphisms \\ $\Gnij{\mbf{G}}{\etaij{i}{j}}^{\op}: \Goij{\mbf{G}}{i}{j} \rightarrow \Goi{\mbf{G}}{i}$ and ${\Gnij{\mbf{G}}{\tauij{i}{j} \! \circ \! \etaji{i}{j}}}^{\op}:  \Goij{\mbf{G}}{i}{j} \rightarrow \Goi{\mbf{G}}{j}$. 
		\end{itemize}
		\item If $\mathrm{I}=\{1,2\}$, then $\Co{\glI{I}}=\{1,2,[1,2],[2,1]\}$, and if $\mathbf{G}$ is a $\mathbf{C}$-gluing data functor from $\glI{I}$, then a glued-up object, if it exists, is the pullback of $\Gnij{\mbf{G}}{\etaij{1}{2}}$ and $\Cn{\mathbf{G}}(\tauij{1}{2}\circ \etaji{1}{2})$. So, pullbacks are a particular case of limits of a certain gluing data functor, and therefore, gluing data functors generalize the concept of pullbacks.

	\end{enumerate}
\end{remark}

\section{Gluing topological spaces categorically}\label{joycat}
Throughout the rest of the paper, we denote $\mathbf{Top}$ as the category of topological spaces, where morphisms correspond to continuous maps. Additionally, we refer to $\mathbf{oTop}$ as a subcategory of $\mathbf{Top}$ where the objects remain the same as those in $\topo$, and the morphisms are restricted to continuous open maps. In the upcoming definition, we provide a tangible representation for the limit over $\mathbf{(o)Top}^{\op}$-gluing data functors. 
\subsection{Correspondence between topological gluing data and Gluing data functors}
We refer to \cite{munkres} for the topological background. In this section, we explore the process of gluing topological spaces using the formal language developed in the previous section. We also introduce alternative perspectives on glued-up objects that shed light on how certain required properties can be viewed as equivalent to the usual glued-up object.

In the context of our discussion on gluing in topological spaces, we now present a formal definition that generalizes the conventional concept of gluing data for a topological space. The following definition introduces the collection of the components that are necessary to allow the gluing of topological spaces.  
\begin{definition}\label{tpsgd}
	A collection $$\left(\mathrm{I},\subsm{({\subsm{U}{i}})}{i\!\in\! \mathrm{I}},\subsm{(\subsm{U}{i,j}, \subsm{\upsilon}{i,j} )}{(i,j)\! \in\! \mathrm{I}^2}, \subsm{(\subsm{\varphi}{i,j}, \phijk{k}{i}{j} )}{(i,j,k)\!\in\! \mathrm{I}^3}\right)$$ where 
	\begin{enumerate}
		\item $\mathrm{I}$ is an index set;
		\item $\subsm{(\subsm{U}{i})}{i\!\in\! \mathrm{I}}$ is a family of topological spaces;
		\item $\subsm{(\subsm{U}{i,j}, \subsm{\upsilon}{i,j} )}{(i,j)\! \in\! \mathrm{I}^2}$ is a family where $\subsm{U}{i,j}$ are topological spaces and $\subsm{\upsilon}{i,j} : \subsm{U}{i,j} \rightarrow \subsm{U}{i}$ are (open) continuous maps, for all $i, j \in \mathrm{I}$; 
		\item $\subsm{(\subsm{\varphi}{i,j}, \subsm{\subsm{\varphi}{k}}{(i,j)})}{(i,j,k)\!\in\! \mathrm{I}^3}$ is a family of (open) continuous maps where $\subsm{\varphi}{i,j} : \subsm{U}{i,j} \rightarrow \subsm{U}{j,i}$ and $\phijk{k}{i}{j}:\subsm{U}{i, j,k} \rightarrow \subsm{U}{j, i,k } $ where $\subsm{U}{i, j,k} =\subsm{U}{i, k,j} $ is an arbitrary chosen pullback of the morphism $\subsm{\upsilon}{i,j}$ and $\subsm{\upsilon}{i,k}$ and we choose $U_{i,j,j}=U_{i,i,j}:= U_{i,j}$ and $U_{i,i,i}:=U_i$, for all $(i,j,k) \in \mathrm{I}^3$.	
	\end{enumerate}
	such that, for each $i,j,k\in \mathrm{I}$, 
	\begin{enumerate}
		\item[(a)] $\subsm{U}{i,i}=\subsm{U}{i}$;
		\item[(b)] $\subsm{\varphi}{i,i}=\subsm{\operatorname{id}}{\subsm{U}{i}};$ 
		
		\item [(c)] $\phijk{j}{i}{k}=\phijk{i}{j}{k}\circc \phijk{k}{i}{j}$.
		\item[(d)] $\subsm{\mathfrak{i}}{\subsm{U}{j, i,k }, \subsm{U}{j,i}}\circ \phijk{k}{i}{j}=\subsm{\varphi}{i,j}\circ \subsm{\mathfrak{i}}{\subsm{U}{i, j,k }, \subsm{U}{i,j}}$.
	\end{enumerate}
	is called a \textbf{\textit{topological space (open) gluing data}}.
\end{definition}

\begin{remark}\label{inve}
	Given $\left(\mathrm{I},\subsm{({\subsm{U}{i}})}{i\!\in\! \mathrm{I}},\subsm{(\subsm{U}{i,j}, \upsilon_{i,j} )}{(i,j)\! \in\! \mathrm{I}^2}, \subsm{(\subsm{\varphi}{i,j}, \phijk{k}{i}{j})}{(i,j,k)\!\in\! \mathrm{I}^3}\right)$ a topological space gluing data, $\subsm{\varphi}{i,j}$ is a homeomorphism from $\subsm{U}{i, j}$ to $\subsm{U}{j,i }$ whose inverse is $\subsm{\varphi}{j,i}$. This is a consequence of $(c)$ applied to $k=i$ and $(b)$.
\end{remark}

In the upcoming lemma, we explore the connection between $\mathbf{(o)Top}^{\op}$-gluing data functors and topological space gluing data. We prove the equivalence between these two constructs, shedding light on how the abstract notion of gluing in category theory aligns with the concrete gluing of topological spaces. This insight bridges the gap between theory and practice, providing a powerful tool for both abstract mathematical exploration and practical applications in topology.
\begin{lemma}\label{equi}
	An $\mathbf{(o)Top}^{\op}$-gluing data functor $\mathbf{G}$ induces the topological space (open) gluing data $$\left(\mathrm{I}, \subsm{(\Goi{\mbf{G}}{i})}{i\!\in\! \mathrm{I}},\subsm{(\Goij{\mbf{G}}{i}{j}, \Gnij{\mbf{G}}{\etaij{i}{j}}^{\op})}{(i,j)\!\in\! \mathrm{I}^2},\subsm{(\Gnij{\mbf{G}}{\tauij{i}{j}}^{\op},\Gnij{\mbf{G}}{\tauijk{k}{i}{j}{k}}^{\op})}{( i,j,k)\! \in\! \mathrm{I}^3}\right).$$ Conversely, a topological space (open) gluing data $$\left(\mathrm{I},\subsm{({\subsm{U}{i}})}{i\!\in\! \mathrm{I}},\subsm{(\subsm{U}{i,j}, \upsilon_{i,j} )}{(i,j)\! \in\! \mathrm{I}^2}, \subsm{(\subsm{\varphi}{i,j},\phijk{k}{i}{j})}{(i,j,k)\!\in\! \mathrm{I}^3}\right)$$ induces the $\mathbf{(o)Top}^{\op}$-gluing data functor $\mathbf{G}$ defined by $\Goi{\mbf{G}}{i}:= \subsm{U}{i}$, $\Goij{\mbf{G}}{i}{j}:=\subsm{U}{i,j}$,  $\Gnij{\mbf{G}}{\etaij{i}{j}}^{\op}:=\upsilon_{i,j}, \Gnij{\mbf{G}}{\tauij{i}{j}}^{\op}:=\subsm{\varphi}{i,j}$ and $\Gnij{\mbf{G}}{\tauijk{k}{i}{j}{k}}^{\op}:=\phijk{k}{i}{j}$ for all $i, j,k\in \mathrm{I}$.
	\begin{proof}
		Consider an $\mathbf{(o)Top}^{\op}$-gluing data functor $\mathbf{G}$. We aim to prove that the collection $\left(\mathrm{I}, \subsm{(\Goi{\mbf{G}}{i})}{i\!\in\! \mathrm{I}},\subsm{(\Goij{\mbf{G}}{i}{j}, \Gnij{\mbf{G}}{\etaij{i}{j}}^{\op})}{(i,j)\!\in\! \mathrm{I}^2},\subsm{(\Gnij{\mbf{G}}{\tauij{i}{j}}^{\op},\Gnij{\mbf{G}}{\tauijk{k}{i}{j}{k}}^{\op})}{( i,j,k)\! \in\! \mathrm{I}^3}\right)$ constitutes a topological space (open) gluing data. According to the definition of a functor from $\glI{I}$ to $\mathbf{(o)Top}^{\op}$, we know that $\subsm{(\Goi{\mbf{G}}{i})}{i\!\in\! \mathrm{I}}$ represents a family of topological spaces, $\subsm{(\Goij{\mbf{G}}{i}{j})}{(i,j)\!\in\! \mathrm{I}^2}$ corresponds to a family of topological spaces such that $\Gnij{\mbf{G}}{\etaij{i}{j}}^{\op}$ is a (open) continuous map, and $\subsm{(\Gnij{\mbf{G}}{\tauij{i}{j}}^{\op},\Gnij{\mbf{G}}{\tauijk{k}{i}{j}{k}}^{\op})}{( i,j,k)\! \in\! \mathrm{I}^3}$ denotes a family of (open) continuous maps where  $\Gnij{\mbf{G}}{\tauij{i}{j}}^{\op}$ is a map from $\Goij{\mbf{G}}{i}{j}$ to $\Goij{\mbf{G}}{j}{i}$ and $\Gnij{\mbf{G}}{\tauijk{k}{i}{j}{k}}^{\op}$ is map from $\Goijk{\mbf{G}}{i}{j}{k}$ to $\Goijk{\mbf{G}}{j}{i}{k}$ for all $i,j,k\in \mathrm{ I}$. Now, we proceed to verify that conditions $(a)$, $(b)$, $(c)$ and $(d)$ of Definition \ref{tpsgd} are fulfilled.
		
		$(a)$ can be deduced when applying the gluing data functor $\mathbf{G}$ to Remark \ref{reglue} $(1)$ $(a)$. Since a functor maps identities to identities, condition $(b)$ is automatically satisfied. To establish condition $(c)$, we apply the gluing data functor $\mathbf{G}$ to Remark \ref{reglue} $(1)$ $(c)$. Finally, condition $(d)$ is obtained by applying $\mathbf{G}$ to Remark \ref{reglue} $(1)$ $(e)$.
		
		Conversely, let $\left(\mathrm{I},\subsm{({\subsm{U}{i}})}{i\!\in\! \mathrm{I}},\subsm{(\subsm{U}{i,j}, \upsilon_{i,j} )}{(i,j)\! \in\! \mathrm{I}^2}, \subsm{(\subsm{\varphi}{i,j}, \phijk{k}{i}{j})}{(i,j,k)\!\in\! \mathrm{I}^3}\right)$ be a topological space (open) gluing data. We want to define $\mathbf{G}$ as an $\mathbf{(o)Top}^{\op}$-gluing data functor such that $\Goi{\mbf{G}}{i}:= \subsm{U}{i}$, $\Goij{\mbf{G}}{i}{j}:=\subsm{U}{i,j}$,  $\Gnij{\mbf{G}}{\etaij{i}{j}}^{\op}:=\upsilon_{i,j}, \Gnij{\mbf{G}}{\tauij{i}{j}}^{\op}:=\subsm{\varphi}{i,j}$ and $\Gnij{\mbf{G}}{\tauijk{k}{i}{j}{k}}^{\op}:=\phijk{k}{i}{j}$ for all $i, j,k\in \mathrm{I}$. To prove that $\mathbf{G}$ is well-defined, it suffices to show that $\mathbf{G}$ preserves the equalities $\subsm{\tau}{[i,i]}=\subsm{\operatorname{id}}{[i,i]}$,  $\tauijk{k}{i}{j}{k} \circ \tauijk{i}{j}{k}{i} = \tauijk{j}{i}{k}{j}$, and $\tauijk{k}{i}{j}{k} \circ \etaijk{i}{j}{i}{k}= \etaijk{j}{i}{j}{k} \circ \tauij{i}{j}$ for all $i, j,k\in \mathrm{I}$. This follows directly from Definition \ref{tpsgd} conditions $(b)$, $(c)$ and $(d)$, respectively. 
		
	\end{proof}
\end{lemma}

\begin{figure}[H]
$$ {\tiny\xymatrixrowsep{0.2in}
\xymatrixcolsep{0.1in} \xymatrix{ 
&&&&& U_i \ar@{..>}[dddddd]|-{\dindsi{\iota}{Q}{i}} &&&&& \\
&&&&&&&&&& \\
&&&U_{i,j}  \ar[rruu]|-{\upsilon_{i,j}} \ar@/^1pc/[ddddll]|-{\subsm{\varphi}{i,j}}  & &\circlearrowleft&   & \ar[lluu]|-{\upsilon_{i,k}}U_{i,k}  \ar@/^1pc/[ddddrr]|-{\subsm{\varphi}{i,k}} &&& \\
&&&&&&&&&& \\
&&&\circlearrowleft&U_{i,j,k} \ar[luu]|-{{\iuv{\subsm{U}{i,j,k}}{\subsm{U}{i,j}}}} \ar@/_1pc/[ldd]|-{\phijk{k}{i}{j}}   \ar@{=}[rr] & & U_{i,k,j}  \ar[ruu]|-{{\iuv{\subsm{U}{i,k,j}}{\subsm{U}{i,k}}}} \ar@/_1pc/[rdd]|-{\phijk{j}{i}{k}}  & \circlearrowleft & &&\\  
&& && &  \circlearrowleft&&&& \\
  & U_{j,i} \ar[dddl]|-{\upsilon_{j,i}} \ar@/^1pc/[uuuurr]|-{\subsm{\varphi}{j,i}} &&\ar@/_1pc/[ruu]|-{\phijk{k}{j}{i}} U_{j,i,k}  \ar[ll]|-{\scalebox{0.6}{${\iuv{\subsm{U}{j,i,k}}{\subsm{U}{j,i}}}$} }  & &Q   & & U_{k,i,j} \ar[rr]|- {\scalebox{0.6}{${\iuv{\subsm{U}{k,i,j}}{\subsm{U}{k,i}}}$} }\ar@/_1pc/[luu]|-{\phijk{j}{k}{i}}  
&&  U_{k,i}  \ar[dddr]|-{\upsilon_{k,i}} \ar@/^1pc/[uuuull]|-{\subsm{\varphi}{k,i}} &  \\
&&&&&&&&&& \\
&&\circlearrowleft&& U_{j,k,i}\ar@{=}[uul]  \ar[ldd]|-{\iuv{\subsm{U}{j,k,i}}{\subsm{U}{j,k}}} \ar@/_1pc/[rr]|-{\phijk{i}{j}{k}} && \ar@/_1pc/[ll]|-{\phijk{i}{k}{j}} U_{k,j,i} \ar[rdd]|-{\iuv{\subsm{U}{k,j,i}}{\subsm{U}{k,j}}} \ar@{=}[uur] &&\circlearrowleft& &\\
\subsm{U}{j}\ar@{..>}[uuurrrrr]|-{\dindsi{\iota}{Q}{j}}  &&&&&\circlearrowleft&&&&&\subsm{U}{k}\ar@{..>}[uuulllll]|-{\dindsi{\iota}{Q}{k}}  \\
&&&U_{j,k} \ar[lllu]|-{\upsilon_{j,k}} \ar@/^1pc/[rrrr]|-{{\begin{array}{@{}c@{}} 
\subsm{\varphi}{j,k} \end{array}}}& &&&\ar@/^1pc/[llll]|-{\subsm{\varphi}{k,j}} U_{k,j} \ar[rrru]|-{\upsilon_{k,j}}  &&&\\
}}$$
\caption{Diagram representation of a gluing data functor $\mathbf{G}$ of type $\mathrm{I}=\{ i,j,k\}$ and a glued object $(Q, \iota)$ over $\mathbf{G}$, where, for all $n\in \{j,k\}$: $\Goi{\mathbf{G}}{i} := \subsm{U}{i}$; $\Goij{\mathbf{G}}{i}{j} := \subsm{U}{i,j}$; $\Goijk{\mathbf{G}}{i}{j}{k} := \subsm{U}{i,j}\subsm{\times}{\subsm{U}{i}} \subsm{U}{i,k}$; $\Gnij{\mbf{G}}{\tauij{i}{j}}^{\op}:=\subsm{\varphi}{i,j}$; $\Gnij{\mbf{G}}{\tauijk{k}{i}{j}{k}}^{\op}:=\phijk{k}{i}{j}$; $\Gnij{\mbf{G}}{\etaij{i}{j}}^{\op}:=\upsilon_{i,j}$; $\Gnij{\mbf{G}}{\etaijk{n}{i}{j}{k}}^{\op}:=\iuv{\subsm{U}{i,j,k}}{\subsm{U}{i,n}}$ where $\iuv{\subsm{U}{i,j,k}}{\subsm{U}{i,n}}$ is the canonical pullback morphism.\\ }
\end{figure}

\subsection{Characterization of glued-up \texorpdfstring{${\protect\mathbf{(o)Top}}^{\protect\op}$} \texorpdfstring{-}objects}

This definition provides the concrete description that can be applied in topological settings. By constructing what we call 'the standard representative of the limit of $\mathbf{G}$,' we create a practical and complete framework for understanding and working with these limits. This representative allow us to manipulate and analyze limits effectively. It serves as a tool for translating theoretical concepts into practical applications in the category of topological spaces.

\begin{definition}[Lemma]\label{deflema}
Let $\mathbf{G}$ be a $\mathbf{(o)Top}^{\op}$-gluing data functor.
We define the \textbf{\textit{ standard representative of the limit of $\mathbf{G}$}} as the pair $(\subsm{Q}{\mathbf{G}}, \dindi{\iota}{Q}{\mathbf{G}}^{\op})$ where         
\begin{itemize}
    \item $\subsm{Q}{\mathbf{G}}:={\subsm{\coprod\nolimits}{i\!\in\! \mathrm{I}} \Goi{\mbf{G}}{i}}/\Rel{\mbf{G}}$ such that $\Rel{\mbf{G}}$ is the equivalence relation on the disjoint union $\subsm{\coprod\nolimits}{i\!\in\! \mathrm{I}} \Goi{\mbf{G}}{i}$ defined by $(x,i)\Rel{\mbf{G}} (y, j)$ if there exists $u\in \Goij{\mbf{G}}{i}{j}$ such that 
    $$x= \Gnij{\mbf{G}}{\etaij{i}{j}}^{\op} (u) \text{ and } y= \Gnij{\mbf{G}}{\etaij{j}{i}}^{\op} \circ \Gnij{\mbf{G}}{\tauij{i}{j}}^{\op}(u),$$ for any $(x,i), (y, j) \in \subsm{\coprod\nolimits}{i\!\in\! \mathrm{I}}\Goi{\mbf{G}}{i}$ where $(i,j)\in \mathrm{I}$. Moreover, $\subsm{Q}{\mathbf{G}}$ is a topological space via the final topology with respect to the family $\subsm{({\dindiv{\iota}{Q}{\mbf{G}}{i}}\!)}{i\! \in\! \mathrm{ I}}$. 
    \item $\dindi{\iota}{Q}{\mathbf{G}}^{\op}=\subsm{(\dindiv{\iota}{Q}{\mbf{G}}{a}: \Goi{\mbf{G}}{a}\rightarrow {\subsm{Q}{\mathbf{G}}})}{a\!\in\! \Co{\glI{I}}}$, with ${\dindiv{\iota}{Q}{\mbf{G}}{i}}:= \pi \circ \subsm{\boldsymbol{\varepsilon}}{\Goi{\mbf{G}}{i},\subsm{\coprod\nolimits}{i\!\in\! \mathrm{I}} \Goi{\mbf{G}}{i}}$, ${\dindgij{\iota}{Q}{\mathbf{G}}{i}{j}} :={\dindiv{\iota}{Q}{\mbf{G}}{i}} \circ \Gnij{\mbf{G}}{\etaij{i}{j}}^{\op} $, and ${\dindgijk{\iota}{Q}{\mathbf{G}}{i}{j}{k}}:=\dindiv{\iota}{Q}{\mbf{G}}{i} \circ \Gnij{\mbf{G}}{\etaiijk{i}{j}{k}}^{\op} $ where $\subsm{\boldsymbol{\varepsilon}}{\Goi{\mbf{G}}{i},\subsm{\coprod\nolimits}{j\!\in\! \mathrm{I}} \Goi{\mbf{G}}{j}}$ is the canonical map from $ \Goi{\mbf{G}}{i}$ to $\subsm{\coprod\nolimits}{j\!\in\! \mathrm{I}} \Goi{\mbf{G}}{j}$ sending $x$ to $(x,i)$, and $\pi: \subsm{\coprod\nolimits}{i\!\in\! \mathrm{I}} \Goi{\mbf{G}}{i} \rightarrow {\subsm{Q}{\mathbf{G}}}$ is the quotient map, for all $(i,j,k)\in \mathrm{I}$.  
\end{itemize} 
For all $i \in \mathrm{ I}$, ${\dindiv{\iota}{Q}{\mbf{G}}{i}}\!\!: \Goi{\mbf{G}}{i}\rightarrow {\subsm{Q}{\mathbf{G}}}$ is a one-to-one (open) continuous map. We will prove in Theorem \ref{gluingtop} that  $(\subsm{Q}{\mathbf{G}}, \dindi{\iota}{Q}{\mathbf{G}}^{\op})$ is indeed a limit over $\mathbf{G}$.
\end{definition}
\begin{proof} We first verify that $\Rel{\mbf{G}}$ is an equivalence relation. 
\begin{itemize}
\item We have that $\Rel{\mbf{G}}$ is reflexive since $\Goij{\mbf{G}}{i}{i}=\Goi{\mbf{G}}{i}$ and $\Gnij{\mbf{G}}{\etaij{i}{i}}^{\op} =\Gnij{\mbf{G}}{\tauij{i}{i}}^{\op}=\subsm{\operatorname{id}}{i}$, for $i\in \mathrm{I}$. 

\item We prove that $\Rel{\mbf{G}}$ is symmetric. Let $(x,i), (y, j) \in \subsm{\coprod\nolimits}{i\!\in\! \mathrm{I}}\Goi{\mbf{G}}{i}$ where $i, j\in \mathrm{I}$. Suppose that $(x,i)\Rel{\mbf{G}} (y,j)$. By definition of $\Rel{\mbf{G}}$, there exists $u\in \Goij{\mbf{G}}{i}{j}$ such that 
$$
x= \Gnij{\mbf{G}}{\etaij{i}{j}}^{\op} (u) \quad \text{and} \quad y= \Gnij{\mbf{G}}{\etaij{j}{i}}^{\op} \circ \Gnij{\mbf{G}}{\tauij{i}{j}}^{\op}(u).
$$
Taking $u':= \Gnij{\mbf{G}}{\tauij{i}{j}}^{\op}(u)$, we have $u'\in \Goij{\mbf{G}}{j}{i}$ and 
$$
y= \Gnij{\mbf{G}}{\etaij{j}{i}}^{\op} (u') \quad \text{and} \quad x= \Gnij{\mbf{G}}{\etaij{i}{j}}^{\op} \circ \Gnij{\mbf{G}}{\tauij{j}{i}}^{\op}(u').
$$

\item  We prove that $\Rel{\mbf{G}}$ is transitive. Let $(x,i),(y,j), (z,k) \in \subsm{\coprod\nolimits}{i\!\in\! \mathrm{I}}\Goi{\mbf{G}}{i}$  where $i, j,k \in \mathrm{I}$. Suppose that $(x,i)\Rel{\mbf{G}} (y,j)$  and $(y,j)\Rel{\mbf{G}} (z,k)$. By definition of $\Rel{\mbf{G}}$, there exist $u\in \Goij{\mbf{G}}{i}{j}, v\in  \Goij{\mbf{G}}{j}{k}$ such that 

\begin{itemize} 
\item $x= \Gnij{\mbf{G}}{\etaij{i}{j}}^{\op} (u)$ and $y= \Gnij{\mbf{G}}{\etaij{j}{i}}^{\op} \circ \Gnij{\mbf{G}}{\tauij{i}{j}}^{\op}(u)$, and
\item  $y= \Gnij{\mbf{G}}{\etaij{j}{k}}^{\op} (v)$ and $z= \Gnij{\mbf{G}}{\etaij{k}{j}}^{\op} \circ \Gnij{\mbf{G}}{\tauij{j}{k}}^{\op}(v)$. 
\end{itemize} 

Using the uniqueness of pullbacks up to cone isomorphism, we know that there is a cone isomorphism 
\begin{itemize} 
\item $\psi_{j,i,k}$ from the pullback 
$$\Goij{\mbf{G}}{j}{i}\subsm{\times}{\Goi{\mbf{G}}{j}}  \Goij{\mbf{G}}{j}{k}=\{ (u, v)
\in \Goij{\mbf{G}}{j}{i}\times  \Goij{\mbf{G}}{j}{k}  | \Gnij{\mbf{G}}{\etaij{j}{i}}^{\op} (u)= \Gnij{\mbf{G}}{\etaij{j}{k}}^{\op} (v)\}$$ 
to $\Goijk{\mbf{G}}{j}{i}{k}=\Goijk{\mbf{G}}{j}{k}{i}$. By the definition of a cone morphism, we obtain $\subsm{\mathfrak{i}}{\Goijk{\mbf{G}}{j}{i}{k}, \Goij{\mbf{G}}{j}{i}}\circ \psi_{j,i,k}= \subsm{\mathfrak{i}}{\Goij{\mbf{G}}{j}{i}\subsm{\times}{\Goi{\mbf{G}}{j}}  \Goij{\mbf{G}}{j}{k}, \Goij{\mbf{G}}{j}{i}}$.
\item $\psi_{i,j,k}$ from the pullback 
$$\Goij{\mbf{G}}{i}{j}\subsm{\times}{\Goi{\mbf{G}}{i}}  \Goij{\mbf{G}}{i}{k}=\{ (u, v)
\in \Goij{\mbf{G}}{i}{j}\times  \Goij{\mbf{G}}{i}{k}  | \Gnij{\mbf{G}}{\etaij{i}{j}}^{\op} (u)= \Gnij{\mbf{G}}{\etaij{i}{k}}^{\op} (v)\}$$
 to $\Goijk{\mbf{G}}{i}{j}{k}=\Goijk{\mbf{G}}{i}{k}{j}$.  By the definition of a cone morphism, we obtain $\subsm{\mathfrak{i}}{\Goijk{\mbf{G}}{i}{j}{k}, \Goij{\mbf{G}}{i}{j}}\circ \psi_{i,j,k}= \subsm{\mathfrak{i}}{\Goij{\mbf{G}}{i}{j}\subsm{\times}{\Goi{\mbf{G}}{i}}  \Goij{\mbf{G}}{i}{k}, \Goij{\mbf{G}}{i}{j}}$.

\end{itemize}
which is possible since by definition $\mbf{G}$ sends pushouts to pushouts, and pushouts in the opposite category are pullbacks. 

By definition of the pullback $ \Goij{\mbf{G}}{j}{i}\subsm{\times}{\Goi{\mbf{G}}{j}}  \Goij{\mbf{G}}{j}{k}$, we have 
$$(\Gnij{\mbf{G}}{\tauij{i}{j}}^{\op}(u),v)\in  \Goij{\mbf{G}}{j}{i}\subsm{\times}{\Goi{\mbf{G}}{j}}  \Goij{\mbf{G}}{j}{k} \text{ and } \psi_{j,i,k} (\Gnij{\mbf{G}}{\tauij{i}{j}}^{\op}(u),v)\in \Goijk{\mbf{G}}{j}{i}{k}.$$

By our assumptions, we know that $y= \Gnij{\mbf{G}}{\etaij{j}{i}}^{\op} \circ \Gnij{\mbf{G}}{\tauij{i}{j}}^{\op}(u)=\Gnij{\mbf{G}}{\etaij{j}{k}}^{\op} (v)$. We set $$\alpha:= \psi_{i,j,k}^{-1} \circ \Gnij{\mbf{G}}{\tauijk{k}{j}{i}{k}}^{\op}\circ  \psi_{j,i,k} (\Gnij{\mbf{G}}{\tauij{i}{j}}^{\op}(u),v),$$ and we know that $\alpha \in \Goij{\mbf{G}}{i}{j}\subsm{\times}{\Goi{\mbf{G}}{i}}  \Goij{\mbf{G}}{i}{k}$. Moreover,

$$
\begin{array}{lll} 
&&\subsm{\mathfrak{i}}{\Goij{\mbf{G}}{i}{j}\subsm{\times}{\Goi{\mbf{G}}{i}}  \Goij{\mbf{G}}{i}{k}, \Goij{\mbf{G}}{i}{j}}(\alpha) \\
& = & \subsm{\mathfrak{i}}{\Goij{\mbf{G}}{i}{j}\subsm{\times}{\Goi{\mbf{G}}{i}}  \Goij{\mbf{G}}{i}{k}, \Goij{\mbf{G}}{i}{j}}\circ \psi_{i,j,k}^{-1}  \circ \Gnij{\mbf{G}}{\tauijk{k}{j}{i}{k}}^{\op}\circ \psi_{j,i,k} (\Gnij{\mbf{G}}{\tauij{i}{j}}^{\op}(u),v)\\
& = & \subsm{\mathfrak{i}}{\Goijk{\mbf{G}}{i}{j}{k}, \Goij{\mbf{G}}{i}{j}}  \circ \Gnij{\mbf{G}}{\tauijk{k}{j}{i}{k}}^{\op}\circ \psi_{j,i,k} (\Gnij{\mbf{G}}{\tauij{i}{j}}^{\op}(u),v)\\
& = & \Gnij{\mbf{G}}{\tauij{j}{i}}^{\op}\circ \subsm{\mathfrak{i}}{\Goijk{\mbf{G}}{j}{i}{k}, \Goij{\mbf{G}}{j}{i}}  (\Gnij{\mbf{G}}{\tauij{i}{j}}^{\op}(u),v)\\
& = & \Gnij{\mbf{G}}{\tauij{j}{i}}^{\op}\circ \subsm{\mathfrak{i}}{\Goij{\mbf{G}}{j}{i}\subsm{\times}{\Goi{\mbf{G}}{j}}  \Goij{\mbf{G}}{j}{k}, \Goij{\mbf{G}}{j}{i}}  (\Gnij{\mbf{G}}{\tauij{i}{j}}^{\op}(u),v)=u,
\end{array}
$$
by Remark \ref{reglue} $(1) (e)$.
Thus $\alpha=(u,\subsm{\mathfrak{i}}{\Goij{\mbf{G}}{i}{j}\subsm{\times}{\Goi{\mbf{G}}{i}}  \Goij{\mbf{G}}{i}{k}, \Goij{\mbf{G}}{i}{k}}(\alpha))$. 

Now, we set $w:=\subsm{\mathfrak{i}}{\Goij{\mbf{G}}{i}{j}\subsm{\times}{\Goi{\mbf{G}}{i}}  \Goij{\mbf{G}}{i}{k}, \Goij{\mbf{G}}{i}{k}}(\alpha)$, and we know that $w\in \Goij{\mbf{G}}{i}{k}$. By definition of the pullback $\Goij{\mbf{G}}{i}{j}\subsm{\times}{\Goi{\mbf{G}}{i}}  \Goij{\mbf{G}}{i}{k}$, we have $x=\Gnij{\mbf{G}}{\etaij{i}{j}}^{\op} (u)=\Gnij{\mbf{G}}{\etaij{i}{k}}^{\op} (w)$. 

Now we want to prove that $z=\Gnij{\mbf{G}}{\etaij{k}{i}}^{\op}\circ \Gnij{\mbf{G}}{\tauij{i}{k}}^{\op}(w)$. We have:
\begin{align*}
&\Gnij{\mbf{G}}{\etaij{k}{i}}^{\op}\circ \Gnij{\mbf{G}}{\tauij{i}{k}}^{\op}(w)\\
=&\Gnij{\mbf{G}}{\etaij{k}{i}}^{\op}\circ \Gnij{\mbf{G}}{\tauij{i}{k}}^{\op}\circ \subsm{\mathfrak{i}}{\Goijk{\mbf{G}}{i}{j}{k}, \Goij{\mbf{G}}{i}{k}}\circ \Gnij{\mbf{G}}{\tauijk{k}{j}{i}{k}}^{\op}\circ \psi_{j,i,k}(\Gnij{\mbf{G}}{\tauij{i}{j}}^{\op}(u),v)\\
=&\Gnij{\mbf{G}}{\etaij{k}{i}}^{\op}\circ \subsm{\mathfrak{i}}{\Goijk{\mbf{G}}{k}{j}{i}, \Goij{\mbf{G}}{k}{i}}\circ \Gnij{\mbf{G}}{\tauijk{j}{i}{k}{j}}^{\op} \circ  \Gnij{\mbf{G}}{\tauijk{k}{j}{i}{k}}^{\op}\circ \psi_{j,i,k}(\Gnij{\mbf{G}}{\tauij{i}{j}}^{\op}(u),v)\\
& \;\text{by Remark \ref{reglue} $(1) (e)$}\\
=&\Gnij{\mbf{G}}{\etaij{k}{i}}^{\op}\circ \subsm{\mathfrak{i}}{\Goijk{\mbf{G}}{k}{j}{i}, \Goij{\mbf{G}}{k}{i}}\circ \Gnij{\mbf{G}}{\tauijk{i}{j}{k}{i}}^{\op}\circ \psi_{j,i,k} (\Gnij{\mbf{G}}{\tauij{i}{j}}^{\op}(u),v)\\
&\;\text{by Remark \ref{reglue} $(1) (c)$}\\
=&\Gnij{\mbf{G}}{\etaij{k}{j}}^{\op}\circ \Gnij{\mbf{G}}{\tauij{j}{k}}^{\op}\circ  \subsm{\mathfrak{i}}{\Goijk{\mbf{G}}{j}{i}{k}, \Goij{\mbf{G}}{j}{k}} \circ \psi_{j,i,k} (\Gnij{\mbf{G}}{\tauij{i}{j}}^{\op}(u),v)\\
&\;\text{by Remark \ref{reglue} $(1) (e)$}\\
=&\Gnij{\mbf{G}}{\etaij{k}{j}}^{\op}\circ \Gnij{\mbf{G}}{\tauij{j}{k}}^{\op}\circ  \subsm{\mathfrak{i}}{\Goij{\mbf{G}}{j}{i}\subsm{\times}{\Goi{\mbf{G}}{j}}  \Goij{\mbf{G}}{j}{k}, \Goij{\mbf{G}}{j}{k}}  (\Gnij{\mbf{G}}{\tauij{i}{j}}^{\op}(u),v)\\
=&\Gnij{\mbf{G}}{\etaij{k}{j}}^{\op} \circ \Gnij{\mbf{G}}{\tauij{j}{k}}^{\op}(v)=z.
\end{align*}
\end{itemize}
\noindent This completes proving that $\Rel{\mbf{G}}$ is an equivalence relation.

Finally, we prove that the map ${\dindiv{\iota}{Q}{\mbf{G}}{i}}\!\!: \Goi{\mbf{G}}{i}\rightarrow {\subsm{Q}{\mathbf{G}}}$ is a one-to-one continuous map. Since $\subsm{Q}{\mathbf{G}}$ is a topological space under the final topology with respect to the family $\subsm{({\dindiv{\iota}{Q}{\mbf{G}}{i}}\!)}{i\!\in\! \mathrm{I}}$, it follows that ${\dindiv{\iota}{Q}{\mbf{G}}{i}}$ is a continuous map, for all $i \in \mathrm{I}$. So, we only need to prove that ${\dindiv{\iota}{Q}{\mbf{G}}{i}}$ is one-to-one and open, for all $i \in \mathrm{I}$.

To establish the one-to-one property, consider $x,y\in \Goi{\mbf{G}}{i}$ such that ${\dindiv{\iota}{Q}{\mbf{G}}{i}}\!\!(x)={\dindiv{\iota}{Q}{\mbf{G}}{i}}\!\!(y)$. According to the definition of ${\dindiv{\iota}{Q}{\mbf{G}}{i}}$, we have $\pi(x,i)=\pi(y,i)$. Hence, by the definition of $\Rel{\mbf{G}}$, it follows that $x=y$, since $\Goij{\mbf{G}}{i}{i}=\Goi{\mbf{G}}{i}$ and $\Gnij{\mbf{G}}{\etaij{i}{i}}^{\op} =\Gnij{\mbf{G}}{\tauij{i}{i}}^{\op}=\subsm{\operatorname{id}}{i}$, for $i\in \mathrm{I}$. 

When $\mathbf{G}$ is an $\mathbf{oTop}^{\op}$-gluing data functor, we prove that $\dindiv{\iota}{Q}{\mbf{G}}{i}$ is also an open map. Consider $V\subsm{\subseteq}{\operatorname{op}} \Goi{\mbf{G}}{i}$. Our goal is to show that ${\dindiv{\iota}{Q}{\mbf{G}}{i}}\!\!(V)\subsm{\subseteq}{\operatorname{op}} \subsm{Q}{\mathbf{G}}$. In other words, we need to prove that ${\dindiv{\iota}{Q}{\mbf{G}}{j}^{-1}}\!({\dindiv{\iota}{Q}{\mbf{G}}{i}}\!\!(V))\subsm{\subseteq}{\operatorname{op}}\Goi{\mbf{G}}{j}$ for all $j\in \mathrm{I}$, considering the definition of the final topology on $\subsm{Q}{\mathbf{G}}$ with respect to $\dindi{\iota}{Q}{\mathbf{G}}$. Let $(i,j) \in \mathrm{I}$. To establish this, we prove that ${\dindiv{\iota}{Q}{\mbf{G}}{j}^{-1}}\!({\dindiv{\iota}{Q}{\mbf{G}}{i}}\!\!(V))= {\Gnij{\mbf{G}}{\etaij{j}{i}}^{\op}}(\Gnij{\mbf{G}}{\tauij{i}{j}}^{\op}( {\Gnij{\mbf{G}}{\etaij{i}{j}}^{\op}}^{-1} (V)))$. We have
	\begin{align*}
		&y\in {\Gnij{\mbf{G}}{\etaij{j}{i}}^{\op}}(\Gnij{\mbf{G}}{\tauij{i}{j}}^{\op}( {\Gnij{\mbf{G}}{\etaij{i}{j}}^{\op}}^{-1} (V)))  \\&\Leftrightarrow y={\Gnij{\mbf{G}}{\etaij{j}{i}}^{\op}}(\Gnij{\mbf{G}}{\tauij{i}{j}}^{\op}( z))\; \text{for some $z \in {\Gnij{\mbf{G}}{\etaij{i}{j}}^{\op}}^{-1} (V)$}  
		\\& \Leftrightarrow  \pi(y,j)=\pi(x,i)\; \text{for some $x\in V$ (take $x={{\Gnij{\mbf{G}}{\etaij{i}{j}}^{\op}}} (z)$)}  \\& \Leftrightarrow y\in {\dindiv{\iota}{Q}{\mbf{G}}{j}^{-1}}\!({\dindiv{\iota}{Q}{\mbf{G}}{i}}\!\!(V)), \; \text{since}\; \pi(x,i)={\dindiv{\iota}{Q}{\mbf{G}}{i}}\!\!(x)\; \\& \quad \quad \text{and}\; \pi(y,j)={\dindiv{\iota}{Q}{\mbf{G}}{j}}\!\!(y)\; \text{for some $x\in V $}.
	\end{align*}
	Since ${\Gnij{\mbf{G}}{\etaij{i}{j}}^{\op}}$ is continuous, ${\Gnij{\mbf{G}}{\etaij{i}{j}}^{\op}}^{-1} (V)$ is open. Given that $\Gnij{\mbf{G}}{\tauij{i}{j}}^{\op}$ is a homeomorphism, it follows that $\Gnij{\mbf{G}}{\tauij{i}{j}}^{\op}({\Gnij{\mbf{G}}{\etaij{i}{j}}^{\op}}^{-1}(V))$ is also open. Since ${\Gnij{\mbf{G}}{\etaij{j}{i}}^{\op}}$ is a morphism in $\Cn{\mathbf{oTop}}$, $ {\Gnij{\mbf{G}}{\etaij{j}{i}}^{\op}}(\Gnij{\mbf{G}}{\tauij{i}{j}}^{\op}( {\Gnij{\mbf{G}}{\etaij{i}{j}}^{\op}}^{-1} (V)))$ is open. Thus, successfully demonstrating that $\dindi{\iota}{Q}{i}$ is a topological embedding when $\mathbf{G}$ is an $\mathbf{oTop}^{\op}$-gluing data functor.	
	
\end{proof}

\begin{remark}\label{remtop}
	Let $\mathbf{G}$ be an $\mathbf{(o)Top}^{\op}$-gluing data functor. 
	
		 For all $q\in \subsm{Q}{\mathbf{G}}$ there exists $i\in \mathrm{ I}$ and $x\in \Goi{\mbf{G}}{i}$ such that $\dindiv{\iota}{Q}{\mbf{G}}{i}\!(x)=q$. Moreover, by definition of $\dindiv{\iota}{Q}{\mbf{G}}{i}$ we obtain $q=\pi(x,i)$.
	
\end{remark}

We introduce the 'Gluing Topological Spaces Theorem'. This theorem establishes the equivalence between a topological space obtained by gluing objects using the specific gluing data functor constructed above and the conventional construction of glued-up spaces. It elucidates conditions that allow us to view a 'glued-up' object as a limit of this functor, catering to a more intuitive understanding of an index gluing category. Furthermore, this theorem determines properties characterizing these composite spaces. As we embark on the journey of proving this theorem, we leave no essential proof component hidden from view.

\begin{theorem} \label{gluingtop}
	Given a $\mathbf{(o)Top}^{\op}$-gluing data functor $\mathbf{G}$. Let $Q$ be a topological space, $\dindi{\iota}{Q}{}^{\op}$ be a family $\subsm{(\dindi{\iota}{Q}{a}^{\op})}{a\! \in\! \Co{\glI{I}}}$ where $\dindi{\iota}{Q}{a}\!:  \Goi{\mbf{G}}{a}\rightarrow Q$ are in $\Cn{\mathbf{(o)Top}}$, for all $a \in \Co{\glI{I}}$. The following assertions are equivalent:
	\begin{enumerate}
		\item  $Q$ is a glued-up $\mathbf{(o)Top}^{\op}$-object along $\mathbf{G}$ through $\dindi{\iota}{Q}{}^{\op}$;
		\item $(Q, \subsm{\iota}{Q}^{\op})$ is a cone over $\mathbf{G}$ isomorphic to $(\subsm{Q}{\mathbf{G}},\dindi{\iota}{Q}{\mathbf{G}}^{\op})$ in the category of cones over $\mathbf{G}$;		
		\item For all $(i,j,k)\in \mathrm{I}$, $(Q, \dindi{\iota}{Q}{}^\op)$ satisfies the following properties:
		\begin{itemize}
			\item[(a)] $\dindij{\iota}{Q}{i}{j} =\dindi{\iota}{Q}{i} \circ  \Gnij{\mbf{G}}{\etaij{i}{j}}^{\op}$;  
			\item[(b)] $\dindijk{\iota}{Q}{i}{j}{k} =\dindi{\iota}{Q}{(i,j)}  \circ \Gnij{\mbf{G}}{\etaijk{n}{i}{j}{k}}^{\op}$;
			\item[(c)]  $\dindi{\iota}{Q}{i} \circ \Gnij{\mbf{G}}{\etaij{i}{j}}^{\op} ={\dindi{\iota}{Q}{j}} \circ \Gnij{\mbf{G}}{\etaji{i}{j}}^{\op}  \circc\Gnij{\mbf{G}}{\tauij{i}{j}}^\op$;
			\item[(d)] $Q=\subsm{\cup}{i\!\in\! \mathrm{I}} {\dindsi{\iota}{Q}{i}}\!\!(\Goi{\mbf{G}}{i})$;
			\item[(e)] $\dindi{\iota}{Q}{i}$ is a one-to-one (open) continuous map.
		\end{itemize}	
	\end{enumerate}

	\begin{proof} 
		We start the proof proving that $(\subsm{Q}{\mathbf{G}},\dindi{\iota}{Q}{\mathbf{G}}^{\op})$ satisfies $(3)$ $(c)$ and $(d)$ as we will need it throughout the proof. Let $(i,j,k)\in \mathrm{I}^3$.
		\begin{itemize} 
		
\item We need to prove that $\dindiv{\iota}{Q}{\mbf{G}}{i} \circ \Gnij{\mbf{G}}{\etaij{i}{j}}^{\op} =\dindiv{\iota}{Q}{\mbf{G}}{j} \circ \Gnij{\mbf{G}}{\etaji{i}{j}}^{\op}  \circc\Gnij{\mbf{G}}{\tauij{i}{j}}$. Let $x \in \Goij{\mbf{G}}{i}{j}$. We have:

\begin{align*}
&\dindiv{\iota}{Q}{\mbf{G}}{i} \circ \Gnij{\mbf{G}}{\etaij{i}{j}}^{\op} (x) = \pi ( \Gnij{\mbf{G}}{\etaij{i}{j}}^{\op} (x), i), \text{ and } \\
&\dindiv{\iota}{Q}{\mbf{G}}{j} \circ \Gnij{\mbf{G}}{\etaji{i}{j}}^{\op}  \circc\Gnij{\mbf{G}}{\tauij{i}{j}} (x) = \pi ( \Gnij{\mbf{G}}{\etaji{i}{j}}^{\op}  \circ\Gnij{\mbf{G}}{\tauij{i}{j}} (x), j)
\end{align*}

Therefore, we obtain the equality wanted by the definition of the equivalence relation $\Rel{\mbf{G}}$. Hence, property $(3)(c)$ is satisfied.
			\item Using the definition of $\subsm{Q}{\mathbf{G}}$, we obtain easily that the property $(3)(d)$ is satisfied.
			
		\end{itemize}
		
Property $(3) (c)$ implies for any $(i,j) \in \mathrm{I}^2$
$$\dindiv{\iota}{Q}{\mbf{G}}{j}\!(\Gnij{\mbf{G}}{\etaji{i}{j}}^{\op} (\Goij{\mbf{G}}{j}{i}))=\dindiv{\iota}{Q}{\mbf{G}}{i}\!(\Gnij{\mbf{G}}{\etaji{j}{i}}^{\op} (\Goij{\mbf{G}}{i}{j}))\!=\dindiv{\iota}{Q}{\mbf{G}}{i}\!(\Goi{\mbf{G}}{i})\cap \dindiv{\iota}{Q}{\mbf{G}}{j}\!(\Goi{\mbf{G}}{j}).$$
The equality $\dindiv{\iota}{Q}{\mbf{G}}{j}\!(\Gnij{\mbf{G}}{\etaji{i}{j}}^{\op} (\Goij{\mbf{G}}{j}{i}))=\dindiv{\iota}{Q}{\mbf{G}}{i}\!(\Gnij{\mbf{G}}{\etaji{j}{i}}^{\op} (\Goij{\mbf{G}}{i}{j}))$ follows easily from property $(3) (c)$ since $\Gnij{\mbf{G}}{\tauij{i}{j}}^\op$ is a homeomorphism from $\Goij{\mbf{G}}{i}{j}$ to $\Goij{\mbf{G}}{j}{i}$.

We now prove that $\dindiv{\iota}{Q}{\mbf{G}}{j}\!(\Gnij{\mbf{G}}{\etaji{i}{j}}^{\op} (\Goij{\mbf{G}}{j}{i}))\!=\dindiv{\iota}{Q}{\mbf{G}}{i}\!(\Goi{\mbf{G}}{i})\cap \dindiv{\iota}{Q}{\mbf{G}}{j}\!(\Goi{\mbf{G}}{j}).$

Let $x\in \dindiv{\iota}{Q}{\mbf{G}}{i}\!(\Gnij{\mbf{G}}{\etaji{j}{i}}^{\op} (\Goij{\mbf{G}}{i}{j})$, we know that there exists $u\in \Goij{\mbf{G}}{i}{j}$ such that $x=\dindiv{\iota}{Q}{\mbf{G}}{i}\!(\Gnij{\mbf{G}}{\etaji{j}{i}}^{\op}(u))$.

We can deduce that $x\in \dindiv{\iota}{Q}{\mbf{G}}{i}(\Goi{\mbf{G}}{i})$ since $\Gnij{\mbf{G}}{\etaji{j}{i}}^{\op}(u)\in \Goi{\mbf{G}}{i}$. 

From the equality $\dindiv{\iota}{Q}{\mbf{G}}{j}\!(\Gnij{\mbf{G}}{\etaji{i}{j}}^{\op} (\Goij{\mbf{G}}{j}{i}))=\dindiv{\iota}{Q}{\mbf{G}}{i}\!(\Gnij{\mbf{G}}{\etaji{j}{i}}^{\op} (\Goij{\mbf{G}}{i}{j}))$, we also have $x\in \dindiv{\iota}{Q}{\mbf{G}}{j} \!(\Gnij{\mbf{G}}{\etaji{i}{j}}^{\op} (\Goij{\mbf{G}}{j}{i}))$ then there exists $v\in \Goij{\mbf{G}}{j}{i}$ such that $x=\dindiv{\iota}{Q}{\mbf{G}}{j}\!(\Gnij{\mbf{G}}{\etaji{i}{j}}^{\op}(v))$

This implies that $x\in \dindiv{\iota}{Q}{\mbf{G}}{j}(\Goi{\mbf{G}}{j})$ since $\Gnij{\mbf{G}}{\etaji{i}{j}}^{\op}(v)\in \Goi{\mbf{G}}{j}$. Therefore, $x\in \dindiv{\iota}{Q}{\mbf{G}}{i}\!(\Goi{\mbf{G}}{i})\cap \dindiv{\iota}{Q}{\mbf{G}}{j}\!(\Goi{\mbf{G}}{j})$. Conversely, suppose that $x\in \dindiv{\iota}{Q}{\mbf{G}}{i}\!(\Goi{\mbf{G}}{i})\cap \dindiv{\iota}{Q}{\mbf{G}}{j}\!(\Goi{\mbf{G}}{j})$. There exists $u\in \Goi{\mbf{G}}{i}$ and $v\in \Goi{\mbf{G}}{j}$ such that $x=\pi(u,i)=\pi(v,j)$. That is, $(u,i)\Rel{\mbf{G}} (v,j)$. By the definition of the relation $\Rel{\mbf{G}}$, we know that there exists $w\in \Goij{\mbf{G}}{i}{j}$ such that $u= \Gnij{\mbf{G}}{\etaij{i}{j}}^{\op} (w)$ and $v= \Gnij{\mbf{G}}{\etaij{j}{i}}^{\op} \circ \Gnij{\mbf{G}}{\tauij{i}{j}}^{\op}(w)$. We conclude that $$x\in \dindiv{\iota}{Q}{\mbf{G}}{i}\!(\Gnij{\mbf{G}}{\etaji{j}{i}}^{\op} (\Goij{\mbf{G}}{i}{j})).$$ 

		\begin{enumerate}
			\item[$(1)\Rightarrow (2)$] 
			Let $i\in\mathrm{I}$. 
			We prove that $(\subsm{Q}{\mathbf{G}},\dindi{\iota}{Q}{\mathbf{G}}^{\op})$ is a terminal cone over $\mathbf{G}$. 
			\begin{itemize} 			
				\item Our aim is to prove that the pair $(\subsm{Q}{\mathbf{G}},\dindi{\iota}{Q}{\mathbf{G}}^{\op})$ satisfies property $(1) (a)$, $(b)$, and $(c)$ of Remark \ref{oppc}. Property $(1) (a)$ holds by definition of $\dindgij{\iota}{Q}{\mbf{G}}{i}{j}$. Property $(1) (b)$ follows directly from the definition of ${\dindgij{\iota}{Q}{\mbf{G}}{i}{j}}$ and $\dindgijk{\iota}{Q}{\mathbf{G}}{i}{j}{k}$. 
				\item Finally suppose that $(Q',{\subsm{\iota}{{Q'}}})$ is another pair  making the following diagram
				\begin{figure}[H] 
					\begin{center}
						{\tiny	\begin{tikzcd}[column sep=normal]
								&\Goij{\mbf{G}}{i}{j}\arrow[]{rrrr}{\Gnij{\mbf{G}}{\etaji{i}{j}}^{\op} \circ \Gnij{\mbf{G}}{\tauij{i}{j}}^{\op}} \arrow[swap]{d}{\Gnij{\mbf{G}}{\etaij{i}{j}}^{\op}}& && &\Goi{\mbf{G}}{j}\arrow[]{d} {{\dindpi{\iota}{{Q}}{\!j}}}\\ 
								& \Goi{\mbf{G}}{i} \arrow[swap]{rrrr}{\dindpi{\iota}{{Q}}{i}}&& & &Q'
					\end{tikzcd}}\end{center} 	 \caption{}\label{topol130}	
				\end{figure}
				\noindent commute, for all $(i,j)\in \mathrm{I}^{2}$. We want to prove that there exists a unique map $\mu:\subsm{Q}{\mathbf{G}}\rightarrow Q'$ in $\Cn{\mathbf{\topo}}$ making the following diagram
				\begin{figure}[H]
					\begin{center}
						{\tiny  \begin{tikzcd}[column sep=normal]
								\Goij{\mbf{G}}{i}{j}\arrow[swap]{d}{\Gnij{\mbf{G}}{\etaij{i}{j}}^{\op}}  \arrow[]{rrrr}{{\Gnij{\mbf{G}}{\etaji{i}{j}}^{\op}\circ \Gnij{\mbf{G}}{\tauij{i}{j}}^{\op}}}  && & &
								\Goi{\mbf{G}}{j}\arrow[bend left=20]{ddr}{{\dindpi{\iota}{{Q}}{j}}} \arrow[]{d}{{\dindiv{\iota}{Q}{\mbf{G}}{j}}}\\ 
								\Goi{\mbf{G}}{i}\arrow[swap,bend right=15]{drrrrr}{\dindpi{\iota}{{Q}}{i}}    \arrow[swap]{rrrr}{{\dindiv{\iota}{Q}{\mbf{G}}{i}}}  &&&&
								\subsm{Q}{\mathbf{G}}\arrow[dashed]{dr}{ \mu} \\ 
								&&&& &Q'
						\end{tikzcd} }
					\end{center}\caption{}\label{topol132}
				\end{figure}
				\noindent commutes, for all $(i,j)\in \mathrm{I}^{ 2}$. 
				\begin{itemize} 
					\item If such a $\mu$ exists we have $\mu\circ {\dindiv{\iota}{Q}{\mbf{G}}{i}}=\dindpi{\iota}{{Q}}{i}$, for all $i\in \mathrm{I}$. Let $q\in {\subsm{Q}{\mathbf{G}}}$. By Remark \ref{remtop}, there is $y\in \subsm{\coprod\nolimits}{i\!\in\! \mathrm{I}} \Goi{\mbf{G}}{i}$, such that $q= \pi (y)$. 
					Thus, we have $q= \dindiv{\iota}{Q}{\mbf{G}}{i}\!\!(x)$. Therefore, if $\mu$ exists, it is uniquely determined by $\mu(q)=\dindpi{\iota}{{Q}}{i}(x)$. 
					\item We will demonstrate the well-definedness of $\mu$. Let $(i,j)\in \mathrm{I}^{2}$, $x \in \Goi{\mbf{G}}{i}$, and $y \in \Goi{\mbf{G}}{j}$ such that $\dindiv{\iota}{Q}{\mbf{G}}{i}\!(x) = \dindiv{\iota}{Q}{\mbf{G}}{j}\!(y)$. We prove that $\dindpi{\iota}{{Q}}{i}(x)= \dindpi{\iota}{{Q}}{i}(y)$.
					
	By definition of $\dindiv{\iota}{Q}{\mbf{G}}{i}$ and $\dindiv{\iota}{Q}{\mbf{G}}{j}$ implies $\pi(x,i)=\pi(y,j)$. This implies that $(x,i)\Rel{\mbf{G}} (y,j)$. By definition of the relation $\Rel{\mbf{G}}$ there exists $u\in \Goij{\mbf{G}}{i}{j}$ such that \begin{equation}\label{rg1}x= \Gnij{\mbf{G}}{\etaij{i}{j}}^{\op} (u)\end{equation} and \begin{equation}\label{rg2}y= \Gnij{\mbf{G}}{\etaij{j}{i}}^{\op} \circ \Gnij{\mbf{G}}{\tauij{i}{j}}^{\op}(u).\end{equation} Applying $\dindpi{\iota}{Q}{i}$ to Equation \ref{rg1} and $\dindpi{\iota}{Q}{j}$ to Equation \ref{rg2} we obtain $$\dindpi{\iota}{Q}{i}(x)=\dindpi{\iota}{Q}{i}(\Gnij{\mbf{G}}{\etaij{i}{j}}^{\op} (u))$$ and $$\dindpi{\iota}{Q}{j}(y)=\dindpi{\iota}{Q}{j}(\Gnij{\mbf{G}}{\etaij{j}{i}}^{\op} \circ \Gnij{\mbf{G}}{\tauij{i}{j}}^{\op}(u)).$$ By commutativity of the diagram in Figure \ref{topol130} we get $\dindpi{\iota}{Q}{i}(x)=\dindpi{\iota}{Q}{j}(y)$.
					
					Hence, we can conclude that the map $\mu$ is indeed well-defined.
					
					\item We now prove that $\mu$ is continuous. Let $U\subsm{\subseteq}{\operatorname{op}} {Q'}$, we prove that \\$\mu^{-1}(U)\subsm{\subseteq}{\operatorname{op}} \subsm{Q}{\mathbf{G}}$. By definition of the final topology with respect to ${\dindiv{\iota}{Q}{\mbf{G}}{i}}$, we need to prove that $\dindiv{\iota}{Q}{\mbf{G}}{i}^{-1}\!(\mu^{-1}(U))\subsm{\subseteq}{\operatorname{op}} \Goi{\mbf{G}}{i}$ for all $i \in \mathrm{I}$. Let $i\in \mathrm{I}$. We know that $\dindiv{\iota}{Q}{\mbf{G}}{i}^{-1}\!(\mu^{-1}(U))=(\mu\circ {\dindiv{\iota}{Q}{\mbf{G}}{i}})^{-1}(U)= \dindpi{\iota}{Q}{i}^{-1} (U)$, by the definition of $\mu$. Since $\dindpi{\iota}{{Q}}{i}$ is a morphism in $\Cn{\topo}$, $\dindpi{\iota}{{Q}}{i}^{-1}(U)\subsm{\subseteq}{\operatorname{op}} \Goi{\mbf{G}}{i}$. Hence this proves that $\mu^{-1}(U)\subsm{\subseteq}{\operatorname{op}} \subsm{Q}{\mathbf{G}}$. Therefore $\mu$ is continuous.

				\end{itemize}
			\end{itemize}
			
			This completes proving that $\subsm{Q}{\mathbf{G}}$ is a glued-up $\topo^{\op}$-object along $\mathbf{G}$ through $\dindi{\iota}{Q}{\mathbf{G}}^{\op}$. Finally, we prove (2) using the uniqueness of limits up to isomorphism. 
			\item[$(2) \Rightarrow (3)$] Suppose $(Q, \subsm{\iota}{Q}^{\op})$ is a cone over $\mathbf{G}$ isomorphic to $(\subsm{Q}{\mathbf{G}}, \dindi{\iota}{Q}{\mathbf{G}}^{\op})$ in the category of cones over $\mathbf{G}$. Since $(Q, \subsm{\iota}{Q}^{\op})\simeq (\subsm{Q}{\mathbf{G}}, \dindi{\iota}{Q}{\mathbf{G}}^{\op})$, we obtain that $(Q, \subsm{\iota}{Q}^{\op})$ satisfies $(3)$ $(a)$-$(f)$ is equivalent to proving that $(\subsm{Q}{\mathbf{G}}, \dindi{\iota}{Q}{\mathbf{G}}^{\op})$ satisfies $(3)$ $(a)$-$(f)$.  Since $ (\subsm{Q}{\mathbf{G}}, \dindi{\iota}{Q}{\mathbf{G}}^{\op})$ is a cone over $\mathbf{G}$, we have that properties $(3)$ $(a)$-$(c)$ are satisfied. Finally, properties $(3)$ $(d)$-$(e)$ have been proven above and this concludes the proof.
			\item[$(3)\Rightarrow (1)$] Suppose that $(Q, \dindi{\iota}{Q}{}^{\op})$ satisfies the properties $(a)$-$(e)$ of the statement $(3)$. We want to prove that $Q$ is glued-up $\topo^{\op}$-space along $\mathbf{G}$ through $\dindi{\iota}{Q}{}^{\op}$. Property $(1) (a)$ and $(1) (b)$ of the Remark \ref{oppc} are precisely the  properties $(a)$ and $(b)$ of statement $(3)$ of this theorem. 
			The following diagram
			\begin{figure}[H]\begin{center}
					{\tiny\begin{tikzcd}[column sep=normal]
							&\Goij{\mbf{G}}{i}{j}\arrow[]{rrrr}{{\Gnij{\mbf{G}}{\etaji{i}{j}}^{\op}\circ \Gnij{\mbf{G}}{\tauij{i}{j}}^{\op}}} \arrow[swap]{d}{\Gnij{\mbf{G}}{\etaij{i}{j}}^{\op}} & &&& \Goi{\mbf{G}}{j}\arrow[]{d} {{\dindi{\iota}{Q}{j}}}\\
							& \Goi{\mbf{G}}{i} \arrow[swap]{rrrr}{\dindi{\iota}{Q}{i}}&& & & Q
				\end{tikzcd}}\end{center}  \caption{}\label{topol122}   \end{figure}
			\noindent 
			commutes, for all $(i,j)\in \mathrm{I}^{2}$, by property $(3) (c)$ of our assumptions.
			Next, let $(Q',{\subsm{\iota}{Q'}})$ be another pair as above, making the following diagram 
			\begin{figure}[H]\begin{center}
					{\tiny	\begin{tikzcd}[column sep=normal]
							&\Goij{\mbf{G}}{i}{j}\arrow[]{rrrr}{{\Gnij{\mbf{G}}{\etaji{i}{j}}^{\op}\circ \Gnij{\mbf{G}}{\tauij{i}{j}}^{\op}}} \arrow[swap]{d}{\Gnij{\mbf{G}}{\etaij{i}{j}}^{\op}} & &&& \Goi{\mbf{G}}{j}\arrow[]{d} {{\dindpi{\iota}{Q}{\!j}}}\\
							& \Goi{\mbf{G}}{i} \arrow[swap]{rrrr}{\dindpi{\iota}{Q}{\!i}}&& &&Q'
				\end{tikzcd}}\end{center}  \caption{}\label{topol123}   \end{figure}
			\noindent commute. We want to prove that there exists a unique map $\mu: Q\rightarrow Q'$ in $\Cn{\mathbf{\topo}}$ making the following diagram 
			\begin{figure}[H]
				\begin{center}
					{\tiny  \begin{tikzcd}[column sep=normal]
							\Goij{\mbf{G}}{i}{j} \arrow[]{rrrr}{{\Gnij{\mbf{G}}{\etaji{i}{j}}^{\op}\circ \Gnij{\mbf{G}}{\tauij{i}{j}}^{\op}}}\arrow[swap]{d}{\Gnij{\mbf{G}}{\etaij{i}{j}}^{\op}} && & &
							\Goi{\mbf{G}}{j}  \arrow[]{d}{{\dindi{\iota}{Q}{j}}} \arrow[bend left=20]{ddr}{{\dindpi{\iota}{{Q}}{j}}}\\ 
							\Goi{\mbf{G}}{i}\arrow[swap]{rrrr}{\dindi{\iota}{Q}{i}} \arrow[swap,bend right=15]{drrrrr}{\dindpi{\iota}{{Q}}{i}}    & && &
							Q\arrow[dashed,near start]{dr}{\exists!\mu} \\
							& & && &Q' 
					\end{tikzcd} }
				\end{center}\caption{}\label{topol14}
			\end{figure}
			\noindent commute, for all $(i,j)\in \mathrm{I}^{2}$. If such $\mu$ exists the commutativity of the diagram gives $\mu\circ \dindi{\iota}{Q}{i}\!(x)=\dindpi{\iota}{Q}{i}(x)$ for all $i\in \mathrm{I}$ and $x\in  \Goi{\mbf{G}}{i}$. Thus, $\mu$ is uniquely determined by $\dindpi{\iota}{Q}{i}$. Indeed, by property $(3) (d)$ of the assumption, given $q \in {Q}$, by Remark \ref{remtop}, there exists $i\in \mathrm{I}$ and $x\in \Goi{\mbf{G}}{i}$ such that $q= \dindi{\iota}{Q}{i}\!(x)$ and $\mu(q)=\dindpi{\iota}{Q}{i}(x)$. 
			
			We will demonstrate the well-definedness of $\mu$. Let $(i,j)\in \mathrm{I}^{2}$, $x \in \Goi{\mbf{G}}{i}$, and $y \in \Goi{\mbf{G}}{j}$ such that $\dindiv{\iota}{Q}{\mbf{G}}{i}\!(x) = \dindiv{\iota}{Q}{\mbf{G}}{j}\!(y)$. We prove that $\dindpi{\iota}{{Q}}{i}(x)= \dindpi{\iota}{{Q}}{i}(y)$. Since we have proven above that $\dindi{\iota}{Q}{i}\!(\Gnij{\mbf{G}}{\etaji{j}{i}}^{\op} (\Goij{\mbf{G}}{i}{j}))\!=\dindi{\iota}{Q}{i}\!(\Goi{\mbf{G}}{i})\cap {\dindi{\iota}{Q}{j}}\!(\Goi{\mbf{G}}{j})$ and using property $(3) (c)$, we can find $u\in \Goij{\mbf{G}}{i}{j}$ such that
			$$\dindi{\iota}{{Q}}{i}(x)= \dindi{\iota}{{Q}}{i}(y)= \dindi{\iota}{Q}{i}\!(\Gnij{\mbf{G}}{\etaji{j}{i}}^{\op} (u)) ={\dindi{\iota}{Q}{j}} ( \Gnij{\mbf{G}}{\etaji{i}{j}}^{\op}  (\Gnij{\mbf{G}}{\tauij{i}{j}}^\op(u))).$$
			Since $\dindi{\iota}{{Q}}{i}$ and $\dindi{\iota}{{Q}}{j}$ are one-to-one, by Property $(3) (f)$, we obtain $x= \Gnij{\mbf{G}}{\etaji{j}{i}}^{\op} (u)$ and $y=\Gnij{\mbf{G}}{\etaji{i}{j}}^{\op}  (\Gnij{\mbf{G}}{\tauij{i}{j}}^\op(u))$. Applying $\dindpi{\iota}{Q}{i}$ to the equality and using the commutativity of the diagram above, we obtain  $\dindpi{\iota}{{Q}}{i}(x)= \dindpi{\iota}{{Q}}{i}(y)$. This completes the proof.

When $\mathbf{G}$ is chosen to be a $\mathbf{oTop}^\op$-gluing data functor. We need to prove that $\mu$ is also an open map. 

Let $V \subop \subsm{Q}{\mathbf{G}}$ be given. We have $V = \subsm{\cup}{i \in \mathrm{I}}(V \cap {\dindiv{\iota}{Q}{\mbf{G}}{i}}\!\!(\Goi{\mbf{G}}{i}))$, since $\subsm{Q}{\mbf{G}} = \subsm{\cup}{i \in \mathrm{I}}{\dindiv{\iota}{Q}{\mbf{G}}{i}}\!\!(\Goi{\mbf{G}}{i})$, and ${\dindiv{\iota}{Q}{\mbf{G}}{i}}\!\!(\Goi{\mbf{G}}{i})$ is open in $\subsm{Q}{\mathbf{G}}$ for all $i \in \mathrm{I}$, as ${\dindiv{\iota}{Q}{\mbf{G}}{i}}$ is an open map.
We write $\mu(V) = \subsm{\cup}{i \in \mathrm{I}} \mu(V \cap {\dindiv{\iota}{Q}{\mbf{G}}{i}}\!\!(\Goi{\mbf{G}}{i})) \subseteq Q'$. Let $i \in \mathrm{I}$. We have $\mu(V \cap {\dindiv{\iota}{Q}{\mbf{G}}{i}}\!\!(\Goi{\mbf{G}}{i})) = \dindpi{\iota}{{Q}}{i}({\dindiv{\iota}{Q}{\mbf{G}}{i}^{-1}}\!(V \cap {\dindiv{\iota}{Q}{\mbf{G}}{i}}\!\!(\Goi{\mbf{G}}{i})))$ by commutativity of the diagram in Figure \ref{topol14} and $\dindiv{\iota}{Q}{\mbf{G}}{i}$ is a surjective map onto $V \cap \dindiv{\iota}{Q}{\mbf{G}}{i}\!\!(\Goi{\mbf{G}}{i})$.

Now, $V \cap {\dindiv{\iota}{Q}{\mbf{G}}{i}}\!\!(\Goi{\mbf{G}}{i})$ is open in $\subsm{Q}{\mathbf{G}}$ as an intersection of two opens. Moreover, ${\dindiv{\iota}{Q}{\mbf{G}}{i}^{-1}}\!(V \cap {\dindiv{\iota}{Q}{\mbf{G}}{i}}\!\!(\Goi{\mbf{G}}{i}))$ is open as ${\dindiv{\iota}{Q}{\mbf{G}}{i}}$ is continuous. Finally, $\mu(V \cap {\dindiv{\iota}{Q}{\mbf{G}}{i}}\!\!(\Goi{\mbf{G}}{i}))$ is open, since $\dindpi{\iota}{{Q}}{i}$ is an open map as it is in $\Cn{\mathbf{oTop}}$.
		 
		\end{enumerate}
	\end{proof}
\end{theorem}

\begin{remark}
\begin{enumerate}
\item The gluing topological space data corresponding to a $\mathbf{(o)Top}^{\op}$-gluing data functor is more general than the typical gluing topological space data found in the litterature, as discussed in \cite[Proposition $12.27$]{wedhorn2016}.
\item In the proof of the Theorem, we have proven that Property $(3) (c)$ implies that for any $(i,j) \in \mathrm{I}^{2}$, 
$$\dindi{\iota}{Q}{j}\!(\Gnij{\mbf{G}}{\etaji{i}{j}}^{\op} (\Goij{\mbf{G}}{j}{i}))=\dindi{\iota}{Q}{i}\!(\Gnij{\mbf{G}}{\etaji{j}{i}}^{\op} (\Goij{\mbf{G}}{i}{j}))\!=\dindi{\iota}{Q}{i}\!(\Goi{\mbf{G}}{i})\cap {\dindi{\iota}{Q}{j}}\!(\Goi{\mbf{G}}{j}).$$
\item Given a $\mathbf{(o)Top}^{\op}$-gluing data functor $\mathbf{G}$ and $Q$ a topological space, let $\dindi{\iota}{Q}{}^{\op}$ be a family $\subsm{\{\dindi{\iota}{Q}{a}^{\op} \}}{a\!\in\! \Co{\glI{I}}}$ where $\dindi{\iota}{Q}{a}: \Goi{\mbf{G}}{a}\rightarrow Q$ is in $\Cn{\mathbf{(o)Top}}$ for all $a\in \Co{\glI{I}}$. If $Q$ is a glued-up $\mathbf{(o)Top}^{\op}$-object along $\mathbf{G}$ through $\dindi{\iota}{Q}{}^{\op}$, then the topology on $Q$ is the final topology with respect to $\dindi{\iota}{Q}{}^{\op}$. This follows directly from $(Q, \subsm{\iota}{Q}^{\op})\simeq (\subsm{Q}{\mathbf{G}}, \dindi{\iota}{Q}{\mathbf{G}}^{\op})$ proved in Theorem \ref{gluingtop}.

\item With the notation of Theorem \ref{gluingtop}, let $i,j ,k  \in \mathrm{I}$, we have 
\begin{align*}&{\dindiv{\iota}{Q}{\mbf{G}}{i}}\!(\Gnij{\mbf{G}}{\etaij{i}{j}}^{\op}\circ \Gnij{\mbf{G}}{\etaijk{j}{i}{j}{k}}^{\op}(\Goijk{\mbf{G}}{i}{j}{k} ))\\&= {\dindiv{\iota}{Q}{\mbf{G}}{j}}\!(\Gnij{\mbf{G}}{\etaij{j}{i}}^{\op}\circ  \Gnij{\mbf{G}}{\etaijk{i}{j}{i}{k}}^{\op}(\Goijk{\mbf{G}}{j}{i}{k} )) 
\\&= {\dindiv{\iota}{Q}{\mbf{G}}{i}}\!(\Gnij{\mbf{G}}{\etaij{i}{j}}^{\op}(\Goij{\mbf{G}}{i}{j}))\cap {\dindiv{\iota}{Q}{\mbf{G}}{i}}\!(\Gnij{\mbf{G}}{\etaij{i}{k}}^{\op}(\Goij{\mbf{G}}{i}{k}))\\&= {\dindiv{\iota}{Q}{\mbf{G}}{j}}\!(\Gnij{\mbf{G}}{\etaij{j}{i}}^{\op}(\Goij{\mbf{G}}{j}{i}))\cap {\dindiv{\iota}{Q}{\mbf{G}}{j}}\!(\Gnij{\mbf{G}}{\etaij{j}{k}}^{\op}(\Goij{\mbf{G}}{j}{k}))
\\&=  {\dindiv{\iota}{Q}{\mbf{G}}{i}}\!( \Goi{\mbf{G}}{i})\cap {\dindiv{\iota}{Q}{\mbf{G}}{j}}\!( \Goi{\mbf{G}}{j})\cap {\dindiv{\iota}{Q}{\mbf{G}}{k}}\!( \Goi{\mbf{G}}{k})
;\end{align*}  
Since $  \Gnij{\mbf{G}}{\tauijk{k}{i}{j}{k}}^{\op} (\Goijk{\mbf{G}}{i}{j}{k} )= \Goijk{\mbf{G}}{j}{i}{k}$, then \begin{align*} &{\dindiv{\iota}{Q}{\mbf{G}}{i}}\!(\Gnij{\mbf{G}}{\etaij{i}{j}}^{\op}\circ \Gnij{\mbf{G}}{\etaijk{j}{i}{j}{k}}^{\op}(\Goijk{\mbf{G}}{i}{j}{k}))\\&={\dindiv{\iota}{Q}{\mbf{G}}{j}}\!\!\circ \Gnij{\mbf{G}}{\etaij{j}{i}}^{\op} \Gnij{\mbf{G}}{\etaijk{i}{j}{i}{k}}^{\op} \circ \Gnij{\mbf{G}}{\tauijk{k}{i}{j}{k}}^{\op} (\Goijk{\mbf{G}}{i}{j}{k})\\&={\dindiv{\iota}{Q}{\mbf{G}}{j}}\!(\Gnij{\mbf{G}}{\etaij{j}{i}}^{\op}\circ  \Gnij{\mbf{G}}{\etaijk{i}{j}{i}{k}}^{\op}(\Goijk{\mbf{G}}{j}{i}{k}).\end{align*} Moreover, we can prove as in the proof of the Theorem \ref{gluingtop} that
\begin{align*}&{\dindiv{\iota}{Q}{\mbf{G}}{i}}\!(\Gnij{\mbf{G}}{\etaij{i}{j}}^{\op}\circ \Gnij{\mbf{G}}{\etaijk{j}{i}{j}{k}}^{\op}(\Goijk{\mbf{G}}{i}{j}{k} ))\\&={\dindiv{\iota}{Q}{\mbf{G}}{i}}\!(\Gnij{\mbf{G}}{\etaij{i}{j}}^{\op}(\Goij{\mbf{G}}{i}{j}))\cap {\dindiv{\iota}{Q}{\mbf{G}}{i}}\!(\Gnij{\mbf{G}}{\etaij{i}{k}}^{\op}(\Goij{\mbf{G}}{i}{k})).\end{align*}  and that 
\begin{align*}&{\dindiv{\iota}{Q}{\mbf{G}}{j}}\!(\Gnij{\mbf{G}}{\etaij{j}{i}}^{\op}\circ  \Gnij{\mbf{G}}{\etaijk{i}{j}{i}{k}}^{\op}(\Goijk{\mbf{G}}{j}{i}{k} ))
\\&={\dindiv{\iota}{Q}{\mbf{G}}{j}}\!(\Gnij{\mbf{G}}{\etaij{j}{i}}^{\op}(\Goij{\mbf{G}}{j}{i})) \cap {\dindiv{\iota}{Q}{\mbf{G}}{j}}(\Gnij{\mbf{G}}{\etaij{j}{k}}^{\op}(\Goij{\mbf{G}}{j}{k})).\end{align*}

\end{enumerate} 
\end{remark}

The following Definition (Lemma) follows directly from Theorem \ref{gluingtop}.
\begin{definition}[Lemma] \label{gluingcover}
Let $U$ be a topological space, $\rm{I}$ be an index set and $\mathcal{U}:=\subsm{\{\subsm{U}{i}\}}{i\!\in\! \rm{I}}$ be an open covering of $U$. We define $\subsm{\mathfrak{i}}{\mathcal{U}}^{\op}:=\subsm{(\iuv{\Coi{\gcov{\subsm{\mathcal{U}}{0}}}{a}}{U})}{a\!\in\!\Co{\glI{I}}}$ and $\gcov{\mathcal{U}}: \glI{I} \rightarrow \mathbf{oTop}^\op$ where for all $i,j, k\in \rm{I}$ and  $n\in \{ j,k\}$ 
\begin{enumerate}
    \item $\Coi{\gcov{\subsm{\mathcal{U}}{0}}}{i}=\subsm{U}{i}$;
    \item $\Coij{\gcov{\subsm{\mathcal{U}}{0}}}{i}{j}=\subsm{U}{i,j}$;
    \item $\Coijk{\gcov{\subsm{\mathcal{U}}{0}}}{i}{j}{k}=\subsm{U}{i,j,k}$; 
    \item $\Cnij{\gcov{\subsm{\mathcal{U}}{1}}}{\etaij{i}{j}}=\iuv{\subsm{U}{i,j}}{\subsm{U}{i}}^{\op}$;
\item $\Cnij{\gcov{\subsm{\mathcal{U}}{1}}}{\tauij{i}{j}}$ is the canonical pullback isomorphism from $\subsm{U}{i,j}$ to $\subsm{U}{j,i}$;
\item 
$\Cnij{\gcov{\subsm{\mathcal{U}}{1}}}{\etaijk{n}{i}{j}{k}}=\iuv{\subsm{ U}{i,j,k}}{\subsm{U}{i,n}}^{\op}$
\end{enumerate}
where  $\subsm{U}{i,j}= U_i \cap U_j$ and $\subsm{U}{i,j,k}:= U_i \cap U_j\cap U_k$, for all $i,j,k \in \mathrm{I}$.
 Then, $\gcov{\mathcal{U}}$ is an $\mathbf{oTop}^\op$-gluing data functor and $U$ is a glued-up $\mathbf{oTop}^\op$-object along ${\gcov{\mathcal{U}}}$ through $\subsm{\mathfrak{i}}{\mathcal{U}}^{\op}$. 
\end{definition}

\section{Gluing sheaves categorically}
For some background materials about sheaf theory, we refer to \cite{Liu}.

\subsubsection{\texorpdfstring{$\mathbf{E}\Psh{\protect\Ops{X}}{Ab}$} \texorpdfstring{and} \texorpdfstring{$\mathbf{E}\Sh{\protect\Ops{X}}{Ab}$}\texorpdfstring{-}gluing data functor} 
We begin this section by defining a category of enriched presheaves of abelian group on $\Ops{X}$, $\mathbf{E}\Psh{\Ops{X}}{Ab}$ (resp. enriched sheaves of abelian group on $\Ops{X}$, $\mathbf{E}\Sh{\Ops{X}}{Ab}$). Then, we define a $\mathbf{E}\Psh{\Ops{X}}{Ab}$ (resp. $\mathbf{E}\Sh{\Ops{X}}{Ab}$)-gluing data functor. Later, we prove that giving such a functor is equivalent to giving a  presheaf (resp. sheaf) gluing data. In this context, we have chosen to operate within an enriched gluing category, which serves as a precursor to the category of ringed topological spaces.

\begin{definition}\label{catshf}
    Let $X$ be a topological space. The category of enriched presheaves (resp. sheaves) on $\Ops{X}$ denoted as $\mathbf{E}\Psh{\Ops{X}}{Ab}$ (resp. $\mathbf{E}\Sh{\Ops{X}}{Ab}$) is the category in which: 
    \begin{enumerate}
        \item the objects are of the form $(U, \mathcal{F})$ where $U\in \dindsi{\mbf{Ops}}{X}{0}$ and $\mathcal{F}$ is a presheaf (resp. enriched sheaf) on $U$;
        \item the morphisms are pairs $(\iuv{V}{U}^{\op}, \alpha): (U, \mathcal{F})\rightarrow (V, \mathcal{G})$ where $\alpha$ is defined as a family $\subsm{(\subsm{\alpha}{W}: \sho{\mathcal{F}}{W} \rightarrow \sho{\mathcal{G}}{W\cap V})}{W\!\in\!  \dindi{\mbf{Ops}}{U}{0}}$ of morphism of abelian groups such that for any inclusion map $\iuv{W'}{W}$ in $\dindi{\mbf{Ops}}{U}{1}$, the following diagram
        \begin{figure}[H]
\begin{center}
		{\tiny	\begin{tikzcd}[column sep=normal]
				&\sho{\mathcal{F}}{W} \arrow{rr}{\Cn{\mathcal{F}}(\iuv{W'}{W})} \arrow[swap]{d}{ \subsm{\alpha}{W}} && \sho{\mathcal{F}}{W'}\arrow{d} {\subsm{\alpha}{W'}}\\ 
				& \sho{\mathcal{G}}{W\cap V} \arrow[swap]{rr}{\:\Cn{\mathcal{G}}(\iuv{W'\cap V}{W\cap V})}& &\sho{\mathcal{G}}{W'\cap V}
			\end{tikzcd}}
		\end{center} \caption{}\label{topol138}  \end{figure}
    \noindent commutes. We also refer to $\alpha$ as a natural transformation from $\mathcal{F}$ to $\mathcal{G}$.
    \end{enumerate}
\end{definition}
For the rest of section $4$, $\mathbf{C}$ refers to an element in the set $\{\mathbf{E}\Psh{\Ops{X}}{Ab},\mathbf{E}\Sh{\Ops{X}}{Ab}\}$.
\begin{definition}\label{sheffnct}
    Let $U$ be a topological space, $\rm{I}$ be an index set and $\mathcal{U}:=\subsm{\{\subsm{U}{i}\}}{i\!\in \!\rm{I}}$ be an open covering of $U$. A $\mathbf{C}$\textbf{\textit{-gluing data functor}} $\mathbf{G}$ along $\mathcal{U}$ is a functor from $\glI{I}$ to $\mathbf{C}$ such that for all $i,j, k \in\rm{I}$ and $n\in\{ j,k \}$ we have 
    \begin{enumerate}
        \item $\Goi{\mbf{G}}{i}=(\subsm{U}{i}, \subsm{\mathcal{F}}{i})$;
        \item $\Goij{\mbf{G}}{i}{j}=\subsm{\Goi{\mbf{G}}{i}|}{\subsm{U}{i,j}}$;
        \item $\Goijk{\mbf{G}}{i}{j}{k} =\subsm{\Goi{\mbf{G}}{i}|}{\subsm{U}{i,j,k}}$;
        \item $\Gnij{\bf{G}}{\etaij{i}{j}}=\subsm{\subsm{\underbar{$\mathfrak{i}$}}{\Goi{\mbf{G}}{i}}}{\subsm{U}{i,j},\subsm{U}{i}}$; 
        \item $\Gnij{\bf{G}}{\etaijk{n}{i}{j}{k}}=\subsm{\subsm{\underbar{$\mathfrak{i}$}}{\Goi{\mbf{G}}{i}}}{\subsm{U}{i,j,k}, \subsm{U}{i,n}}$.      
    \end{enumerate}
    
\end{definition}
The following definition is the usual definition for a presheaf (resp. sheaf) gluing data on a topological space $X$.
\begin{definition}\label{sheafgludata}
Let $X$ be a topological space, $\rm{I}$ be an index set and  $\mathcal{U}:=\subsm{\{\subsm{U}{i}\}}{i\!\in\! \rm{I}}$ be a open covering of $X$.  A collection $$\big(\subsm{(\subsm{\mathcal{F}}{i})}{i\!\in \!\rm{I}}, \subsm{(\subsm{\Phi}{i,j})}{(i,j)\!\in\! \mathrm{I}^2}\big)$$ 
where 
\begin{enumerate}

\item $\subsm{\mathcal{F}}{i}$ is a presheaf (resp. sheaf) of abelian groups on $\subsm{U}{i}$ for each $i\in \rm{I}$;
\item $\subsm{\Phi}{i,j}$ is an isomorphism of presheaf (resp. sheaf) of abelian groups from $ \subsm{\subsm{\mathcal{F}}{i}|}{\subsm{U}{i,j}}$ to $ \subsm{\subsm{\mathcal{F}}{j}|}{\subsm{U}{i,j}}$ for each $i,j\in \mathrm{I}$;
 \end{enumerate}
such that 
\begin{enumerate}
\item[a)] for each $i\in \mathrm{I}$, $\subsm{\Phi}{i,i}=\subsm{\operatorname{id}}{\subsm{\mathcal{F}}{i}}$;  
\item [b)] for each $i,j,k\in \mathrm{I}$, $\subsm{\Phi}{i,k}=\subsm{\Phi}{j,k}\circc \subsm{\Phi}{i,j}$ on $\Uij{i,j}{i,k}$;
 \end{enumerate}
is called a \textbf{\textit{presheaf}} (resp. \textbf{\textit{sheaf}}) \textbf{\textit{gluing data}} on $X$ with respect to the open covering $\mathcal{U}$.
\end{definition}

\begin{remark}\label{shinv}
Given $\big(\subsm{(\subsm{\mathcal{F}}{i})}{i\!\in \!\rm{I}}, \subsm{\subsm{(\Phi}{i,j})}{(i,j)\!\in\! \mathrm{I}^2}\big)$ a presheaf (resp. sheaf) gluing data, $\subsm{\Phi}{i,j}$ is a natural correspondence whose inverse is $\subsm{\Phi}{j,i}$. This is a consequence of $b)$ applied to $k=i$ and $a)$.
\end{remark} 

In the following lemma, we establish the equivalence between being given a $\mathbf{C}$-gluing data functor and being provided with a presheaf (or sheaf) gluing data. 
\begin{lemma}\label{equishf}
    Let $U$ be a topological space and $\mathcal{U}:=\subsm{\{\subsm{U}{i}\}}{i\!\in\! \rm{I}}$ be an open covering of $U$. A $\mathbf{C}$-gluing data functor $\mathbf{G}$ along $\mathcal{U}$ induces the presheaf (resp. sheaf) gluing data $\big(\subsm{(\subsm{\Goi{\mbf{G}}{i}}{\mathbf{Sh}})}{i\!\in\! \rm{I}}, \subsm{(\subsm{\Gnij{\bf{G}}{\tauij{j}{i}}}{\mathbf{Sh}})}{( i,j)\! \in\! \rm{I}^2}\big)$ with respect to $\mathcal{U}$. 
   On the other hand, a presheaf (resp. sheaf) gluing data $\big(\subsm{(\subsm{\mathcal{F}}{i})}{i\!\in \!\mathrm{I}}, \subsm{\subsm{(\Phi}{i,j})}{(i,j)\!\in\! \mathrm{I}^2}\big)$ induces the 
    $\mathbf{C}$-gluing data functor $\mathbf{G}$ defined by \\$\Goi{\mbf{G}}{i}:=(\subsm{U}{i}, \subsm{\mathcal{F}}{i})$ and $\Gnij{\bf{G}}{\tauij{j}{i}}:=(\subsm{\operatorname{id}}{\Uij{i}{j}}^{\op} ,\subsm{\Phi}{i,j})$ for all $i,j\in \rm{I}$.
    \begin{proof}
      Let $\mathbf{G}$ be a $\mathbf{C}$-gluing data functor on $U$ along $\mathcal{U}$. We want to prove that $$\big(\subsm{(\subsm{\Goi{\mbf{G}}{i}}{\mathbf{Sh}})}{i\!\in\! \rm{I}}, \subsm{(\subsm{\Gnij{\bf{G}}{\tauij{j}{i}}}{\mathbf{Sh}})}{( i,j)\! \in\! \rm{I}^2}\big)$$ is a presheaf (resp. sheaf) gluing data on $U$ along $\mathcal{U}$. We prove that conditions $a)$ and $b)$ of Definition \ref{sheafgludata} are satisfied. Since any functor preserves identities, we have that condition $a)$ is satisfied. For condition $b)$, since any functor preserves compositions, then applying $\mathbf{G}$ to the equality $\tauijk{k}{i}{j}{k} \circc \tauijk{i}{j}{k}{i} =\tauijk{j}{i}{k}{j}$ in $\Cn{\glI{I}}$  we obtain $\subsm{\Phi}{i,k}=\subsm{\Phi}{j,k}\circc \subsm{\Phi}{i,j}$ on $\Uij{i,j}{i,k}$ for all $i,j,k\in \rm{I}$. Conversely, let $\big(\subsm{(\subsm{\mathcal{F}}{i})}{i\!\in \!\rm{I}}, \subsm{(\subsm{\Phi}{i,j})}{(i,j)\!\in\! \mathrm{I}^2}\big)$ be a presheaf (resp. sheaf) gluing data. We define $\mathbf{G}$ to be a functor from $\glI{I}$ to $\mathbf{C}$ by setting $\Goi{\mbf{G}}{i}:=(\subsm{U}{i}, \subsm{\mathcal{F}}{i})$ and $\Gnij{\bf{G}}{\tauij{j}{i}}:=(\subsm{\operatorname{id}}{U_{i,j}}^{\op} ,\subsm{\Phi}{i,j})$ for all $i,j\in\rm{I}$, $\mathbf{G}$ is well-defined. Indeed, it is enough to prove that $\mathbf{G}$ preserves the equalities $\tauij{i}{i}=\subsm{\operatorname{id}}{[i,i]}$, and $ \tauijk{k}{i}{j}{k}  \circc   \tauijk{i}{j}{k}{i} = \tauijk{j}{i}{k}{j}$, and this follows from Definition \ref{sheafgludata} $a)$ and $b)$, respectively.
    \end{proof}
\end{lemma}

\subsubsection{Characterization of glued-up \texorpdfstring{$\mathbf{E}\Psh{\protect\Ops{X}}{Ab}$} \texorpdfstring{and} \texorpdfstring{$\mathbf{E}\Sh{\protect\Ops{X}}{Ab}$}
	\texorpdfstring{-}objects}\text{         }
We provide a definition for a representative of the limit over a $\mathbf{C}$-gluing data functor. The proof in Theorem \ref{rivegluedshf} utilizes the categorical framework discussed in section 2. We would like to acknowledge David Smith for sharing his comprehensive write-up of the classical proof for the gluing of sheaves from \cite[Chapter II, Exercise 1.22]{harts}. Certain aspects of the classical proof remain pertinent and have been integrated into the proof presented in Theorem \ref{rivegluedshf}.

\begin{definition}\label{stlemma}
Let $U$ be a topological space and $\mathcal{U}:=\subsm{\{\subsm{U}{i}\}}{i\!\in\! \rm{I}}$ be a open covering of $U$. Let $\mathbf{G}$ be a $\mathbf{C}$-gluing data functor.
We define \textbf{\textit{the standard representative of the limit of $\mathbf{G}$}} as the pair $((U ,\subsm{\mathcal{L}}{\mathbf{G}}) , {(\subsm{\mathfrak{i}}{\mathcal{U}}^{\op},\dindi{\pi}{\mathcal{L}}{\mathbf{G}}}))$ where
\begin{itemize}
    \item $\subsm{\mathcal{L}}{\mathbf{G}}$ is the presheaf (resp. sheaf) on $U$ defined by 
  $$
	\dindsi{\mathcal{L}}{\mathbf{G}}{0}\!(V):=\begin{Bmatrix}\subsm{(\subsm{s}{i})}{i\!\in \!\rm{I}} \in \subsm{\prod\nolimits}{i\! \in\! \rm{I} }{\dindsi{\Goi{\mbf{G}}{i}}{\mathbf{Sh}}{0}}\!\big(V\cap U_i\big) |  \dindsi{\Gnij{\bf{G}}{\tauij{j}{i}}}{\mathbf{Sh}}{V}\!\!\big(\subsm{\subsm{s}{i}|}{V \cap \subsm{U}{i,j} }\big)=\subsm{\subsm{s}{j}|}{V\! \cap \subsm{U}{i,j}}, \forall i,j\in\rm{I}\end{Bmatrix}$$
 
 and $\dindsi{\mathcal{L}}{\mathbf{G}}{1}\!(\iuv{W}{V}^{\op}):= \subsm{(\dindsi{\Goi{\mbf{G}}{i}}{\mathbf{Sh}}{1}\!\!\big(\iuvop{W\cap \subsm{U}{i}}{V\cap \subsm{U}{i}}{\op}\big))}{i\!\in\! \rm{I}} $ for all $W\subsm{\subseteq}{\operatorname{op}} V \subsm{\subseteq}{\operatorname{op}} U$ ; 
    \item $(\subsm{\mathfrak{i}}{\mathcal{U}}^{\op}, \dindi{\pi}{\mathcal{L}}{\mathbf{G}})\!:= \subsm{(\iuv{\Coi{\gcov{\subsm{\mathcal{U}}{0}}}{a}}{U}, \dindiv{\pi}{\mathcal{L}}{\mathbf{G}}{a})}{a\! \in\! \Co{\glI{I}}}$ such that $\subsm{\mathfrak{i}}{\mathcal{U}}^{\op}$ is defined in Definition \ref{gluingcover} and given $V\subsm{\subseteq}{\operatorname{op}} U$, for each $i\in \rm{I}$, 
    $\dindgiv{\pi}{\mathcal{L}}{\mathbf{G}}{i}{V}\!\!:  \dindsi{\mathcal{L}}{\mathbf{G}}{0}\!(V) \rightarrow \dindsi{\Goi{\mbf{G}}{i}}{\mathbf{Sh}}{0}\!\big(V\cap U_i\big) $ sends $ \subsm{(\subsm{s}{k})}{k\!\in \!\rm{I}}$ to $\subsm{s}{i},$
      for each $(i,j,k) \in \rm{I}^3$ and $n\in \{j,k\}$,  
   $\dindgijv{\pi}{\mathcal{L}}{\mbf{G}}{i}{j}{V}\!\!:=\dindsi{\Gnij{\bf{G}}{\etaij{i}{j}}}{\mathbf{Sh}}{V}\circc \dindgiv{\pi}{\mathcal{L}}{\mbf{G}}{i}{V}$ and $\dindgijkv{\pi}{\mathcal{L}}{\mbf{G}}{i}{j}{k}{V}\!\!:=\dindsi{\Gnij{\bf{G}}{\etaijk{n}{i}{j}{k}}}{\mathbf{Sh}}{V}\circc \dindgijv{\pi}{\mathcal{L}}{\mbf{G}}{i}{n}{V}.$ 
   
\end{itemize}
\end{definition}

\begin{remark}
   Let $i \in \rm{I}$. We can define $\dindiv{\rho}{\mathcal{L}}{\mathbf{G}}{i}\!\!: \dindsi{\Goi{\mbf{G}}{i}}{\mathbf{Sh}}{}\rightarrow \subsm{\mathcal{L}}{\mbf{G}}$ such that $\dindiv{\rho}{\mathcal{L}}{i}{V}\!\!: \dindsi{\Goi{\mbf{G}}{i}}{\mathbf{Sh}}{0}\!(V\cap\subsm{U}{i})\rightarrow \dindsi{\mathcal{L}}{\mbf{G}}{0}\!(V)$ sending $\subsm{s}{i}$ to $\subsm{(\subsm{\subsm{s}{i}|}{V\cap \subsm{U}{j}})}{j\!\in \! \rm{I}}$ for all $V\subsm{\subseteq}{\operatorname{op}} X$. Moreover, $\dindiv{\pi}{\mathcal{L}}{\mathbf{G}}{i}\circc \dindiv{\rho}{\mathcal{L}}{\mathbf{G}}{i}=\subsm{\operatorname{id}}{\dindsi{\Goi{\mbf{G}}{i}}{\mathbf{Sh}}{}}$
\end{remark}

The following theorem describes a glued-up object up to isomorphism. Since $U$ is a glued-up $\mathbf{(o)Top}^{\op}$-object along $\gcov{\mathcal{U}}$ through $\subsm{\mathfrak{i}}{\mathcal{U}}^{\op}$ by Definition \ref{gluingcover} (Lemma) and an isomorphism in $\mathbf{(o)Top}^{\op}$ is an equality, the first component of the limit over $\mathbf{G}$ is $U$.
\begin{theorem}\label{rivegluedshf}
Let $U$ be a topological space and $\mathcal{U}=\subsm{\{\subsm{U}{i}\}}{i\!\in\! \rm{I}}$ be a open covering of $U$. Let $\mathbf{G}$ be a $\mathbf{C}$-gluing data functor, $\mathcal{F}$ be a sheaf on $U$, $\subsm{\pi}{\mathcal{F}}$ be a family $\subsm{(\dindi{\pi}{\mathcal{F}}{a})}{a\!\in\! \Co{\glI{I}}}$ of  morphism of abelian groups $\dindi{\pi}{\mathcal{F}}{a}: \mathcal{F}\rightarrow \Goi{\mbf{G}}{a}$ for all $a\in \Co{\glI{I}}$.   The following assertions are equivalent:
\begin{enumerate} 
\item $(U,\mathcal{F})$ is a glued-up $\mathbf{C}$-object along $\mathbf{G}$ through $(\subsm{\mathfrak{i}}{\mathcal{U}}^{\op},\subsm{\pi}{\mathcal{F}})$;
\item $\left((U ,\mathcal{F}) , {(\subsm{\mathfrak{i}}{\mathcal{U}}^{\op}, \subsm{\pi}{\mathcal{F}}})\right)$ is a cone over $\mathbf{G}$ isomorphic to $\left((U ,\subsm{\mathcal{L}}{\mathbf{G}}) , {(\subsm{\mathfrak{i}}{\mathcal{U}}^{\op}, \dindi{\pi}{\mathcal{L}}{\mbf{G}}})\right)$ in the category of cones over $\mathbf{G}$.
\end{enumerate} 
\begin{proof}

   In order to prove that $(1)\Leftrightarrow (2)$ it is enough to show that $(U ,\subsm{\mathcal{L}}{\mathbf{G}})$ is a glued-up $\mathbf{C}$-object along $\mathbf{G}$ through $(\subsm{\mathfrak{i}}{\mathcal{U}}^{\op}, \dindi{\pi}{\mathcal{L}}{\mbf{G}})$. 
 
   We begin by verifying that $\subsm{\mathcal{L}}{\mathbf{G}}$ is in $\Co{\mathbf{C}}$ and $\dindiv{\pi}{\mathcal{L}}{\mbf{G}}{i} $ is in $\Cn{\mathbf{C}}$. Let  $V\subsm{\subseteq}{\operatorname{op}} U$, $\dindsi{\mathcal{L}}{\mathbf{G}}{0}\!(V)$ is an abelian group since ${\dindsi{\Goi{\mbf{G}}{i}}{\mathbf{Sh}}{0}}\!\big(V\cap \subsm{U}{i}\big)$ is an abelian group for all $i\in \rm{I}$. Moreover, by definition of ${{\dindgiv{\pi}{\mathcal{L}}{\mathbf{G}}{i}{V}}}$, we deduce easily that ${{\dindgiv{\pi}{\mathcal{L}}{\mathbf{G}}{i}{V}}}$ is a morphism of abelian groups for each $i \in \rm{I}$.
   
 Next, we have ${\subsm{\mathcal{L}}{\mathbf{G}}}$ is a presheaf since $\subsm{\Goi{\mbf{G}}{i}}{\mathbf{Sh}}$ is a presheaf, for all $i\in \rm{I}$.    
When $\mathbf{C}= \mathbf{E}\Sh{\Ops{X}}{Ab}$, we prove that $\subsm{\mathcal{L}}{\mathbf{G}}$ is a sheaf.   
Let $\subsm{\{\subsm{V}{k} \}}{k\!\in\! \mathrm{K}}$ be an open covering of $V$.
Suppose that $s=\subsm{(\subsm{s}{i})}{i\!\in \!\rm{I}}\in \dindsi{\mathcal{L}}{\mathbf{G}}{0}\!(V)$ and $\subsm{s|}{\subsm{V}{k}}=0$, for all $k \in\mathrm{ K}$ where 
$\subsm{s|}{\subsm{V}{k}}=\subsm{\left(\subsm{\subsm{s}{i}|}{\subsm{V}{k}\cap \subsm{U}{i}}\right)}{i\!\in \!\rm{I}}.$ Since $\subsm{\Goi{\mbf{G}}{i}}{\mathbf{Sh}}$ is a sheaf for each $i\in \rm{I}$, from the uniqueness property of a sheaf, we get $\subsm{s}{i}=0$ for all $i\in \rm{I}$. Hence $s=\subsm{(\subsm{s}{i})}{i\!\in\! \rm{I}}=0$. 
This proves the uniqueness condition for $\subsm{\mathcal{L}}{\mathbf{G}}$. To prove the gluing condition, let $\subsm{t}{k}=(\subsm{{\subsm{t}{k_i}})}{i\!\in\! \rm{I}} \!\in\! \dindsi{\mathcal{L}}{\mathbf{G}}{0}\!(\subsm{V}{k})$ for $k \in \mathrm{K}$ be a family of sections such that
\begin{equation}\label{secn}
\subsm{\subsm{t}{k}|}{\subsm{V}{k}\cap \subsm{V}{k'}}=\subsm{\subsm{t}{k'}|}{\subsm{V}{k}\cap\subsm{V}{k'}}
\end{equation}
 for all $k,k'\in \mathrm{K}$. We want to prove that there exists a section $s=\subsm{(\subsm{s}{i})}{i\!\in\! \rm{I}}\!\in\! \dindsi{\mathcal{L}}{\mathbf{G}}{0}\!(V)$ such that $\subsm{s|}{\subsm{V}{k}}=\subsm{t}{k}$ for all $k\in \mathrm{K}$. 
 Let $i\in \mathrm{I}$. Equation (\ref{secn}) is equivalent to saying that, for $k,k'\in \mathrm{K}$
\begin{equation}\nonumber
\subsm{{\subsm{t}{k_{i}}}|}{\subsm{V}{k}\cap\subsm{V}{k'}\cap \subsm{U}{i}}=\subsm{{\subsm{t}{k'_{ i}}}|}{\subsm{V}{k}\cap \subsm{V}{k'}\cap \subsm{U}{i}}.
\end{equation}
Thus, since $\subsm{\Goi{\mbf{G}}{i}}{\mathbf{Sh}}$ is a sheaf and $\subsm{\{\subsm{V}{k}\cap \subsm{U}{i}\}}{i\!\in\! \rm{I}}$ is an open covering of $\subsm{U}{i}$, there exist $\subsm{s}{i} \!\in\! \dindsi{\Goi{\mbf{G}}{i}}{\mathbf{Sh}}{0}\!(V\cap \subsm{U}{i})$, such that
\begin{equation}\label{E6}
 \subsm{{\subsm{s}{i}}|}{\subsm{V}{k}\cap \subsm{U}{i}}=  {\subsm{t}{k_i}}
\end{equation}
for all $k\in \mathrm{K}$. We set $s:=\subsm{(\subsm{s}{i})}{i\!\in\! \rm{I}}$. We want to prove $s\in \dindsi{\mathcal{L}}{\mathbf{G}}{0}\!(V)$. That is, 
\begin{equation}\label{E.7}
\dindsi{\Gnij{\bf{G}}{\tauij{i}{j}}}{\mathbf{Sh}}{V}\!(\subsm{\subsm{s}{i}|}{V\cap \subsm{U}{i,j}})=\subsm{\subsm{s}{j}|}{V\cap \subsm{U}{i,j}}
\end{equation}
for all $i,j\in \rm{I}$. Let $i,j\in \rm{I}$. 
Since $\subsm{t}{k} \in \dindsi{\mathcal{L}}{\mathbf{G}}{0}\!(\subsm{V}{k})$, Equation (\ref{E6}) implies 
\begin{equation} \label{Em}\nonumber
\dindsi{\Gnij{\bf{G}}{\tauij{i}{j}}}{\mathbf{Sh}}{\subsm{V}{k}}\!(\subsm{\subsm{s}{i}|}{\subsm{V}{k}\cap \subsm{U}{i,j}})=(\subsm{\subsm{s}{j}|}{\subsm{V}{k}\cap \subsm{U}{i,j}}).
\end{equation}
By the functorial property of restriction maps, the previous equality can be rewritten as 
$$
\dindsi{\Gnij{\bf{G}}{\tauij{i}{j}}}{\mathbf{Sh}}{V} \!\subsm{(\subsm{\subsm{s}{i}|}{V\cap \subsm{U}{i,j}})|}{\subsm{V}{k}\cap \Uij{i}{j}}=\subsm{(\subsm{\subsm{s}{j}|}{V\cap \subsm{U}{i,j}})|}{\subsm{V}{k}\cap \subsm{U}{i,j}}
$$
for all $k\in \mathrm{K}$. Since $\subsm{\{\subsm{V}{k}\cap \subsm{U}{i,k}\}}{k\!\in \!\mathrm{K}}$ is an open covering of $V\cap \subsm{U}{i,j}$, by the uniqueness property of the sheaf $\subsm{\Goi{\mbf{G}}{j}}{\mathbf{Sh}}$, we obtain that Equation \ref{E.7} is satisfied. Hence $s\in \dindsi{\mathcal{L}}{\mathbf{G}}{0}\!(V)$.

 Now, for $\mathbf{C}$ an arbitrary element in the set $\{\mathbf{E}\Psh{\Ops{X}}{Ab},\mathbf{E}\Sh{\Ops{X}}{Ab}\}$, we want to prove that $\left((U ,\subsm{\mathcal{L}}{\mathbf{G}}) , {(\subsm{\mathfrak{i}}{\mathcal{U}}^{\op}, \dindi{\pi}{\mathcal{L}}{\mbf{G}}})\right)$ satisfies condition $(1)$-$(3)$ of the Theorem  \ref{preglued}. Let $V\subsm{\subseteq}{\operatorname{op}} U$, $s\in \dindsi{\mathcal{L}}{\mathbf{G}}{0}\!(V)$, $i,j ,k \in \rm{I}$ and $n \in \{ j,k\}$. We have  
 \begin{align*}
\dindsi{ \Gnij{\bf{G}}{\etaij{j}{i}}}{\mbf{Sh}}{V}\!({{ {\dindgiv{\pi}{\mathcal{L}}{\mathbf{G}}{j}{V}}}}\!\!(s))
&=\subsm{\subsm{s}{j}|}{V \cap \subsm{U}{j,i}} , 
  \text{by definition of ${{\dindgiv{\pi}{\mathcal{L}}{\mathbf{G}}{j}{V}}}$ and $\dindsi{ \Gnij{\bf{G}}{\etaij{j}{i}}}{\mbf{Sh}}{V}$}
     \\&= \dindgijv{\pi}{\mathcal{L}}{\mathbf{G}}{j}{i}{V}\! \!(s).
 \end{align*}

 Hence property $(1)$ of the Theorem \ref{preglued} is satisfied. Property $(2)$ is satisfied as follows:
 \begin{align*}
 \dindsi{\Gnij{\bf{G}}{\etaijk{n}{i}{j}{k}}}{\mathbf{Sh}}{V}\!\!\circ {\dindgijv{\pi}{\mathcal{L}}{\mathbf{G}}{i}{n}{V}\! \!}\!(s)
 &=\subsm{(\subsm{\subsm{s}{i}|}{V\cap \subsm{U}{i,n}})|}{V\cap \subsm{U}{i,j,k}}, \quad \text{by definition of ${{ \dindgijv{\pi}{\mathcal{L}}{\mathbf{G}}{i}{n}{V}\! \!}}$ and $\dindsi{\Gnij{\bf{G}}{\etaijk{n}{i}{j}{k}}}{\mathbf{Sh}}{V}$} \\&=\subsm{\subsm{s}{i}|}{V\cap \subsm{U}{i,j,k}}=\dindgijkv{\pi}{\mathcal{L}}{\mathbf{G}}{i}{j}{k}{V}\!\! \!(s). 
\end{align*}
Hence property $(2)$ of the Theorem \ref{preglued} is satisfied.
 Next, the following diagram 
  \begin{figure}[H] 
   \begin{center}
		{\tiny	\begin{tikzcd}[column sep=normal]
				&{\dindsi{\Goij{\mbf{G}}{j}{i}}{\mathbf{Sh}}{0}}\!\big(V\cap \subsm{U}{j,i}\big)  &&& {\dindsi{\Goi{\mbf{G}}{i}}{\mathbf{Sh}}{0}}\!(V\cap \subsm{U}{i})\arrow[swap]{lll}{\dindsi{{\Gnij{\bf{G}}{\tauij{j}{i}\circ \etaji{j}{i}}}}{\mbf{Sh}}{V}}\\
				 &{\dindsi{\Goi{\mbf{G}}{j}}{\mathbf{Sh}}{0}}\!\big(V\cap \subsm{U}{j}\big) \arrow[]{u}{\dindsi{\Gnij{\bf{G}}{\etaji{i}{j}}}{\mbf{Sh}}{V}}& & &\dindsi{\mathcal{L}}{\mathbf{G}}{0}\!(V) \arrow[]{lll}{ {\dindgiv{\pi}{\mathcal{L}}{\mathbf{G}}{j}{V}} }\arrow[swap]{u} {{ \dindgiv{\pi}{\mathcal{L}}{\mathbf{G}}{i}{V}}}
			\end{tikzcd}}\end{center} 	 \caption{}\label{topol140}	
		 \end{figure}
 \noindent commutes. Indeed, 
 \begin{align*}
     \dindsi{{\Gnij{\bf{G}}{\tauij{j}{i}\circc \etaji{j}{i}}}}{\mbf{Sh}}{V}\circc {{ {\dindgiv{\pi}{\mathcal{L}}{\mathbf{G}}{i}{V}}}}\!\!(s)
     &=\dindsi{\Gnij{\bf{G}}{\tauij{j}{i}}}{\mathbf{Sh}}{V}\!\!\left(\subsm{\subsm{s}{i}|}{V \cap \subsm{U}{i,j} }\right) \text{by definition of ${{ {\dindgiv{\pi}{\mathcal{L}}{\mathbf{G}}{i}{V}}}}$ and ${\dindsi{\Gnij{\bf{G}}{\etaij{i}{j}}}{\mbf{Sh}}{V}}$} \\&=\subsm{\subsm{s}{j}|}{V \cap \subsm{U}{i,j}}, \text{by definition of $\dindsi{\mathcal{L}}{\mathbf{G}}{0}\!(V)$}\\&={\dindsi{\Gnij{\bf{G}}{\etaij{j}{i}}}{\mbf{Sh}}{V}}\circc {{ {\dindgiv{\pi}{\mathcal{L}}{\mathbf{G}}{j}{V}}}}\!\!\!(s).
\end{align*} 
Since $U$ is a glued-up $\mathbf{(o)Top}^{\op}$-object along $\gcov{\mathcal{U}}$ through $\subsm{\mathfrak{i}}{\mathcal{U}}^{\op}$ by Definition \ref{gluingcover} (Lemma) and morphism are (open) continuous maps thus the first component of the limit is $U$.
 Suppose that $\left((U ,\mathcal{L}') , {(\subsm{\mathfrak{i}}{\mathcal{U}}^{\op}, \subsm{\pi}{\mathcal{L}'}})\right)$ is cone over $\mathbf{G}$ making the following diagram
\begin{figure}[H] 
   \begin{center}
		{\tiny	\begin{tikzcd}[column sep=normal]
				&{\dindsi{\Goij{\mbf{G}}{j}{i}}{\mathbf{Sh}}{0}}\!\big(V\cap \subsm{U}{j,i}\big)   &&& {\dindsi{\Goi{\mbf{G}}{i}}{\mathbf{Sh}}{0}}\!(V\cap \subsm{U}{i})\arrow[swap]{lll}{\dindsi{{\Gnij{\bf{G}}{\tauij{j}{i}\circ \etaji{j}{i}}}}{\mbf{Sh}}{V}}\\
				 &{\dindsi{\Goi{\mbf{G}}{j}}{\mathbf{Sh}}{0}}\!\big(V\cap \subsm{U}{j}\big) \arrow[]{u}{\dindsi{\Gnij{\bf{G}}{\etaji{i}{j}}}{\mbf{Sh}}{V}}&& & \Co{\mathcal{L}}'(V) \arrow[]{lll}{\dindpi{\pi}{\mathcal{L}}{\mkern+1.5mu \dindsi{}{j}{V}}}\arrow[swap]{u} {\dindpi{\pi}{\mathcal{L}}{\mkern+1.5mu \dindsi{}{i}{V}}}
			\end{tikzcd}}\end{center} 	 \caption{}\label{topol141}	
		 \end{figure} 
 \noindent commute, for all $i,j\in \rm{I}$. 
 We want to prove that there exists a unique map $\subsm{\mu}{V}: \Co{\mathcal{L}}'(V)\rightarrow \dindsi{\mathcal{L}}{\mathbf{G}}{0}\!(V)$ making the following diagram 
 \begin{figure}[H]
   \begin{center}
     {\tiny\begin{tikzcd}[column sep=normal]
{\dindsi{\Goij{\mbf{G}}{j}{i}}{\mathbf{Sh}}{0}}\!\big(V\cap \subsm{U}{j,i}\big)   & & &
{\dindsi{\Goi{\mbf{G}}{i}}{\mathbf{Sh}}{0}}\!(V\cap \subsm{U}{i}) \arrow[swap]{lll}{\dindsi{{\Gnij{\bf{G}}{\tauij{j}{i}\circ \etaji{j}{i}}}}{\mbf{Sh}}{V}}  \\
{\dindsi{\Goi{\mbf{G}}{j}}{\mathbf{Sh}}{0}}\!\big(V\cap \subsm{U}{j}\big) \arrow[]{u}{\dindsi{\Gnij{\bf{G}}{\etaji{i}{j}}}{\mbf{Sh}}{V}}     & & &
\dindsi{\mathcal{L}}{\mathbf{G}}{0}\!(V)\arrow[swap]{u}{ \dindgiv{\pi}{\mathcal{L}}{\mathbf{G}}{i}{V} }\arrow[]{lll}{{{ {\dindgiv{\pi}{\mathcal{L}}{\mathbf{G}}{j}{V}}}} } \\
&&&& \Co{\mathcal{L}}'(V) \arrow[swap,dashed]{ul}{\exists ! \subsm{\mu}{V}}\arrow[bend left=10]{ullll}{\dindpi{\mkern+1.5mu\pi}{\mathcal{L}}{\mkern+1.5mu\dindsi{}{j}{V}}} \arrow[swap,bend right=15]{uul}{\dindpi{\pi}{\mathcal{L}}{\mkern+1.5mu\dindsi{}{i}{V}}}
\end{tikzcd} }
   \end{center}\caption{}\label{topol142}
\end{figure}
 \noindent commute, for all $i,j\in \rm{I}$ and $V\subsm{\subseteq}{\operatorname{op}} U$. Let $i,j\in \rm{I}$ and $V\subsm{\subseteq}{\operatorname{op}} U$. When such a $\mu_V$ exists, given $s\in \Co{\mathcal{L}}'(V)$, the diagram gives 
\begin{equation}\label{g16}\nonumber
 {{ {\dindgiv{\pi}{\mathcal{L}}{\mathbf{G}}{i}{V}}}}\!\circc \subsm{\mu}{V}(s)= {{ {\dindgiv{\pi}{\mathcal{L}}{\mathbf{G}}{i}{V}}}}\!\!(\subsm{\mu}{V}(s))=\dindpi{\pi}{\mathcal{L}}{\mkern+1.5mu \dindsi{}{i}{V}}\!(s).
\end{equation}
 That is equivalent to having 
 \begin{equation}\nonumber
\subsm{\mu}{V}(s)=\subsm{(\dindpi{\pi}{\mathcal{L}}{\mkern+1.5mu \dindsi{}{i}{V}}\!(s))}{i\!\in\! \rm{I}} .
 \end{equation}
This shows that $\subsm{\mu}{V}$ is uniquely determined. 
 By the commutativity of the diagram in Figure \ref{topol142}, we obtain that $\subsm{(\dindpi{\pi}{\mathcal{L}}{\mkern+1.5mu \dindsi{}{i}{V}}\!(s))}{i\!\in\! \rm{I}}$ is an element of $\dindsi{\mathcal{L}}{\mathbf{G}}{0}\!(V)$. This proves that such a $\mu$ exists. 
 Finally, we prove that $\mu$ is a natural transformation. Let $W\subsm{\subseteq}{\operatorname{op}} V$. 
We have 
\begin{align*}
    \dindsi{\mathcal{L}}{\mathbf{G}}{1}\!(\iuv{W}{V}^{\op})\circ \mu_V(s)&= \dindsi{\mathcal{L}}{\mathbf{G}}{1}\!(\iuv{W}{V}^{\op})\circc \subsm{(\dindpi{\pi}{\mathcal{L}}{\mkern+1.5mu \dindsi{}{i}{V}}\!(s))}{i\!\in\! \rm{I}},\text{by definition of $\mu_V$}\\
    &= \subsm{\big(\dindsi{\Goi{\mbf{G}}{i}}{\mathbf{Sh}}{1}\!\!\big(\iuvop{W\cap \subsm{U}{i}}{V\cap \subsm{U}{i}}{\op}\big)\circc \dindpi{\pi}{\mathcal{L}}{\mkern+1.5mu \dindsi{}{i}{V}}\!(s) \big)}{i\!\in\! \rm{I}}, \text{by definition of $\dindsi{\mathcal{L}}{\mathbf{G}}{1}\!(\iuv{W}{V}^{\op})$}\\
        &= \subsm{\big(  \dindpi{\pi}{\mathcal{L}}{\mkern+1.5mu \dindsi{}{i}{W}}  \circc \Cn{\mathcal{L}}'(\iuvop{W}{V}{\op})(s) \big)}{i\!\in\! \rm{I}}, \text{since $\dindpi{\pi}{\mathcal{L}}{\mkern+1.5mu i}$ is a natural tranformation}
    \\&=  \mu_W\circ \Cn{\mathcal{L}}'(\iuvop{W}{V}{\op})(s).   \end{align*}
This completes the proof.

\end{proof}

\end{theorem}

\section{Gluing ringed spaces and schemes categorically}
For further insights on ringed topological spaces and scheme theory, we recommend referring to \cite{harts} and \cite{Liu}. 

\subsubsection{\texorpdfstring{$\rt$ (} \texorpdfstring{resp.} \texorpdfstring{$\lrt$} \texorpdfstring{} resp. \texorpdfstring{$\sch$)}\texorpdfstring{-}gluing data functors} 
We initiate the discussion by introducing the concept of an $\rt$ (resp. $\lrt$, resp. $\sch$)-gluing data functor. To enhance clarity and ease of notation throughout this section, we will use the symbol $\mathbf{C}$ to denote an element belonging to the set $\{\rt, \lrt, \sch\}$.

\begin{definition}\label{shcmct}
	A $\mathbf{C}$-\textbf{\textit{gluing data functor}} $\mathbf{G}$ is a functor from $\glI{I}$ to $\mathbf{C}$ such that for all $i,j, k \in\! \rm{I}$ and $n \in\! \{ j,k \}$ we have 
	\begin{enumerate}
		\item  $\Goij{\mbf{G}}{i}{j}=(\subsm{U}{i,j}, \subsm{\subsm{\Goi{\mbf{G}}{i}}{\mbf{Sh}}|}{\subsm{U}{i,j}});$
		\item $\Goijk{\mbf{G}}{i}{j}{k} =\subsm{\Goi{\mbf{G}}{i}|}{\Uij{i,j}{i,k}}$;
		\item $\Gnij{\bf{G}}{\etaij{i}{j}}=\subsm{\subsm{\underbar{$\mathfrak{i}$}}{\Goi{\mbf{G}}{i}}}{\subsm{U}{i,j},\!\subsm{\Goi{\mbf{G}}{i}}{\mathbf{Top}}}$;
		
		\item $\Gnij{\bf{G}}{\etaijk{n}{i}{j}{k}}=\subsm{\subsm{\underbar{$\mathfrak{i}$}}{\Goi{\mbf{G}}{i}}}{\Uij{i,j}{i,k},\!\subsm{U}{i,n}}$.       
		
	\end{enumerate} 
\end{definition}

The $\mathbf{C}$-gluing data functor induces naturally an $\topo^{\op}$-gluing data functor and a $\mathbf{E}\mathbf{Sh}$-gluing data functor as follows:
\begin{definition}[Lemma]\label{gtop} Let $\mathbf{G}$ be a $\mathbf{C}$-gluing data functor. We define the $\topo^{\op}$-\textbf{\textit{gluing data functor induced by $\mathbf{G}$}} denoted by $\subsm{\mathbf{G}}{\mathbf{Top}}$ as the functor such that $\dindsi{\mbf{G}}{\mbf{Top}}{0}$ sends $a\in \Co{\glI{I}}$ to $\subsm{\Goi{\mbf{G}}{a}}{\mathbf{Top}}$ and $\dindsi{\mathbf{G}}{\mbf{Top}}{1}$ sends $f$ in $\Cn{\glI{I}}$ to $\subsm{\Gnij{\bf{G}}{f}}{\mathbf{Top}}$. 
\end{definition}

\begin{definition} [Lemma]\label{gshef}
	Let $\mathbf{G}$ be a $\mathbf{C}$-gluing data functor and $\iota:=\dindi{\iota}{Q}{\subsm{\mbf{G}}{\mbf{Top}}}$. We define the {\bf $\mathbf{E}\Sh{\dindi{Q}{\mbf{G}}{\mbf{Top}}}{\mathbf{\mathbf{Ab}}}$-\textbf{\textit{gluing data functor along }} $\mathcal{U}:=\subsm{\{\subsm{\iota}{i}(\subsm{\Goi{\mbf{G}}{i}}{\mathbf{Top}})\}}{i\!\in\! \rm{I}}$ induced by $\mathbf{G}$} denoted $\mathbf{G}_{\mathbf{Sh}}$ as the functor such that
	\begin{enumerate}
		\item $\dindsi{\mathbf{G}}{\mathbf{Sh}}{0}\!(i):=\big(\subsm{\iota}{i}(\subsm{\Goi{\mbf{G}}{i}}{\mathbf{Top}}), \subsm{\iota}{i}\sta\subsm{\Goi{\mbf{G}}{i}}{\mathbf{Sh}} \big)$;
		\item $\dindsi{\mathbf{G}}{\mathbf{Sh}}{0}\!(i,j):=(\iotaij{i}{j}(\subsm{\Goij{\mbf{G}}{i}{j}}{\mbf{Top}}),\iotaij{i}{j}\sta\subsm{\Goij{\mbf{G}}{i}{j}}{\mbf{Sh}})$; 
		\item $ \subsm{{\dindsi{\mathbf{G}}{\mathbf{Sh}}{1}}\!(\tauij{i}{j})}{\mbf{Top}}:=\left( \iotaij{j}{i}\circ \subsm{\Gnij{\bf{G}}{\tauij{i}{j}}}{\mbf{Top}}^{\op}\circ \widetilde{\iotaij{i}{j}}^{-1}\right)^{\op}$   where $ \widetilde{\iotaij{i}{j}}$ is the map  $\iotaij{i}{j}$ corestricted to $\iotaij{i}{j}(\subsm{\Goij{\mbf{G}}{i}{j}}{\mbf{Top}})$ and
		\item $\subsm{{\dindsi{\mathbf{G}}{\mathbf{Sh}}{1}}\!(\tauij{i}{j})}{\mbf{Sh}}$ is a family of maps defined as $\dindsi{{\dindsi{\mathbf{G}}{\mathbf{Sh}}{1}}\!(\tauij{i}{j})}{\mbf{Sh}}{V}\!\!\!:=\dindsi{\Gnij{\bf{G}}{\tauij{i}{j}}}{\mathbf{Sh}}{\iotaij{j}{i}^{-1}(V)}$, for any $V\subsm{\subseteq}{\operatorname{op}} \;{\iotaij{j}{i}}(\subsm{\Goij{\mbf{G}}{j}{i}}{\mbf{Top}})$
	\end{enumerate}  
	\begin{proof} Since $\subsm{Q}{\gdt{G}}$ is a glued-up $\topo^{\op}$-object along $\subsm{\mathbf{G}}{\mathbf{Top}}$ through ${\iota}$, by Theorem \ref{gluingtop} we know that $\subsm{\{\subsm{\iota}{i}(\subsm{\Goi{\mbf{G}}{i}}{\mathbf{Top}})\}}{i\!\in\! \rm{I}}$ is an open covering of $\subsm{Q}{\gdt{G}}$. Moreover, 
		$$\iotaij{i}{j}(\subsm{\Goij{\mbf{G}}{i}{j}}{\mbf{Top}})=\iotaij{j}{i}(\subsm{\Goij{\mbf{G}}{j}{i}}{\mbf{Top}})=\subsm{\iota}{i}(\subsm{\Goi{\mbf{G}}{i}}{\mathbf{Top}})\cap\subsm{\iota}{j}(\subsm{\Goi{\mbf{G}}{j}}{\mathbf{Top}})$$ by Theorem \ref{gluingtop}. In addition, for any $V\subsm{\subseteq}{\operatorname{op}} \;\iotaij{j}{i}(\subsm{\Goij{\mbf{G}}{j}{i}}{\mbf{Top}})$ we have 
		\begin{align*}
			\Co{\dindsi{\mathbf{G}}{\mathbf{Sh}}{0}\!(j,i)}(V)={\iotaij{j}{i}}\sta  \subsm{\subsm{\Goij{\mbf{G}}{j}{i}}{\mathbf{Sh}}|}{{\dindsi{\Goij{\mbf{G}}{i}{j}}{\mbf{Top}}{0}}}\!(V)=\dindsi{\Goij{\mbf{G}}{j}{i}}{\mathbf{Sh}}{0}(\iotaij{j}{i}^{-1}(V))
		\end{align*}
		and 
		\begin{align*}
			\Co{\Big( \subsm{{\dindsi{\mathbf{G}}{\mathbf{Sh}}{1}}\!(\tauij{i}{j})}{\mbf{Top}}^{\op}\sta \dindsi{\mathbf{G}}{\mathbf{Sh}}{0}\!(i,j)\Big)} (V)
			&= \dindsi{\Goij{\mbf{G}}{i}{j}}{\mathbf{Sh}}{0}\left( \left( \iotaij{i}{j}^{-1} \circ \widetilde{\iotaij{i}{j}}\circ {\subsm{\Gnij{\bf{G}}{\tauij{i}{j}}}{\mbf{Top}}^{\op}}^{-1}\circ {\iotaij{j}{i}}^{-1}\right) (V)\right)\\
			&= \dindsi{\Goij{\mbf{G}}{i}{j}}{\mathbf{Sh}}{0}\left( \left( {\subsm{\Gnij{\bf{G}}{\tauij{i}{j}}}{\mbf{Top}}^{\op}}^{-1}\circ {\iotaij{j}{i}}^{-1}\right) (V)\right)\\
			&= \Co{\Big( \subsm{{\dindsi{\mathbf{G}}{\mathbf{Sh}}{1}}\!(\tauij{i}{j})}{\mbf{Top}}^{\op}\sta \subsm{\Goij{\bf{G}}{i}{j}}{\mathbf{Sh}}\Big)} ({\iotaij{i}{j}^{-1}(V)}).
		\end{align*}
		Thus, the constructions of the statements are well-defined and turn $\subsm{\mathbf{G}}{\mathbf{Sh}}$ into a $\mathbf{E}\Sh{\dindi{Q}{\mbf{G}}{\mbf{Top}}}{\mathbf{\mathbf{Ab}}}$-gluing data functor. 
	\end{proof}
\end{definition}
Below, we present the usual definition of a ringed topological space (resp. locally ringed topological space resp. scheme) gluing data.
\begin{definition}
	A collection $$\left(\mathrm{I},\subsm{( \subsm{\mathfrak{U}}{i})}{i\!\in\! \mathrm{I}}, \subsm{(\subsm{U}{i,j})}{(i,j)\!\in\! \mathrm{I}^2}, \subsm{(\subsm{\alpha}{i,j})}{(i,j)\!\in \!\mathrm{I}^2}\right)$$ 
	where
	\begin{enumerate}
		\item $\mathrm{I}$ is a set; 
		\item $\subsm{( \subsm{\mathfrak{U}}{i})}{i\!\in\! \mathrm{I}}$ is a family of ringed topological space (resp. locally ringed topological space resp. scheme);
		\item $\subsm{(\subsm{U}{i,j})}{(i,j)\!\in\! \mathrm{I}^2}$ is a family of open subsets of $\subsm{U}{i}$;
		\item $\subsm{(\subsm{\alpha}{i,j})}{(i,j)\!\in \!\mathrm{I}^2}$ is a family of isomorphism of ringed topological space (or locally ringed topological space or scheme) from $\subsm{\subsm{\mathfrak{U}}{i}|}{\subsm{U}{i,j}}$ to $\subsm{\subsm{\mathfrak{U}}{j}|}{\subsm{U}{j,i}}$.
	\end{enumerate}
	such that  
	\begin{enumerate}
		\item[a)] for each $i\in \rm{I}$, $\subsm{\alpha}{i,i}=\subsm{\operatorname{id}}{\subsm{\subsm{\mathfrak{U}}{i}|}{\subsm{U}{i,i}}}$;
		
		\item[b)] for each $i,j,k\in \rm{I}$, $\subsm{\alpha}{i,j}(\subsm{\subsm{\mathfrak{U}}{i}|}{\subsm{U}{i, j \wedge i,k }})=\subsm{\subsm{\mathfrak{U}}{j}|}{\subsm{U}{j, i \wedge j,k }}$ and
		\item[c)]  for each $i,j,k\in \rm{I}$, $\subsm{\alpha}{i,k}=\subsm{\alpha}{j,k}\circc \subsm{\alpha}{i,j}$ on $\subsm{U}{i, j \wedge i,k }$ 
	\end{enumerate}
	is called a \textbf{\textit{ringed topological space (resp. locally ringed topological space resp. scheme) gluing data}}.
	
\end{definition}

The following lemma describes the equivalence between $\mathbf{C}$-gluing data functor and the corresponding gluing data.  The proof of the following lemma follows directly from Lemma \ref{equi} and Lemma \ref{equishf}.

\begin{lemma}
	A $\rt$ (resp. $\lrt$, resp. $\sch$)-gluing data functor $\mathbf{G}$ along $\mathcal{U}$ induces  a ringed topological space (resp. locally ringed topological space, resp. scheme) gluing data $$\left(\mathrm{I}, \subsm{(\Goi{\mbf{G}}{i})}{i\!\in\! \mathrm{I}}, \subsm{(\subsm{\Goij{\mbf{G}}{i}{j}}{\mathbf{Top}})}{(i,j)\!\in\! \mathrm{I}^2}, \subsm{(\Gnij{\bf{G}}{\tauij{j}{i}})}{(i,j)\!\in\! \mathrm{I}^2}\right)$$ along $\mathcal{U}$. 
	Conversely, a ringed topological space (resp. locally ringed topological space resp. scheme) gluing data $\left(\mathrm{I},\subsm{( \subsm{\mathfrak{U}}{i})}{i\!\in\! \mathrm{I}}, \subsm{(\subsm{U}{i,j})}{(i,j)\!\in\! \mathrm{I}^2}, \subsm{(\subsm{\alpha}{i,j})}{(i,j)\!\in \!\mathrm{I}^2}\right)$ induces the $\rt$ (resp. $\lrt$, resp. $\sch$)-gluing data functor $\mathbf{G}$ defined by $\Goi{\mbf{G}}{i}=\subsm{\mathfrak{U}}{i}$, $\Goij{\mbf{G}}{i}{j}=\subsm{\Goi{\mbf{G}}{i}|}{\subsm{U}{i,j}}$ and $\Gnij{\bf{G}}{\tauij{j}{i}}=\subsm{\alpha}{i,j}$ for all $i,j\in \mathrm{I}$. 
	
\end{lemma}

\subsubsection{Characterization of glued-up \texorpdfstring{$\rt$ (} \texorpdfstring{resp.} \texorpdfstring{$\lrt$} \texorpdfstring{} resp. \texorpdfstring{$\sch$)}\texorpdfstring{-}objects}
We define and characterize $\mathbf{C}$-glued-up objects. To introduce the notation we will use, we begin by recalling the definition and some results about the stalk of a presheaf. The proofs of these results are well known, so we refer the reader to \cite[Chapter 2, p.35]{Liu}.
\begin{definition}[Lemma]\label{stalk}
	Let $\mathfrak{R}$ be a ringed topological space and $x \in \subsm{\mathfrak{R}}{\mathbf{Top}}$.
\end{definition}
\begin{enumerate}
	\item We define the stalk diagram at $x$ as follows:
	$$ \subsmt{\mathbf{S}}{\mathfrak{R}}{x}: \subsmt{\mathbf{Ops}}{\subsm{\mathfrak{R}}{\mathbf{Top}}}{x}^{\op} \rightarrow \mathbf{Ab} $$
	such that $\Cos{\subsmt{\mathbf{S}}{\mathfrak{R}}{x}}\mkern-1mu(U) := \dindsi{\mathfrak{R}}{\mathbf{Sh}}{0}\!(U)$ and $\Cns{\subsmt{\mathbf{S}}{\mathfrak{R}}{x}}\!(\iuvop{V}{U}{\op}):=\dindsi{\mathfrak{R}}{\mathbf{Sh}}{1}\!(\iuvop{V}{U}{\op})$, for any $U, V \in \Cos{\subsmt{\mathbf{Ops}}{\subsm{\mathfrak{R}}{\mathbf{Top}}}{x}}$ where $ V \subsm{\subseteq}{\operatorname{op}} U$.
	
	\item  We denote 
	\begin{itemize}  
		\item $\dindsi{\mathfrak{R}}{\mathbf{Sh}}{x}:=\left(\subsm{\coprod\nolimits}{U \!\! \in \!  \Cos{\subsmt{\mathbf{Ops}}{\subsm{\mathfrak{R}}{\mathbf{Top}}}{x}}} \! \Cos{\subsm{\mathfrak{R}}{\mathbf{Sh}}} \!(U)\right)/\Rel{\subsmt{\mathbf{S}}{\mathfrak{R}}{x}}$ where $\Rel{\subsmt{\mathbf{S}}{\mathfrak{R}}{x}}$ is the equivalence relation on $\subsm{\coprod\nolimits}{U\!\!\in\! {\Cos{\subsmt{\mathbf{Ops}}{\subsm{\mathfrak{R}}{\mathbf{Top}}}{x}}}}\dindsi{\mathfrak{R}}{\mathbf{Sh}}{0}\!(U)$ defined as follows: \\
		for any $(s ,U ),(t ,V )\in \subsm{\coprod\nolimits}{U\!\!\in\! {\Cos{\subsmt{\mathbf{Ops}}{\subsm{\mathfrak{R}}{\mathbf{Top}}}{x}}}}\dindsi{\mathfrak{R}}{\mathbf{Sh}}{0}(U)$  we have $(s ,U )\Rel{\subsmt{\mathbf{S}}{\mathfrak{R}}{x}} (t ,V )$ if and only if there exists $W\in \Cos{\subsmt{\mathbf{Ops}}{\subsm{\mathfrak{R}}{\mathbf{Top}}}{x}}$ with
		$W\subsm{\subseteq}{\operatorname{op}} U$ and  $W\subsm{\subseteq}{\operatorname{op}} V$ such that $\dindsi{\mathfrak{R}}{\mathbf{Sh}}{1}\!(\iuvop{W}{U}{\op})(s )= \dindsi{\mathfrak{R}}{\mathbf{Sh}}{1}\!(\iuvop{W}{V}{\op})(t)$.  We denote $[s ,U ]$ the class of $(s ,U )$ in the quotient $\dindsi{\mathfrak{R}}{\mathbf{Sh}}{x}$;
		\item $\subsm{\iota}{\mathbf{Sh},\! x}:=\subsm{(\subsms{\subsm{\iota}{\mathbf{Sh},\! x}}{U})}{U\!\!\in\! \Cos{\subsmt{\mathbf{Ops}}{\subsm{\mathfrak{R}}{\mathbf{Top}}}{x}}}$ is the family of maps  $\subsms{\subsm{\iota}{\mathbf{Sh},\! x}}{U}: \dindsi{\mathfrak{R}}{\mathbf{Sh}}{0}\!(U) \rightarrow \dindsi{\mathfrak{R}}{\mathbf{Sh}}{x}$ sending $s$ to $[s ,U ]$, for each $U\in\Cos{\subsmt{\mathbf{Ops}}{\subsm{\mathfrak{R}}{\mathbf{Top}}}{x}}$.
	\end{itemize} 
\end{enumerate}
We have the following:
\begin{itemize}
	\item 
	$\operatorname{lim} \subsmt{\mathbf{S}}{\mathfrak{R}}{x} \simeq (\dindsi{\mathfrak{R}}{\mathbf{Sh}}{x},\subsm{\iota}{\mathbf{Sh},x})$;
	\item  Given a morphism of locally ringed topological space $\Phi: \mathfrak{R}\rightarrow \mathfrak{S}$ and $\Psi: \mathfrak{S}\rightarrow \mathfrak{T}$ then there is a canonical morphism of local rings $\dindsi{\Phi}{\mathbf{Sh}}{x}: \dindi{\mathfrak{R}}{\mathbf{Sh}}{\subsm{\Phi}{\operatorname{Top}}(x)}\rightarrow \dindsi{\mathfrak{S}}{\mathbf{Sh}}{x}$ sending $[s, U]$ to $[\dindsi{\Phi}{\mathbf{Sh}}{U}\!(s),\subsm{\Phi}{\mbf{Top}}^{-1}(U)]$ where $U\!\!\in\! \Cos{\subsmt{\mathbf{Ops}}{\subsm{\mathfrak{R}}{\mathbf{Top}}}{x}}$ and $ s\in \dindsi{\mathfrak{R}}{\mathbf{Sh}}{0}\!(U)$. Moreover, we have $\dindsi{(\Psi\circc \Phi)}{\mathbf{Sh}}{x}=\dindsi{\Psi}{\mathbf{Sh}}{x} \circc \dindsi{\Phi}{\mathbf{Sh}}{x}$.
\end{itemize}
\begin{lemma}\label{stalky}
	Let $\mathbf{G}$ be a $\mathbf{C}$-gluing data functor, $\pi :=\subsm{\pi}{\subsm{\mathcal{L}}{\mathbf{G}}}$, $\mathcal{L}:=\subsm{\mathcal{L}}{\mathbf{G}}$, $i\in \rm{I}$ and $x\in \subsm{\Goi{\mbf{G}}{i}}{\mathbf{Top}}$. Then,
	\begin{enumerate}
		\item $\dindi{\mathcal{L}}{\mathbf{Sh}}{\subsms{\iota}{\mathbf{Top}}\!(x)}\simeq \dindsi{\Goi{\mbf{G}}{i}}{\mathbf{Sh}}{x}$ through the map ${\dindiv{\pi}{\mathbf{Sh}}{\mkern+1.7mu i}{\mkern+1.7mu x}}$;
		\item  the following diagram 
		\begin{figure}[H]
			\begin{center}
				{\tiny	\begin{tikzcd}[column sep=large]
						&\dindi{\mathcal{L}}{\mathbf{Sh}}{\subsms{\iota}{\mathbf{Top}}\!(x)} \arrow{rrr}{\dindsi{\mathcal{L}}{\mathbf{Sh}}{1}\!\!\subsm{(\iuvop{W}{V}{\op}\!)}{\! x}} \arrow[swap,"{\resizebox{0.5cm}{0.1cm}{$\sim$}}" labl1]{d}{\dindiv{\pi}{\mathbf{Sh}}{\mkern+1.7mu i}{\mkern+1.7mu x}} &&&\dindi{\mathcal{L}}{\mathbf{Sh}}{\subsms{\iota}{\mathbf{Top}}\!(x)}\arrow["{\resizebox{0.5cm}{0.1cm}{$\sim$}}" labl1]{d} {\: \;\;\dindiv{\pi}{\mathbf{Sh}}{\mkern+1.7mu i}{\mkern+1.7mu x} }\\ 
						&\dindsi{\Goi{\mbf{G}}{i}}{\mathbf{Sh}}{x} \arrow[swap]{rrr}{\dindsi{{\Goi{\mbf{G}}{i}}}{\mathbf{Sh}}{1}\!\!\subsm{\big(\iuvop{\subsm{\iota}{\mathbf{Top}}^{-1}\!(W)\cap \subsm{U}{i}}{\subsm{\iota}{\mathbf{Top}}^{-1}\!(V)\cap \subsm{U}{i}}{\op} \big)}{\! \! \! x}}& &&\dindsi{\Goi{\mbf{G}}{i}}{\mathbf{Sh}}{x}
				\end{tikzcd}}
				
			\end{center} \caption{}\label{topol148}  \end{figure} 
		\noindent commutes.
		
	\end{enumerate}
	\begin{proof} 
		\begin{enumerate}
			\item   
			To simplify notations we set $\iota:= \dindi{\iota}{Q}{\mbf{G}}$, $U_i:= \subsm{\Goi{\mbf{G}}{i}}{\mbf{Top}}$ and $y:=\subsms{\iota}{\mathbf{Top}}\!(x)$. 
			By Definition \ref{stalk} (Lemma), $\dindsi{\pi}{\mathbf{Sh}}{i}$  induces a stalk morphism  $\dindiv{\pi}{\mathbf{Sh}}{\mkern+1.7mu i}{\mkern+1.7mu x}:\dindsi{\mathcal{L}}{\mathbf{Sh}}{y}\rightarrow \dindsi{\Goi{\mbf{G}}{i}}{\mathbf{Sh}}{x}$ sending $[s,V]$ to $[\dindiv{\pi}{\mathbf{Sh}}{\mkern+1.7mu i}{\mkern+1.7mu V}\!\!(s), \subsm{\iota}{\mathbf{Top}}^{-1}(V)]$ for any $V\in \Co{\subsmt{\mathbf{Ops}}{\subsm{\iota}{\mathbf{Top}}(\subsm{U}{i}) }{y}}$ and $s\in \dindsi{\mathcal{L}}{\mathbf{Sh}}{0}(V)$. We want to prove that $\dindiv{\pi}{\mathbf{Sh}}{\mkern+1.7mu i}{x}$ is an isomorphism. We first prove that $\dindiv{\pi}{\mathbf{Sh}}{\mkern+1.7mu i}{\mkern+1.7mu x}$ is one-to-one. Let $[s,V]\in \dindsi{\mathcal{L}}{\mathbf{Sh}}{y}$ such that $\dindiv{\pi}{\mathbf{Sh}}{\mkern+1.7mu i}{\mkern+1.7mu x}\!([s,V])=[0,\subsm{\iota}{\mathbf{Top}}^{-1}(V)]$. We have $[\dindiv{\pi}{\mathbf{Sh}}{\mkern+1.7mu i}{\mkern+1.7mu V}\!\!(s), \subsm{\iota}{\mathbf{Top}}^{-1}(V)]=[0,\subsm{\iota}{\mathbf{Top}}^{-1}(V)]$. It implies that $(\dindiv{\pi}{\mathbf{Sh}}{\mkern+1.7mu i}{\mkern+1.7mu V}\! (s), \subsm{\iota}{\mathbf{Top}}^{-1}(V)) \Rel{\subsmt{\mathbf{S}}{\mathfrak{R}}{x}}  (0,\subsm{\iota}{\mathbf{Top}}^{-1}(V))$. Then, by definition of $\Rel{\subsmt{\mathbf{S}}{\mathfrak{R}}{x}} $ there exists $W\in \Co{\subsmt{\mathbf{Ops}}{\subsm{U}{i}}{x}}$ such that $ W\subsm{\subseteq}{\operatorname{op}} \subsm{\iota}{\mathbf{Top}}^{-1}(V)$ and $\dindsi{\Goi{\mbf{G}}{i}}{\mathbf{Sh}}{1}\!\!(\iuvop{W}{\:\subsm{\iota}{\mathbf{Top}}^{-1}(V)}{\op})(\dindiv{\pi}{\mathbf{Sh}}{\mkern+1.7mu i}{\mkern+1.7mu V}\! (s))=0$. Since $\dindi{\pi}{\mathbf{Sh}}{\mkern+1.7mu i}$ is a natural transformation, we have 
			\begin{align*}
				\dindsi{\Goi{\mbf{G}}{i}}{\mathbf{Sh}}{1}\!(\iuvop{W}{\:\subsm{\iota}{\mathbf{Top}}^{-1}(V)}{\op})(\dindiv{\pi}{\mathbf{Sh}}{\mkern+1.7mu i}{\mkern+1.7mu V}\! (s))&=\dindiv{\pi}{\mathbf{Sh}}{\mkern+1.7mu i}{\mkern+1.7mu W}\!\!(\dindi{\mathcal{L}}{\mbf{Sh}}{1}\!(\iuvop{\subsm{\iota}{\mbf{Top}}(W)}{\:V}{\op})(s))\\&=\dindiv{\pi}{\mathbf{Sh}}{\mkern+1.7mu i}{\mkern+1.7mu W}\!\!(\subsm{(\subsm{\subsm{s}{i}|}{\subsm{\iota}{\mbf{Top}}(W)\cap \subsm{U}{i}})}{i\!\in\! \rm{I}})=\subsm{\subsm{s}{i}|}{\subsm{\iota}{\mbf{Top}}(W)\cap \subsm{U}{i}}=0.
			\end{align*}   
			By definition of $\Cns{\subsm{\mathcal{L}}{\mbf{Sh}}}$, $\subsm{s|}{\subsm{\iota}{\mbf{Top}}(W)}=\Cns{\subsm{\mathcal{L}}{\mbf{Sh}}}\!(\iuvop{\subsm{\iota}{\mbf{Top}}(W)}{\:U}{\op})(s)=\subsm{(\subsm{\subsm{s}{j}|}{\subsm{\iota}{\mbf{Top}}(W)\cap \subsm{U}{j}})}{j\!\in \! \rm{I}}\in \dindsi{\mathcal{L}}{\mbf{Sh}}{0}\!(\subsm{\iota}{\mathbf{Top}}(W))$. 
			Then, since $\subsm{\iota}{\mbf{Top}}(W) \subsm{\subseteq}{\operatorname{op}} \subsm{U}{i}$, $\subsm{{\dindsi{\mathbf{G}}{\mathbf{Sh}}{1}}\!(\tauij{i}{j})}{\mbf{Sh}}\!(\subsm{\subsm{s}{i}|}{\subsm{\iota}{\mbf{Top}}(W)})=\subsm{\subsm{s}{j}|}{\subsm{\iota}{\mbf{Top}}(W)\cap \subsm{U}{j}}=0$, for all $j \in \rm{I}$. That is, $\subsm{s|}{\subsm{\iota}{\mbf{Top}}(W)}=0$. Proving that $[s, V] = [0,V]$.
			
			Next, we prove that  $\dindiv{\pi}{\mathbf{Sh}}{\mkern+1.7mu i}{\mkern+1.7mu x}$ is onto. Let $[t,V]\in \dindsi{\Goi{\mbf{G}}{i}}{\mathbf{Sh}}{x}$. We have $V\!\!=\subsm{\iota}{\mathbf{Top}}^{-1}\!\big( \subsm{\iota}{\mathbf{Top}}(V)\!\big)$ because $\subsm{\iota}{\mathbf{Top}}\in \Cn{\topo}$. We set $s:=\subsm{(\subsm{t|}{V\cap \subsm{U}{j}})}{j\!\in\! \rm{I}}$. It is clear that $s\in \dindsi{\mathcal{L}}{\mathbf{Sh}}{0}\!(\subsm{\iota}{\mathbf{Top}}(V))$   
			and $\dindiv{\pi}{\mathbf{Sh}}{\mkern+1.7mu i}{\mkern+1.7mu x}\!([s,\subsm{\iota}{\mathbf{Top}}(V)])=[t,V]$.
			
			\item  It follows easily from the fact that $\dindsi{\pi}{\mathbf{Sh}}{i}$ is a natural transformation. 
		\end{enumerate} 
	\end{proof}
\end{lemma}

The following definition describes a representative for the limit over a $\mathbf{C}$-gluing data functor.
\begin{definition}\label{glusch}
	Let $\mathbf{G}$  be a $\mathbf{C}$-gluing data functor and $\mathcal{U}:=\subsm{\{\subsm{\iota}{\subsm{Q}{\subsm{\mathbf{G}}{\mathbf{Top}}}}(\subsm{\Goi{\mbf{G}}{i}}{\mathbf{Top}})\}}{i\!\in\! \mathrm{I}}$. 
	We define \textbf{\textit{the standard representative of the limit of $\mathbf{G}$}} as the pair $$(( \subsm{Q}{\gdt{G}} , \subsm{\mathcal{L}}{\gds{G}}) , {(\subsm{\mathfrak{i}}{\mathcal{U}}^{\op},    \subsm{\pi}{\subsm{\mathcal{L}}{\gds{G}}}})).$$
	
\end{definition}

The following theorem describes a glued-up object up to isomorphism.
\begin{theorem}\label{gluedsch}
	Let $\mathbf{G}$ be a $\mathbf{C}$-gluing data functor, $\mathfrak{R}$ be a ringed topological space, $\subsm{\Pi}{\mathfrak{R}}$ be a pair $(\dindi{\Pi}{\mathfrak{R}}{\mathbf{Top}}, \dindi{\Pi}{\mathfrak{R}}{\mathbf{Sh}})$ where $\dindi{\Pi}{\mathfrak{R}}{\mathbf{Top}}$ is a family $\subsm{(\dindiv{\Pi}{\mathfrak{R}}{\mbf{Top}}{a})}{a\!\in\!\Co{\glI{I}}}$ with $\dindiv{\Pi}{\mathfrak{R}}{\mbf{Top}}{a}\!\! : \subsm{\mathfrak{R}}{\mathbf{Top}}\rightarrow \subsm{\Goi{\mbf{G}}{a}}{\mathbf{Top}}$ and $\dindi{\Pi}{\mathfrak{R}}{\mathbf{Sh}}$ is a family  $\subsm{(\dindiv{\Pi}{\mathfrak{R}}{\mbf{Sh}}{a})}{a\!\in\!\Co{\glI{I}}}$ with $\dindiv{\Pi}{\mathfrak{R}}{\mbf{Sh}}{a}\!\! : \subsm{\mathfrak{R}}{\mathbf{Sh}}\rightarrow \subsm{\Goi{\mbf{G}}{a}}{\mathbf{Sh}}$. 
	The following assertions are equivalent:
	\begin{enumerate} 
		\item $\mathfrak{R}$ is a glued-up $\mathbf{C}$-object along $\mathbf{G}$ through $\subsm{\Pi}{\mathfrak{R}}$;
		\item $\left(\mathfrak{R} , \subsm{\Pi}{\mathfrak{R}}\right)$ is a cone over $\mathbf{G}$ isomorphic to $((\subsm{Q}{\gdt{G}} ,\subsm{\mathcal{L}}{\gds{G}}) , {(\subsm{\iota}{\mathcal{U}}^{\op},  \subsm{\pi}{\subsm{\mathcal{L}}{\gds{G}}}}))$ in the category of cones over $\mathbf{G}$. 
		
	\end{enumerate} 
	\begin{proof} 
		
		To simplify notations let $\iota:= \dindi{\iota}{Q}{\mbf{G}}$, $\mathcal{L}:=\subsm{\mathcal{L}}{\mbf{G}}$ and $\pi :=\subsm{\pi}{\subsm{\mathcal{L}}{\mathbf{G}}}$. We first prove that $\subsm{\mathcal{L}}{\mbf{Sh}}$ is in $\Co{\mathbf{C}}$ and $\subsm{\pi}{\mbf{Sh}}$ is in $\Cn{\mathbf{C}}$. Let $V\in \dindi{\mbf{Ops}}{X}{0}$, $\dindsi{\mathcal{L}}{\mbf{Sh}}{0}\!(V)$ is a ring since $\dindsi{\mathbf{G}}{\mathbf{Sh}}{0}\! (V\cap \subsm{\iota}{\mathbf{Top}}(\subsm{\Goi{\mbf{G}}{i}}{\mathbf{Top}}))$ is a ring for each $i\in \rm{I}$. Furthermore, we can define $\dindiv{\pi}{\mathbf{Sh}}{\mkern+1.7mu i}{\mkern+1.7mu V} \!\!: \dindsi{\mathcal{L}}{\mbf{Sh}}{0}\!(V) \rightarrow \dindsi{\mathbf{G}}{\mathbf{Sh}}{0}\! (V\cap \subsm{\iota}{\mathbf{Top}}(\subsm{\Goi{\mbf{G}}{i}}{\mathbf{Top}}))$ to be a map sending $\subsm{(\subsm{s}{k})}{k\!\in \rm{I}}$ to $\subsm{s}{i}$. We deduce easily that $\dindi{\pi}{\mathbf{Sh}}{\mkern+1.7mu i}$ is a ring morphism 
		for each $i\in \rm{I}$. When $\mathbf{C}=\lrt $ or $\sch$, $\dindi{\mathcal{L}}{\mathbf{Sh}}{\subsms{\iota}{\mathbf{Top}}\!(x)}$ is a local ring and $\dindiv{\pi}{\mathbf{Sh}}{\mkern+1.7mu i}{\mkern+1.7mu x}$ is a local ring morphism by Lemma \ref{stalky}, for all $x\in \subsm{\Goi{\mbf{G}}{i}}{\mathbf{Top}}$.
		In order to prove that $(1)\Leftrightarrow (2)$ it is enough to show that  $(\subsm{Q}{\gdt{G}} ,\subsm{\mathcal{L}}{\mbf{Sh}})$ is a glued-up $\mathbf{C}$-object along $\mathbf{G}$ through $(\subsm{\mathfrak{i}}{\mathcal{U}}^{\op}, \subsm{\pi}{\mbf{Sh}})$. By Definition \ref{gtop} (Lemma), $\subsm{\mathbf{G}}{\mathbf{Top}}$ is an $\topo^{\op}$-gluing data functor induced by $\mathbf{G}$. Thus, by Theorem \ref{gluingtop} (2), we know that  $\subsm{Q}{\gdt{G}}$ is a glued-up $\topo^{\op}$-object along $\subsm{\mathbf{G}}{\mathbf{Top}}$ through $\subsm{\mathfrak{i}}{\mathcal{U}}^{\op}$. On the other hand, by Definition \ref{gshef}, $\subsm{\mathbf{G}}{\mathbf{Sh}}$ is a $\mathbf{E}\Sh{\dindi{Q}{\mbf{G}}{\mbf{Top}}}{\mathbf{\mathbf{Ab}}}$-gluing data functor induced by $\mathbf{G}$ along $\mathcal{U}$. Hence by Theorem \ref{rivegluedshf} $(2)$, we have that $\subsm{\mathcal{L}}{\mbf{Sh}}$ is a glued-up $\mathbf{E}\Sh{\dindi{Q}{\mbf{G}}{\mbf{Top}}}{\mathbf{\mathbf{Ab}}}$-object along $\subsm{\mathbf{G}}{\mathbf{Sh}}$ through $ \subsm{\pi}{\mbf{Sh}}$. Therefore, $(\subsm{Q}{\gdt{G}} ,\subsm{\mathcal{L}}{\mbf{Sh}})$ is a glued-up $\mathbf{C}$-object along $\mathbf{G}$ through $(\subsm{\mathfrak{i}}{\mathcal{U}}^{\op},  \subsm{\pi}{\mbf{Sh}})$.
		
	\end{proof}
	
\end{theorem}

\end{document}